\documentclass{article}%
\usepackage{amsmath}
\usepackage{amsfonts}
\usepackage{amssymb}
\usepackage{graphicx}
\usepackage[top=2.5cm, bottom=2.5cm, left=2cm, right=2cm, headheight=0.5cm,
footskip=60pt]{geometry}%
\providecommand{\U}[1]{\protect\rule{.1in}{.1in}}
\newtheorem{theorem}{Theorem}

\newtheorem{corollary}{Corollary}

\newtheorem{lemma}{Lemma}
\newtheorem{remarka}{Remark}

\newtheorem{proposition}{Proposition}
\newtheorem{remark}{Remark}

\allowdisplaybreaks
\numberwithin{lemma}{subsection}
\numberwithin{proposition}{subsection}
\numberwithin{equation}{section}
\numberwithin{remark}{subsection}
\numberwithin{definition}{subsection}
\numberwithin{theorem}{section}
\numberwithin{corollary}{subsection}

\begin{document}

\title{Vanishing capillarity limit of the Navier-Stokes-Korteweg system in one
dimension with degenerate viscosity coefficient and discontinuous initial density}
\author{Cosmin Burtea \thanks{Universit\'{e} de Paris and Sorbonne Universit\'{e},
CNRS, IMJ-PRG, F-75006 Paris, France. }, Boris Haspot \thanks{Universit\'{e}
Paris Dauphine, PSL Research University, Ceremade, Umr Cnrs 7534, Place du
Mar\'{ }echal De Lattre De Tassigny 75775 Paris cedex 16 (France),
haspot@ceremade.dauphine.fr }}
\date{}
\maketitle

\begin{abstract}
In the first main result of this paper we prove that one can approximate
discontinious solutions of the $1d$ Navier Stokes system with solutions of the
$1d$ Navier-Stokes-Korteweg system as the capilarity parameter tends to $0$.
Moreover, we allow the viscosity coefficients $\mu=\mu\left(  \rho\right)  $
to degenerate near vaccum. In order to obtain this result, we propose two main
technical novelties. First of all, we provide an upper bound for the density
verifing NSK that does not degenerate when the capillarity coefficient tends
to $0$. Second of all, we are able to show that the positive part of the
effective velocity is bounded uniformly w.r.t. the capillary coefficient. This
turns out to be crucial in providing a lower bound for the density.

The second main result states the existene of unique finite-energy global
strong solutions for the $1d$ Navier-Stokes system assuming only that
$\rho_{0},1/\rho_{0}\in L^{\infty}$. This last result finds itself a natural
application in the context of the mathematical modeling of multiphase flows.

\end{abstract}

\section{Introduction}

\subsection{Presentation of the models}

One classical model used to study mixtures of two or more compressible fluids
with different densities is the Navier-Stokes system (referred as NS in the
sequel):%
\begin{equation}
\left\{
\begin{array}
[c]{l}%
\partial_{t}\rho+\partial_{x}\left(  \rho u\right)  =0,\\
\partial_{t}\left(  \rho u\right)  +\partial_{x}\left(  \rho u^{2}\right)
-\partial_{x}\left(  \mu\left(  \rho\right)  \partial_{x}u\right)
+\partial_{x}(a\rho^{\gamma})=0,\\
\left(  \rho,u\right)  _{|t=0}=\left(  \rho_{0},u_{0}\right)  ,
\end{array}
\right.  \label{Navier_Stokes_1d}%
\end{equation}

The mixture is supposed to be itself a compressible barotropic fluid and, in
order to simplify matters, we suppose that it has constant temperature.
Practical examples include gas bubbles in water flowing in pipes, fluids
containing a melted substance, polluted air/water ..., see \cite{Ishii} for
more relevant examples. The unknowns are the velocity of the fluid
$u=u(t,x)\in\mathbb{R}$ and the density of the fluid $\rho=\rho(t,x)\in
\mathbb{R}^{+}$. We consider here that the pressure $P(\rho)=a\rho^{\gamma}$
is defined throughout a $\gamma$ state-law with $\gamma>1$ and $a>0$. In this
paper, we consider that the viscosity coefficient $\mu=$ $\mu(\rho)\geq0$ is
density-dependent, with $\mu\left(  0\right)  =0$.

The fact that the dissipative operator depends on the density itself is
motivated by physical considerations. Keeping in mind the fact that $\left(
\text{\ref{Navier_Stokes_1d}}\right)  $ governs flows of mixtures, it is
reasonable to suppose that each of the components has its own viscosity. We
also mention that system $\left(  \text{\ref{Navier_Stokes_1d}}\right)  $ can
be obtained via the Boltzmann equation through the Chapman-Enskog expansion to
the second order (see \cite{CC70}) which ensures that the viscosity
coefficient is then a function of the temperature. If we consider the case of
isentropic fluids, this dependence is expressed by a dependence on the density
function (we refer in particular to \cite{HoffSerre}). Let us also mention
that the case $\mu(\rho)=\rho$ is related to the so called viscous shallow
water system. This system with friction has been derived by Gerbeau and
Perthame in \cite{GP} from the Navier-Stokes system with a free moving
boundary in the shallow water regime at the first order. This derivation
relies on the hydrostatic approximation where the authors follow the role of
viscosity and friction on the bottom. It transpires that the mathematical
analysis of this situation is more involved than the constant viscosity case
since the operator $-\partial_{x}\left(  \mu\left(  \rho\right)  \partial
_{x}u\right)  $ is not strongly elliptic.

The system $\left(  \text{\ref{Navier_Stokes_1d}}\right)  $ is a \textit{Sharp
Interface Model}: assuming that at initial time the constituents of the
mixture are separated by a jump in the initial density $\rho_{0}$, then this
discontinuity persists for all later times, and it is advected by the fluid's
velocity see the early work of D. Hoff and J. Smoller \cite{HS85}. This allows
one to distinguish at any time the zones occupied by the "pure" constituents.

This picture changes if one takes in consideration capillary forces. In this
case the two fluids are not anymore separated by a sharp interface but rather
by a thin layer where the density, although passing continuously from one
fluid to another, it experiences large variations. This type of models are
called \textit{Diffusive Interface Models}. One classical model used to
describe this situation is the Navier-Stokes-Korteweg (NSK in the sequel)
system which, in its $1d$ version reads
\begin{equation}
\left\{
\begin{array}
[c]{l}%
\partial_{t}\rho+\partial_{x}\left(  \rho u\right)  =0,\\
\partial_{t}\left(  \rho u\right)  +\partial_{x}\left(  \rho u^{2}\right)
-\partial_{x}\left(  \mu\left(  \rho\right)  \partial_{x}u\right)
+\partial_{x}(a\rho^{\gamma})=c\partial_{x}K,\\
\left(  \rho,u\right)  _{|t=0}=\left(  \rho_{0},u_{0}\right)  .
\end{array}
\right.  \label{NSK_intro}%
\end{equation}
In the above system $K$ stands for the Korteweg capillarity tensor and its
general formula is
\begin{equation}
K=\rho\kappa\left(  \rho\right)  \partial_{xx}\rho+\frac{1}{2}\left(
\rho\kappa^{\prime}\left(  \rho\right)  -\kappa\left(  \rho\right)  \right)
\left(  \partial_{x}\rho\right)  ^{2}, \label{cappilar_intro}%
\end{equation}
for some positive function $\kappa=\kappa\left(  \rho\right)  $. Moreover,
$c>0$ is a positive constant measuring the "strength of the capillarity" or in
an other way "the thickness of the diffusive interface". Indeed the size of
$c$ is directly linked with the thickness of the transition zone separating
the pure phases. For the underlying physical principles behind the governing
equations, see the pioneering paper by J.-E. Dunn and J. Serrin in \cite{fDS}
but also \cite{VW,fK,fDS,fTN,fA,fC}.

Let us observe that formally, when $c\rightarrow0$, we recover the NS system
$\left(  \text{\ref{Navier_Stokes_1d}}\right)  $. Justifying mathematically
such a limit in a physically interesting setting, namely, where interfaces
could be tracked down, raises the following questions:

\begin{enumerate}
\item[Q1] Can we construct solutions for the system $\left(
\text{\ref{Navier_Stokes_1d}}\right)  $ with initial data $\rho_{0}$ allowing
discontinuities? How "wildly" can these discontinuities behave?

\item[Q2] Can we recover solutions of $\left(  \text{\ref{Navier_Stokes_1d}%
}\right)  $ with discontinuous initial data $\rho_{0}$ as the limit of
solutions of $\left(  \text{\ref{NSK_intro}}\right)  $?
\end{enumerate}

In this paper, we give positive answers to the above questions. Of course, we
cannot solve the problem for general $\mu\left(  \rho\right)  $ and
$\kappa\left(  \rho\right)  $, additional hypothesis which ensure some extra
algebraic structure to the underlying PDEs systems are needed. One typical
example we can treat is
\begin{equation}
\mu\left(  \rho\right)  =\rho^{\alpha}\text{ with }\alpha\in(0,\dfrac{1}%
{2})\text{ and }\kappa\left(  \rho\right)  =\frac{\mu^{2}\left(  \rho\right)
}{\rho^{3}}. \label{visc_capilarity_intro}%
\end{equation}

In order to obtain such results, we propose two new technical features in the
context of the NSK system. Before, explaining these features, let us mention
that, it is by now well-established, that the crucial point for obtaining
global existence results is to be able to obtain apriori estimates which
ensure that the density $\rho$ stays bounded and bounded away from vacuum
(i.e. $\rho>0$) all along the time. Indeed, explosion criteria shows that this
is the possible cause for the breakdown of strong solutions, see for instance
Theorem $1.1.$ from \cite{Constantin} and Theorem $3.3.$ from
\cite{BurHas2020,BurHas2021}. Even, for less regular solutions in the spirit
of Hoff-Serre \cite{Hof87,Ser86a,Ser86b}, the former qualitative information
for the densities ensure that one can carry out estimation in the Hoff-class
of regularity (referred to as intermediate-regularity solutions).

Regarding the second question Q$2$: in the case of the NSK system, first of
all, we are able to use a maximum principle in order \textit{to obtain an
}$L^{\infty}$\textit{-bound on the density that does not degenerate in the
vanishing capillarity coefficient limit }$c\rightarrow0$. Secondly, we show
that in order to obtain\textit{ a lower bound for the density, it is
sufficient to assume only a one sided inequality on the so-called effective
velocity}.\textit{ }At the level of the density, loosely speaking, this
amounts to ask a one sided inequality:%
\begin{equation}
\partial_{x}\varphi(\rho_{0})\leq M_{0}, \label{ineg_1}%
\end{equation}
for some $M_{0}\in\mathbb{R}$ and where $\varphi\left(  \rho\right)  $ is a
primitive of $\frac{\mu\left(  \rho\right)  }{\rho^{2}}$. Obviously, the
expression in $\left(  \text{\ref{ineg_1}}\right)  $ has a meaning in the
sense of measures. As these estimates are uniform w.r.t. capillarity parameter
$c,$ \textit{we are able to recover in the vanishing capillarity limit
}$c\rightarrow0$\textit{, solutions for the NS system where the initial data
needs only to verify }$\left(  \text{\ref{ineg_1}}\right)  $. We can observe
that functions $\varphi(\rho_{0})\in BV_{loc}$ with discontinuities with
negative jump enter in the scope of $\left(  \text{\ref{ineg_1}}\right)  $. As
far as the authors are aware, this is the first result in this direction, the
only somehow related results known are obtained for the zero
viscosity-capillarity limit by F. Charve and the second author in \cite{CH}
and by P. Germain and P.G. LeFloch in \cite{Germain}. The relative weak
assumptions that are required on the initial density are in sharp contrast
with what is usually assumed in the literature where the so-called
Bresch-Desjardin entropy \cite{BD} is used in order to ensure that the density
stays away from vacuum (see \cite{MV}). However the assumption that the
BD-entropy is finite at initial time allows only to consider continuous
initial densities.

The method of proof reveals that in the case of the NS system $\left(
\text{\ref{Navier_Stokes_1d}}\right)  $, the situation is much more better: we
can ensure that the density stays away from vacuum using only the fact that
this is the case at the initial time $t=0$. More precisely, under hypothesis
$\left(  \text{\ref{visc_capilarity_intro}}\right)  $ for $\mu\left(
\rho\right)  $ we are able to construct global strong solutions assuming only
that $\rho_{0}$ and $\frac{1}{\rho_{0}}\in L^{\infty}$. This gives an
affirmative answer to Q$1$ at least for mildly degenerate density-dependent
viscosities. Let us also point out that this result is also relevant with
respect to recent results obtained by Bresch and Hillairet in
\cite{BreschHillairetENS} in the context of multiphase modeling. Indeed the
authors show that solutions for a Baer-Nunziato type-system can be obtained as
week limits of highly-oscilating solutions of the NS system $\left(
\text{\ref{Navier_Stokes_1d}}\right)  $. The only norm that one can hope to
keep uniformly bounded in the presence of high-oscilations is the $L^{\infty}%
$-norm. Our existence result stated in Theorem \ref{theo3} thus enlarges the
class of viscosity coefficients for which Theorem $1$ from
\cite{BreschHillairetENS} holds true. In particular, the aforementioned result
holds true for viscosities that degenerate close to vacuum.

\subsection{A short review of known results}

Before, giving the formal statement of our main results we propose to the
reader a short review of the results concerning the existence of global strong
solutions for the NS and NSK systems. We mention here that we use the term
\textit{weak solution }to designate a distributional solution for which
uniqueness is not know to hold true. The term \textit{strong solution} will be
used whenever we may ensure the uniqueness property even thought the
derivatives appearing in the equations may not have a meaning a.e..

\subsubsection{The Navier-Stokes system}

We start mentioning some known results for the system $\left(
\text{\ref{Navier_Stokes_1d}}\right)  $.

\paragraph{The constant viscosity case}

The study of the well-posedness of the Cauchy problem for the compressible
Navier-Stokes equations with \textit{constant viscosity} \textit{coefficients}
can be tracked back to the pioneering works of the Russian school of PDEs, see
Ya. I. Kanel in \cite{Ka} and Kazhikhov and Shelukhin \cite{KS77} where they
considered the case of regular initial data (in particular, the densities are
always continuous). In the $80^{\prime}s$ D. Hoff and J. Smoller started a
program investigating the well-posedness issues for parabolic equations with
rough initial data. Up to the knowledge of the authors, they were the first to
observe that the discontinuities of the density are advected by the flow and
they persist all along the time owing to the regularization of the so called
effective viscous flux \cite{HS85}. The first results tackling the question of
rough initial densities for the NS system were obtained by D. Hoff
\cite{Hof87} and D. Serre \cite{Ser86a,Ser86b}. Qualitative properties for
weak-solutions like the non-formation of vacuum states were adressed in
\cite{HS01}. We mention also the result of D. Hoff and D. Serre
\cite{HoffSerre} where they proved that the assumption of constant viscosity
leads to a failure of continuous dependence on the initial data. Finally D.
Hoff extended the result of \cite{Hof87} by showing in \cite{Hof98} the
existence of global weak solution with initial density admitting shocks
(roughly speaking the density $\rho_{0}$ satisfies $\rho_{0},\frac{1}{\rho
_{0}}\in L^{\infty}$) and with regularizing effects on the velocity,

\paragraph{The density dependent case}

For initial densities that feature vacuum zones, existence of global
\textit{weak} solution has been obtained by Q.Jiu and Z. Xin in \cite{Jiu} for
viscosity coefficients verifying $\mu\left(  \rho\right)  =\rho^{\alpha}$ with
$\alpha>\frac{1}{2}$.

In \cite{mesure}, the second author constructed global weak solutions for
general viscosity coefficients with initial density admitting shock type
discontinuities and with initial velocity belonging to the set of finite
measures. In opposite to \cite{Hof98}, the initial data satisfy the BD entropy
but not the classical energy which allows in particular to obtain regularizing
effects for the density inasmuch as the density becomes instantaneously continuous.

We recall that the basic energy estimate, obtained formally by multiplying the
velocity's equation with $u$ and integrating by parts states that%
\begin{equation}
\int_{\mathbb{R}}\left(  \frac{1}{2}\rho u^{2}(t,x)+\rho e\left(  \rho\right)
(t,x)\right)  dx+\int_{0}^{t}\int_{\mathbb{R}}\mu(\rho)(\partial_{x}%
u)^{2}(s,x)dsdx\leq\int_{\mathbb{R}}\left(  \frac{1}{2}\rho_{0}u_{0}%
^{2}(t,x)+\rho_{0}e\left(  \rho_{0}\right)  (x)\right)  dx, \label{energy_NS}%
\end{equation}
with%
\[
\rho e^{\prime}\left(  \rho\right)  -e\left(  \rho\right)  =\rho^{\gamma}.
\]
In the context of the multidimensional NS system, D. Bresch and B. Desjardins
\cite{BD} discovered that $\left(  \text{\ref{Navier_Stokes_1d}}\right)  $ has
some "hidden" algebraic structure. We can observe (see \cite{MN}) that the
so-called \textit{effective velocity}:%
\[
v=u+\partial_{x}\varphi\left(  \rho\right)
\]
verifies the equation
\[
\partial_{t}(\rho v)+\partial_{x}(\rho uv)+\partial_{x}\rho^{\gamma}=0,
\]
and that multiplying with $v$ the previous relation the following functional,
named \textit{BD entropy}, is formally controlled:%
\begin{equation}
\int_{\mathbb{R}}\left(  \frac{1}{2}\rho v^{2}(t,x)+\rho e\left(  \rho\right)
(t,x)\right)  dx+\gamma\int_{0}^{t}\int_{\mathbb{R}}\mu\left(  \rho\right)
\rho^{\gamma-3}(\partial_{x}\rho)^{2}\leq\int_{\mathbb{R}}\left(  \frac{1}%
{2}\rho_{0}v_{0}^{2}(t,x)+\rho_{0}e\left(  \rho_{0}\right)  (x)\right)  dx,
\label{BD_entropy_NS}%
\end{equation}
We observe that if $\mu\left(  \rho\right)  =\rho^{\alpha}$ then the BD
entropy shows that%
\[
\int_{\mathbb{R}}\frac{1}{2}\rho(v-u)^{2}(t,x)dx=C_{\alpha}\int_{\mathbb{R}%
}(\partial_{x}\rho^{\alpha-\frac{1}{2}})^{2}(t,x)dx
\]
Thus, supposing that $\partial_{x}\rho_{0}^{\alpha-\frac{1}{2}}\in L^{2}$, and
taking $0<\alpha<\frac{1}{2}$, the BD entropy yields that the $L^{\infty}%
$-norm of $\frac{1}{\rho}$ remains bounded. This was the main argument that A.
Mellet and A. Vasseur \cite{MV06} used in order to construct global strong
solutions for the case of degenerate viscosity coefficients with $\alpha$
satisfying $0<\alpha<\frac{1}{2}$.

In \cite{MN}, the second author has proved similar results for the case
$1/2<\alpha\leq1$ where he exploited the fact that the effective velocity $v$
satisfies a damped transport equation. This allows him to obtain an
$L^{\infty}$-estimates for $v$ and using a maximum principle in order to
obtain a $L^{\infty}$ control on $\frac{1}{\rho}$. \newline

More recently Constantin \textit{et al} in \cite{Constantin} have extended the
previous result to the case $\alpha>1$ with $\gamma$ belonging to
$[\alpha,\alpha+1]$ provided that the initial data satisfy:
\begin{equation}
\partial_{x}u_{0}\leq\rho_{0}^{\gamma-\alpha}. \label{flux5}%
\end{equation}
We point out that the condition (\ref{flux5}) amounts to considering a
negative effective flux (see for example \cite{Hof87,Lio98}) at initial time.
The main idea of their proof consists in proving via a maximum principle that
the effective flux i.e. the function $\mu\left(  \rho\right)  \partial
_{x}u-\rho^{\gamma}$, remains negative for all time $t\geq0$ if this is the
case initially. This is sufficient to control the $L^{\infty}$ norm of
$\frac{1}{\rho}$.

In \cite{BurHas2020} we proved that if $\alpha>\frac{1}{2},$ $\gamma\geq
\max\{1,\alpha\}$ and if the so-called effective velocity satisfies initially
an Oleinik type inequality (see \cite{Oleinik}) then we have the existence of
global strong solution provided that $(\rho_{0},\frac{1}{\rho_{0}})$ are
bounded. In particular there is no restriction on the sign of the effective flux.

Let us mention that all the above results, assume that the BD entropy is
bounded initially, it allow sat least when $\alpha>\frac{1}{2}$ to control the
$L^{\infty}$ norm of the density $\rho$ all along the time by using Sobolev
embedding on the quantity $\partial_{x}\rho^{\alpha-\frac{1}{2}}\in L^{2}$.
Since this involves an $L^{2}$ information for the derivative of the density,
\textit{all the above results deal with continuous initial densities}.

Concerning discontinious initial densities we are only aware of two results:
Fang and Zhang \cite{fang2006global} dealing with the free boundary problem
and Ruxu et al \cite{ruxu2012cauchy} dealing with the problem on the whole
real line but assuming small initial energy. The estimates that the authors
obtain depend on the number of points of discontinuity and, in particular,
these estimates blow-up as the number of discontinuity points tends to
$+\infty$. For this reason, they are not appropiate in order to study
$1d-$homogeneisaition as in \cite{BreschHillairetENS}.

As it was mentioned previously, in this paper we will provide a new method to
control the $L^{\infty}$ norm of the density, without using the BD entropy. By
doing so, our initial density $\rho_{0}$ is just $L^{2}\left(  \mathbb{R}%
\right)  \cap L^{\infty}\left(  \mathbb{R}\right)  $ 

\subsubsection{The Navier-Stokes-Korteweg system \newline}

Let us pass now in review some results for system $\left(
\text{\ref{NSK_intro}}\right)  $. F. Charve and the second author in \cite{CH}
proved the global existence of strong solution for the system (\ref{NSK_intro}%
) when $\mu(\rho)=\varepsilon\rho$ and $\kappa(\rho)=\frac{\varepsilon^{2}%
}{\rho}$. In addition, they show that the global strong solutions converge
when $\varepsilon$ goes to $0$ to a global weak entropy solution of the
compressible Euler system with initial data of finite energy. P. Germain and
P.G. LeFloch in \cite{Germain} studied recently the global existence of vacuum
and non-vacuum weak solutions for the Korteweg system. It is important to
point out that they need to impose a tame condition on the viscosity and
capillary coefficients which takes the form:
\begin{equation}
\kappa(\rho)\lesssim\frac{\mu(\rho)^{2}}{\rho^{3}}\;\;\mbox{and}\;\;\delta
(\varepsilon)\lesssim\varepsilon^{2}, \label{tame}%
\end{equation}
if we consider the vanishing viscosity capillary limit for viscosity and
capillary coefficients $\mu_{\varepsilon}(\rho)=\varepsilon\mu(\rho)$,
$\kappa_{\varepsilon}(\rho)=\delta(e)\kappa(\rho)$ when $\varepsilon>0$ goes
to $0$. Roughly speaking, the previous tame condition implies in some sense
that the parabolic behavior, governed by the viscosity tensor dominates the
dispersive effects that are induced by the capillarity tensor. We mention that
this is important because the dispersive effects tend to create strong
oscillations which prevent obtaining strong convergence informations (they
allow only to obtain weak-convergence results which are not well-suited to
treat the compressible setting). In particular, when $\varepsilon$ goes to
$0$, by assuming $\left(  \text{\ref{tame}}\right)  $ we can expect recovering
strong convergence in suitable functional spaces. In \cite{Germain} the
authors study then the zero viscosity-capillarity limit associated with the
Navier-Stokes-Korteweg system generalizing the results of \cite{CH}. They need
in particular for this analysis to assume that $\kappa_{\varepsilon}%
(\rho)=\frac{\mu_{\varepsilon}(\rho)^{2}}{\rho^{3}}$ when the viscosity
coefficient is degenerate with $\mu_{\varepsilon}(\rho)=\varepsilon
\rho^{\alpha}$. We recall that we will also consider this algebraic relation
in the remaining part of this paper. \newline Recently, Chen et al. in
\cite{Chen2} and Chen in \cite{Chen1} have proved for the first time some
results of existence of global strong solutions for initial density far away
from the vacuum in Lagrangian coordinates. More precisely they consider
viscosity and capillary coefficients of the form $\mu(\rho)=\rho^{\alpha_{1}}$
and $\kappa(\rho)=\rho^{\beta_{1}}$ with $(\alpha_{1},\beta_{1})\in
\mathbb{R}^{2}$. In comparison with the present work, there is no relation a
priori between $\alpha_{1}$ and $\beta_{1}$ (furthermore there is no
restriction on $c>0$). They manage essentially to show such result when
$\beta_{1}<-2$ which allows in a direct way to control the $L^{\infty}$ norm
of $\frac{1}{\rho}$ by using the energy estimate (indeed roughly speaking the
energy estimate ensure that $\partial_{x}\rho^{\frac{\beta_{1}}{2}+1}$ is
bounded in $L_{T}^{\infty}(L^{2})$ for any $T>0$). They deal also with the
case $\beta_{1}\geq-2$ but in this case $\alpha_{1}<0$, in particular the
viscosity coefficient explodes near vacuum. The main ideas of the proof is to
obtain $L^{2}$-estimates for the effective velocity $v=u+\frac{\mu\left(
\rho\right)  }{\rho^{2}}\partial_{x}\rho$ by using an energy method combined
with Sobolev embedding in the spirit of Kanel (see \cite{Ka}). Furthermore the
authors show also the existence of global strong solution when the initial
data is a perturbation of a Riemann problem associated to a rarefaction wave
for the compressible Euler problem.

In \cite{BurHas2021} the authors extend the results of \cite{Chen1,Chen2} to
the case of the NSK system with strongly degenerate viscosity coefficients
with initial density far away from vacuum. More precisely, the viscosity
coefficients take the form $\mu(\rho)=\rho^{\alpha}$ with $\alpha>1$ and the
capillarity coefficient satisfies the algebraic relation (\ref{relatie}). The
main difficulty of the proof consists in estimating globally in time the
$L^{\infty}$ norm of $\frac{1}{\rho}$. The method of proof relies on fine
algebraic properties of the NSK system. In \cite{BurHas2021} we introduce two
new effective velocities endowed with weight functions depending both on the
viscosity and the capillarity coefficients as some power laws of the density.
For these two quantities we show some Oleinik-type estimate which provide the
control of the $L^{\infty}$ norm of $\frac{1}{\rho}$ by applying a maximum
principle. It is interesting to point out that the two effective pressure
introduced in this paper depending on the capillary coefficient generalize to
the NSK systems those introduced for the NS system in
\cite{BurHas2020,Constantin}.

\subsection{Main results}

\subsubsection{The algebraic structure of the NSK
system\label{algebrai_structure_Section}}

For the rest of the paper, as in \cite{Germain,BrGisLac}, we will assume
that:
\begin{equation}
\kappa\left(  \rho\right)  =\frac{\mu^{2}\left(  \rho\right)  }{\rho^{3}}.
\label{relatie}%
\end{equation}
With such a choice for $\kappa$ we can rewrite the capillarity tensor $K$ as%
\[
K=\mu\left(  \rho\right)  \partial_{xx}\varphi\left(  \rho\right)
\]
with
\begin{equation}
\varphi^{\prime}\left(  \rho\right)  =\frac{\mu\left(  \rho\right)  }{\rho
^{2}}. \label{def_varphi}%
\end{equation}
We observe also that:
\begin{equation}
\partial_{x}K=\rho\partial_{x}\left(  G^{\prime}\left(  \rho\right)
\partial_{xx}^{2}G\left(  \rho\right)  \right)  \text{ with }G^{\prime}\left(
\rho\right)  =\frac{\mu\left(  \rho\right)  }{\rho^{3/2}}. \label{def_G}%
\end{equation}
We will study the system (\ref{NSK_intro}) on the real line $\mathbb{R}$ with
the following far field assumption:
\begin{equation}
\rho\left(  t,x\right)  \rightarrow1\text{ and }u\left(  t,x\right)
\rightarrow0\text{ when }\left\vert x\right\vert \rightarrow\infty,
\label{conditions_at_infinity}%
\end{equation}
for all $t\geq0$. Let $r\in\left(  0,1\right)  $ and observe that using the
second equation of the system $\left(  \text{\ref{NSK_intro}}\right)  $ we can
write that:%
\begin{align*}
&  \rho\partial_{t}\left(  u+r\partial_{x}\varphi\left(  \rho\right)  \right)
+\rho u\partial_{x}\left(  u+r\partial_{x}\varphi\left(  \rho\right)  \right)
-\left(  1-r\right)  \partial_{x}\left(  \mu(\rho)\partial_{x}\left(
u+r\partial_{x}\varphi\left(  \rho\right)  \right)  \right)  +a\partial
_{x}\rho^{\gamma}\\
&  =\left(  r^{2}-r+c\right)  \partial_{x}K.
\end{align*}
For $c\in(0,\frac{1}{4}]$ the equation
\[
r^{2}-r+c=0
\]
admits the two positive roots%
\[
r_{1}\left(  c\right)  =\frac{1+\sqrt{1-4c}}{2}\leq1\text{ and }r_{0}\left(
c\right)  =\frac{1-\sqrt{1-4c}}{2}\leq1.
\]
In the sequel we will assume that $c\in(0,\frac{1}{4}]$ and we observe that
\begin{equation}
r_{1}\left(  c\right)  \underset{c\rightarrow0}{\longrightarrow}1\text{ and
}r_{0}\left(  c\right)  \underset{c\rightarrow0}{\longrightarrow}0.
\label{limite}%
\end{equation}
We introduce the following two extra variables
\[
v_{i}=u+r_{i}\partial_{x}\varphi\left(  \rho\right)  \text{ for }i\in\left\{
0,1\right\}  ,
\]
which will be referred as \textit{effective velocities}. The effective
velocities verify the following equations:%
\begin{equation}
\rho\partial_{t}v_{i}+\rho u\partial_{x}v_{i}-\left(  1-r_{i}\left(  c\right)
\right)  \partial_{x}\left(  \mu(\rho)\partial_{x}v_{i}\right)  +a\partial
_{x}\rho^{\gamma}=0\text{ for }i\in\left\{  0,1\right\}  .
\label{equations_for_v}%
\end{equation}
Let us discuss now the dissipation of energy. Using (\ref{def_G}) the natural
energy associated to $\left(  \text{\ref{NSK_intro}}\right)  $ is%
\begin{equation}
\int_{\mathbb{R}}\big(\frac{1}{2}\rho u^{2}(t,x)+\rho e\left(  \rho\right)
(t,x)+\frac{1}{2}c\left(  \partial_{x}G\left(  \rho\right)  \right)
^{2}(t,x)\big)dx+\int_{0}^{t}\int_{\mathbb{R}}\mu(\rho)(\partial_{x}%
u)^{2}(s,x)dsdx\leq E_{c}\left(  \rho_{0},u_{0}\right)  , \label{basic_energy}%
\end{equation}
where%
\begin{equation}
E_{c}\left(  \rho_{0},u_{0}\right)  =\frac{1}{2}%
{\displaystyle\int_{\mathbb{R}}}
\big( \rho_{0}u_{0}^{2}(x)+\rho_{0}e\left(  \rho_{0}\right)  (x)+cE_{cap}%
\left(  \partial_{x}\rho_{0}\right)  \big)dx ,\text{ } \label{Ec_t=0}%
\end{equation}
with
\[
E_{cap}\left(  \partial_{x}\rho_{0}\right)  =\dfrac{1}{2}%
{\displaystyle\int_{\mathbb{R}}}
\left(  \partial_{x}G\left(  \rho_{0}\right)  \right)  ^{2}(x)dx=\frac{1}{2}%
{\displaystyle\int_{\mathbb{R}}}
\rho_{0}\left(  \partial_{x}\varphi\left(  \rho_{0}\right)  \right)
^{2}(x)dx
\]
and%
\begin{equation}
e\left(  \rho\right)  =\frac{\left(  \rho^{\gamma}-1-\gamma\left(
\rho-1\right)  \right)  }{\left(  \gamma-1\right)  \rho}=\frac{\rho^{\gamma
-1}}{\gamma-1}+\frac{1}{\rho}-\frac{\gamma}{\gamma-1}.
\label{potential_energy}%
\end{equation}
In the context of the Navier-Stokes-Korteweg system $\left(
\text{\ref{NSK_intro}}\right)  $, since there are two effective-velocities we
have two BD-type entropies obtained by multiplying the equation $\left(
\text{\ref{equations_for_v}}\right)  $ with $v_{i}$ and integrating by parts:
\begin{align}
&  \int_{\mathbb{R}}\left(  \frac{1}{2}\rho v_{i}^{2}(t,x)+\rho e\left(
\rho\right)  (t,x)\right)  dx+\left(  1-r_{i}\left(  c\right)  \right)
\int_{0}^{t}\int_{\mathbb{R}}\mu(\rho)(\partial_{x}v_{i})^{2}(s,x)dsdx+r_{i}%
\left(  c\right)  \gamma\int_{0}^{t}\int_{0}^{t}\mu\left(  \rho\right)
\rho^{\gamma-3}(\partial_{x}\rho)^{2}(s,x)dsdx\\
&  \leq\int_{\mathbb{R}}\big(\frac{1}{2}\rho v_{i,0}^{2}(x)+\rho_{0}e\left(
\rho_{0}\right)  (x)\big)dx=\int_{\mathbb{R}}\big(\frac{1}{2}\rho_{0}u_{0}%
^{2}(x)+\rho_{0}e\left(  \rho_{0}\right)  (x)\big)dx+r_{i}\int_{\mathbb{R}%
}\rho_{0}u_{0}\partial_{x}\varphi\left(  \rho_{0}\right)  (x)dx+r_{i}%
^{2}E_{cap}\left(  \partial_{x}\rho_{0}\right)  . \label{BD_entropies}%
\end{align}

If one ignores the relative size with respect to $c$ of the coefficients
$r_{i}\left(  c\right)  ,$ both estimates $\left(  \text{\ref{BD_entropies}%
}\right)  $ give the same qualitative information. Moreover, the spaces where
the initial data has to be drawn from, in order to give a mathematical meaning
to these functionals are such that $E_{c}\left(  \rho_{0},u_{0}\right)  $
defined in $\left(  \text{\ref{Ec_t=0}}\right)  $ is finite. In particular,
\textit{assuming finite energy, the density in the NSK system is always a
continuous function}. Of course, this is normal since, as discussed in the
Introduction, NSK is a Diffuse Interface Model. However, observe that as
$c\rightarrow0$, the estimate for $v_{0}$ degenerates in the basic energy
estimate for Navier-Stokes $\left(  \text{\ref{energy_NS}}\right)  $, while
the estimate for $v_{1}$ degenerates to the BD-entropy estimate $\left(
\text{\ref{BD_entropy_NS}}\right)  $.\ This observation reveals that if one is
interested in recovering a "rough" solution of the NS\ system by a vanishing
capillarity limit, \textit{the BD entropy estimate for }$v_{1}$\textit{ must
explode as }$c\rightarrow0$.

\subsubsection{Assumptions for the viscosity coefficient}

We wish now to specify the form of the viscosity coefficient with which we
will work. We recall that%
\begin{equation}
\varphi^{\prime}\left(  \rho\right)  =\frac{\mu\left(  \rho\right)  }{\rho
^{2}}\text{ and }\psi^{\prime}\left(  \rho\right)  =\frac{\mu\left(
\rho\right)  }{\rho} \cdot\label{definition_varphi_psi}%
\end{equation}
Also, we introduce
\begin{equation}
\Xi^{\prime}\left(  \rho\right)  =\frac{\left(  \rho^{\gamma}\mu\left(
\rho\right)  \right)  ^{\frac{1}{2}}}{\rho^{2}} \cdot\label{definition_Xi}%
\end{equation}
We suppose that $\mu$ is a $C^{1}$ positive function on $(0,+\infty)$ such
that:
\begin{equation}
\lim_{\rho\rightarrow+\infty}\mu\left(  \rho\right)  =\infty\tag{$\mu
_1$}\label{H1}%
\end{equation}
We assume that there exists positive constants $d_{1},d_{2},d_{3},d_{4},d_{5}$
such that the following relations hold true. First of all, for all $\rho\geq0$
we suppose that%
\begin{equation}
\frac{1}{d_{1}}\mu(\rho)\leq\psi\left(  \rho\right)  \leq d_{1}\mu\left(
\rho\right)  . \tag{$\mu_2$}\label{H2}%
\end{equation}
Next, we suppose that%

\begin{equation}
\mu\left(  \rho\right)  \leq d_{2}\left(  1+\rho^{\gamma}\right)  \text{ for
all }\rho\geq0\text{.} \tag{$\mu_3$}\label{H3}%
\end{equation}
We also require that $\varphi:(0,\infty)\rightarrow\mbox{Im}\varphi$ is
invertible and that
\begin{equation}
\lim_{\rho\rightarrow0}\varphi\left(  \rho\right)  =-\infty. \tag{$\mu
_4$}\label{H4}%
\end{equation}
We impose that
\begin{equation}
\lim_{\rho\rightarrow+\infty}\Xi\left(  \rho\right)  =+\infty. \tag{$\mu
_5$}\label{H5}%
\end{equation}
Finally, we denote
\begin{equation}
\Phi\left(  \tau\right)  =-\varphi\left(  \frac{1}{\tau}\right)  >0\text{ ,
}\Lambda\left(  \tau\right)  =\frac{1}{\mu\left(  \dfrac{1}{\tau}\right)  }>0
\label{Var_Phi}%
\end{equation}
and we ask that%
\begin{equation}
\Lambda\left(  \tau\right)  \leq d_{4}\left(  1+\Phi\left(  \tau\right)
\right)  ^{2} \tag{$\mu_6$}\label{H6}%
\end{equation}
along with
\begin{equation}
\tau^{\frac{1}{2}}\leq d_{5}\left(  1+\Phi\left(  \tau\right)  \right)
^{1-\eta} \tag{$\mu_7$}\label{H7}%
\end{equation}
for some positive $\eta\in(0,1)$.

\begin{remarka}
\label{rem1}Let us investigate how the above hypothesis translate in the case
of viscosity coefficient of the form $\mu\left(  \rho\right)  =\rho^{\alpha}$.
Observe that in this case we have that%
\[
\varphi(\rho)=\frac{\rho^{\alpha-1}}{\alpha-1}\;\;\mbox{and}\;\;\psi
(\rho)=\frac{\rho^{\alpha}}{\alpha}\text{ if }\alpha\not =1
\]
while%
\[
\varphi(\rho)=\ln\rho\text{ and }\psi(\rho)=\rho\text{ if }\alpha=1.
\]
We see that, $\left(  \text{\ref{H1}}\right)  $, $\left(  \text{\ref{H2}%
}\right)  $, $\left(  \text{\ref{H3}}\right)  $ imply that%
\[
0<\alpha\leq\gamma.
\]
Hypothesis $\left(  \text{\ref{H4}}\right)  $ implies that
\[
\alpha\leq1.
\]
Hypothesis $\left(  \text{\ref{H5}}\right)  $ implies that%
\[
\alpha+\gamma\geq2.
\]
Finally, since
\[
\Lambda\left(  \tau\right)  =\tau^{\alpha}\text{ and }\Phi\left(  \tau\right)
=\tfrac{1}{1-\alpha}\tau^{1-\alpha}\text{ if }\alpha\not =1,
\]
the hypothesis $\left(  \text{\ref{H6}}\right)  $ implies that
\[
\alpha\leq\frac{2}{3}
\]
while $\left(  \text{\ref{H7}}\right)  $ implies that%
\[
\alpha<\frac{1}{2}.
\]

\end{remarka}

\subsubsection{Statement of the main results}

Finally, after a rather long introduction, we are in the position of stating
our first main result which addresses the question of global existence of
strong solutions with rough initial densities.

\begin{theorem}
\label{theo3} Consider $\gamma>1$ and $\mu$ verifying $\left(  \text{\ref{H1}%
}\right)  $, $\left(  \text{\ref{H2}}\right)  $, $\left(  \text{\ref{H3}%
}\right)  $, $\left(  \text{\ref{H4}}\right)  $,$\left(  \text{\ref{H6}%
}\right)  $ and $\left(  \text{\ref{H7}}\right)  $. Consider $\left(  \rho
_{0},u_{0}\right)  \in L^{\infty}\left(  \mathbb{R}\right)  \times
L^{2}\left(  \mathbb{R}\right)  $ such that%
\[
\left\Vert \rho_{0}\right\Vert _{L^{\infty}(\mathbb{R})}+\left\Vert \frac
{1}{\rho_{0}}\right\Vert _{L^{\infty}(\mathbb{R})}+E_{0}\left(  \rho_{0}%
,u_{0}\right)  \leq M.
\]
Then, there exists a unique global solution $\left(  \rho,u\right)  $ for the
Navier-Stokes system $\left(  \text{\ref{Navier_Stokes_1d}}\right)  $.
Moreover, there exists $C=C\left(  T,E_{0}\left(  \rho_{0},u_{0}\right)
,\left\Vert \rho_{0},\frac{1}{\rho_{0}}\right\Vert _{L^{\infty}}\right)  $
such that the pair $\left(  \rho,u\right)  $ satisfies for all $T>0$:%
\begin{equation}
E_{0}\left(  \rho\left(  T\right)  ,u\left(  T\right)  \right)  =%
{\displaystyle\int_{\mathbb{R}}}
\left(  \rho u^{2}+\rho e\left(  \rho\right)  \right)  \left(  T\right)  +%
{\displaystyle\int_{0}^{T}}
{\displaystyle\int_{\mathbb{R}}}
\mu(\rho)(\partial_{x}u)^{2}\leq E_{0}\left(  \rho_{0},u_{0}\right)  ,
\end{equation}%
\begin{equation}
\frac{1}{C}\leq\rho(T,\cdot)\leq C,
\end{equation}%
\begin{equation}
\int_{0}^{T}\int_{\mathbb{R}}\sigma\rho\left\vert \dot{u}\right\vert
^{2}+\frac{1}{2}\sigma\left(  T\right)  \int_{\mathbb{R}}\mu\left(
\rho\left(  T\right)  \right)  \left(  \partial_{x}u\left(  T\right)  \right)
^{2}\leq C,
\end{equation}%
\begin{equation}
\frac{1}{2}\int_{0}^{T}\int_{\mathbb{R}}\sigma^{2}\rho\left\vert \dot
{u}\right\vert ^{2}+\left(  1-r_{0}\right)  \int_{0}^{T}\int_{\mathbb{R}%
}\sigma^{2}\mu\left(  \rho\right)  \left\vert \partial_{x}\dot{u}\right\vert
^{2}\leq C,
\end{equation}%
\begin{equation}
\int_{0}^{T}\sigma^{\frac{1}{2}}\left(  \tau\right)  \left\Vert \partial
_{x}u\left(  \tau\right)  \right\Vert _{L^{\infty}}^{2}d\tau+\sup_{0<t\leq
T}\sigma(t)\Vert\partial_{x}u(t,\cdot)\Vert_{L^{\infty}}\leq C,
\end{equation}
where $\sigma\left(  t\right)  =\min\left\{  1,t\right\}  $ and $\dot
{u}=\partial_{t} u+u\partial_{x} u$.
\end{theorem}

\begin{remarka}
Up to your knowledge this is the first result assuring the existence of global
unique solutions for the compressible Navier-Stokes system in one dimension
with degenerate viscosity coefficient and with rough initial density $\rho
_{0}\in L^{\infty}$. Of course, this includes the case of discontinuous
initial data, with an arbitrary number of discontinuities and arbitrary jumps.
In particular, this extends the works \cite{BurHas2020,Constantin,MV} where
the authors consider initial densities which are continuous.
\end{remarka}

\begin{remarka}
According to the Remark \ref{rem1}, we recall that the choice $\mu(\rho
)=\rho^{\alpha}$ with $\alpha\in(0,1/2)$ enters in the framework of Theorem
\ref{theo3}.\newline Since we do not need to assume $\left(  \text{\ref{H5}%
}\right)  $, any $p\left(  \rho\right)  =\rho^{\gamma}$ with $\gamma>1$ fits
in the framework of our result.
\end{remarka}

\begin{remarka}
At this point it is important to remark that, using Proposition
\ref{rho_is_bounded_in_Linfty} proved below, the result obtained by Constantin
et al. \cite{Constantin} can be adapted to obtain an existence result with
discontinuous initial densities $\rho_{0}$ such that%
\begin{equation}
\partial_{x}u_{0}-\rho_{0}^{\gamma-\alpha}\leq0. \label{restriction}%
\end{equation}
Let us remark however that, at least if $\partial_{x}u_{0}$ is continuous,
condition that although cannot be ensured for later times, it nevertheless
implies that the jump of $\rho_{0}^{\gamma-\alpha}$ at a discontinuity must be
positive. This observation shows that general initial data relevant for
multifluids at a mesoscopic scale, as explained in \cite{BrDeGhGrHi2018},
Section $3.2.$ cannot be considered assuming $\left(  \text{\ref{restriction}%
}\right)  $ initially.
\end{remarka}

The second main result states that, if we impose a sign condition on the
effective velocity $v_{1}$ then, we can construct solutions for the NSK system
$\left(  \text{\ref{NSK_intro}}\right)  $ which verify uniform estimates for
the density and the effective velocity $v_{0}$ \textit{that do not degenerate
as }$c\rightarrow0$. Recall that formally, when $c\rightarrow0,$ $v_{0}$
degenerates to $u,$ the velocity of the NS system.

\begin{theorem}
Let $c\in(0,\frac{1}{4})$, $\gamma>1$ and $\mu,\varphi,\psi$ verifying
$\left(  \text{\ref{definition_varphi_psi}}\right)  $, $\left(  \text{\ref{H1}%
}\right)  $-$\left(  \text{\ref{H7}}\right)  $. Consider the initial data
$\left(  \rho_{0},u_{0}\right)  $ which is uniformly bounded in $c$ in the
following space%
\[
u_{0},\sqrt{c}\rho_{0}^{\frac{1}{2}}\partial_{x}\varphi\left(  \rho
_{0}\right)  \in L^{2}\left(  \mathbb{R}\right)  ,\rho_{0}e\left(  \rho
_{0}\right)  \in L^{1}\left(  \mathbb{R}\right)  .
\]
Moreover, suppose that there exists $M_{0}\in\mathbb{R}$ such that for all
$c\in\left(  0,\frac{1}{4}\right)  $ we have
\[
\text{a.e. }x\in\mathbb{R}\text{ \ }:\text{ \ }v_{1|t=0}=u_{0}\left(
x\right)  +r_{1}\left(  c\right)  \partial_{x}\varphi(\rho_{0}\left(
x\right)  )\leq M_{0}.
\]
Then, there exists a unique global solution for $\left(  \text{\ref{NSK_intro}%
}\right)  $ that satisfies the following estimates, uniformly with respect to
the parameter $c$:
\begin{equation}
\left.
\begin{array}
[c]{l}%
{\displaystyle\int_{\mathbb{R}}}
\frac{1}{2}\rho u^{2}+\rho e\left(  \rho\right)  +c\rho\left(  \partial
_{x}\varphi\left(  \rho\right)  \right)  ^{2}+%
{\displaystyle\int_{0}^{t}}
{\displaystyle\int_{\mathbb{R}}}
\mu(\rho)(\partial_{x}u)^{2}\leq E_{c}\left(  \rho_{0},u_{0}\right)  ,\\
\\%
{\displaystyle\int_{\mathbb{R}}}
\frac{1}{2}\rho v_{0}^{2}+\rho e\left(  \rho\right)
{\displaystyle\int_{0}^{t}}
{\displaystyle\int_{\mathbb{R}}}
+\left(  1-r_{0}\left(  c\right)  \right)
{\displaystyle\int_{0}^{t}}
{\displaystyle\int_{\mathbb{R}}}
\mu(\rho)(\partial_{x}v_{i})^{2}(s,x)dsdx+r_{0}\left(  c\right)  \gamma%
{\displaystyle\int_{0}^{t}}
{\displaystyle\int_{\mathbb{R}}}
\mu\left(  \rho\right)  \rho^{\gamma-3}(\partial_{x}\rho)^{2}(s,x)dsdx\leq
2E_{c}\left(  \rho_{0},u_{0}\right)  ,\\
\\
\rho\left(  t,x\right)  \leq C\left(  t,E_{c},\left\Vert \rho_{0}\right\Vert
_{L^{\infty}}\right)  ,\\
\\
u\left(  t,x\right)  +r_{1}\left(  c\right)  \partial_{x}\varphi\left(
\rho\left(  t,x\right)  \right)  \leq M_{0}+C\left(  t,E_{c},\left\Vert
\rho_{0}\right\Vert _{L^{\infty}}\right)  ,\\
\\
\Vert\dfrac{1}{\rho\left(  t,\cdot\right)  }\Vert_{L^{\infty}}\leq C\left(
t,E_{c},M_{0},\left\Vert \rho_{0},\frac{1}{\rho_{0}}\right\Vert _{L^{\infty}%
}\right) \\
\\
\left\Vert \varphi\left(  \rho\right)  \right\Vert _{BV(\left[  0,T\right]
\times\left[  -L,L\right]  )}+\left\Vert \rho\right\Vert _{BV(\left[
0,T\right]  \times\left[  -L,L\right]  )}\leq C\left(  t,L,E_{c}%
,M_{0},\left\Vert \rho_{0},\frac{1}{\rho_{0}}\right\Vert _{L^{\infty}}\right)
\text{ }\forall L>0.\\
\\%
{\displaystyle\int_{0}^{T}}
{\displaystyle\int_{\mathbb{R}}}
\sigma\rho\left\vert \dot{v}_{0}\right\vert ^{2}+\dfrac{1-r_{0}\left(
c\right)  }{2}\sigma\left(  T\right)
{\displaystyle\int_{\mathbb{R}}}
\mu\left(  \rho\left(  T\right)  \right)  \left(  \partial_{x}v_{0}\left(
T\right)  \right)  ^{2}\leq C\left(  t,E_{c},M_{0},\left\Vert \rho_{0}%
,\frac{1}{\rho_{0}}\right\Vert _{L^{\infty}}\right)  ,\\
\\
\frac{1}{2}%
{\displaystyle\int_{0}^{T}}
{\displaystyle\int_{\mathbb{R}}}
\sigma^{2}\rho\left\vert \dot{v}_{0}\right\vert ^{2}+\left(  1-r_{0}\left(
c\right)  \right)
{\displaystyle\int_{0}^{T}}
{\displaystyle\int_{\mathbb{R}}}
\sigma^{2}\mu\left(  \rho\right)  \left\vert \partial_{x}\dot{v}%
_{0}\right\vert ^{2}\leq C\left(  t,E_{c},M_{0},\left\Vert \rho_{0},\frac
{1}{\rho_{0}}\right\Vert _{L^{\infty}}\right)  ,
\end{array}
\right.  \label{estimtheo1}%
\end{equation}
where $\sigma\left(  t\right)  =\min\left\{  1,t\right\}  $ and $\dot{v_{0}%
}=\partial_{t}v_{0}+u\partial_{x}v_{0}$. Above, $C$ are generic functions
depending continuously on their arguments and they are increasing in $t$.
\label{theo1}
\end{theorem}

As we pointed out in the introduction, comparing with \cite{Germain} and
\cite{Chen1,Chen2}, our contribution consists in the uniform estimates with
respect to $c$ for $\left\Vert \rho\right\Vert _{L^{\infty}}$ and the fact
that we only use the one sided bound for $v_{1|t=0}=u_{0}+r_{1}\partial
_{x}\varphi\left(  \rho_{0}\right)  $ in order to show that $\rho$ is bounded
from bellow. Of course this is crucial if one is interested in the
capillarity-vanishing limit problem (the viscosity remaining constant).

As a corollary to Theorem \ref{theo1} we immediately obtain

\begin{theorem}
\label{theo2} Let $c\in(0,\frac{1}{4})$, $\gamma>1$ and $\mu,\varphi,\psi$
verifying $\left(  \text{\ref{definition_varphi_psi}}\right)  $, $\left(
\text{\ref{H1}}\right)  $-$\left(  \text{\ref{H6}}\right)  $. Let $M_{0}%
\in\mathbb{R}$, $M>0$ and consider $\left(  \rho_{0}^{c},u_{0}^{c}\right)  $
such that%
\[
\left\Vert \rho_{0}^{c}\right\Vert _{L^{\infty}(\mathbb{R})}+\left\Vert
\frac{1}{\rho_{0}^{c}}\right\Vert _{L^{\infty}(\mathbb{R})}+E_{c}\left(
\rho_{0}^{c},u_{0}^{c}\right)  \leq M,
\]
and%
\[
u_{0}^{c}\left(  x\right)  +r_{1}\left(  c\right)  \partial_{x}\varphi
(\rho_{0}^{c})\left(  x\right)  \leq M_{0}\text{ a.e. on }\mathbb{R}\text{.}%
\]
Moreover, suppose that $(\rho_{0},u_{0})$ is such that $\rho_{0}e(\rho_{0})\in
L^{1}(\mathbb{R})$, $u_{0}\in L^{2}(\mathbb{R})$, $(\rho_{0},1/\rho_{0})\in
L^{\infty}(\mathbb{R})$, $\varphi(\rho_{0})\in BV_{loc}\left(  \mathbb{R}%
\right)  $ such that
\begin{equation}
u_{0}+\partial_{x}\varphi(\rho_{0})\leq M_{0}\;\;\;\text{in the sense of
measures } \label{mesure}%
\end{equation}
and that:
\begin{equation}%
\begin{cases}
\rho_{0}^{c}-1\rightarrow\rho_{0}-1\text{, }u_{0}^{c}\rightarrow u_{0}\text{
in }L^{2}\left(  \mathbb{R}\right)  ,\\
u_{0}^{c}+r_{1}\partial_{x}\varphi(\rho_{0}^{c})\rightharpoonup u_{0}%
+\partial_{x}\varphi(\rho_{0})\text{ weakly in the sense of measures.}%
\end{cases}
\end{equation}
Then, there exists $\left(  \rho,u\right)  $ such that for all $T>0$,
\[
\lim_{c\rightarrow0}\left(  \rho^{c},u^{c}\right)  _{c>0}=\left(
\rho,u\right)  \text{ weakly}-\star\text{ in }L^{\infty}\left(  \left[
0,T\right]  \times\mathbb{R}\right)  \text{ }%
\]
and $\left(  \rho,u\right)  $ is the unique global solution for the NS system
(\ref{Navier_Stokes_1d}) with initial data $(\rho_{0},u_{0})$ that verifies
for all $T>0$:%
\begin{equation}
E_{0}\left(  \rho\left(  T\right)  ,u\left(  T\right)  \right)  =%
{\displaystyle\int_{\mathbb{R}}}
\left(  \rho u^{2}+\rho e\left(  \rho\right)  \right)  \left(  T\right)  +%
{\displaystyle\int_{0}^{T}}
{\displaystyle\int_{\mathbb{R}}}
\mu(\rho)(\partial_{x}u)^{2}\leq E_{0}\left(  \rho_{0},u_{0}\right)  ,
\end{equation}%
\begin{equation}
C\left(  t,E_{0}\left(  \rho_{0},u_{0}\right)  ,\left\Vert \rho_{0},\frac
{1}{\rho_{0}}\right\Vert _{L^{\infty}}\right)  ^{-1}\leq\rho(t,x)\leq C\left(
t,E_{0}\left(  \rho_{0},u_{0}\right)  ,\left\Vert \rho_{0},\frac{1}{\rho_{0}%
}\right\Vert _{L^{\infty}}\right)  \;\;\mbox{for any}\; x\in\mathbb{R}.
\end{equation}%
\begin{equation}
\left\Vert \varphi\left(  \rho\right)  \right\Vert _{BV(\left[  0,T\right]
\times\left[  -L,L\right]  )}+\left\Vert \rho\right\Vert _{BV(\left[
0,T\right]  \times\left[  -L,L\right]  )}\leq C\left(  t,L,E_{c}%
,M_{0},\left\Vert \rho_{0},\frac{1}{\rho_{0}}\right\Vert _{L^{\infty}}\right)
\text{ }\forall L>0.\\
\end{equation}%
\begin{equation}
u+\partial_{x}\varphi\left(  \rho\right)  \leq M_{0}+C\left(  t,E_{0}%
,\left\Vert \rho_{0}\right\Vert _{L^{\infty}}\right)  \text{ in the sense of
measures}%
\end{equation}%
\begin{equation}
\int_{0}^{T}\int_{\mathbb{R}}\sigma\rho\left\vert \dot{u}\right\vert
^{2}+\frac{1}{2}\sigma\left(  T\right)  \int_{\mathbb{R}}\mu\left(
\rho\left(  T\right)  \right)  \left(  \partial_{x}u\left(  T\right)  \right)
^{2}\leq C=C\left(  T,E_{0}\left(  \rho_{0},u_{0}\right)  ,\left\Vert \rho
_{0},\frac{1}{\rho_{0}}\right\Vert _{L^{\infty}}\right)  ,
\end{equation}%
\begin{equation}
\frac{1}{2}\int_{0}^{T}\int_{\mathbb{R}}\sigma^{2}\rho\left\vert \dot
{u}\right\vert ^{2} \int_{0}^{T}\int_{\mathbb{R}}\sigma^{2}\mu\left(
\rho\right)  \left\vert \partial_{x}\dot{u}\right\vert ^{2}\leq C=C\left(
T,E_{0}\left(  \rho_{0},u_{0}\right)  ,\left\Vert \rho_{0},\frac{1}{\rho_{0}%
}\right\Vert _{L^{\infty}}\right)  ,
\end{equation}%
\begin{equation}
\int_{0}^{T}\sigma^{\frac{1}{2}}\left(  \tau\right)  \left\Vert \partial
_{x}u\left(  \tau\right)  \right\Vert _{L^{\infty}}^{2}d\tau+\sup_{0<t\leq
T}\sigma(t)\Vert\partial_{x}u(t,\cdot)\Vert_{L^{\infty}}\leq C=C\left(
T,E_{0}\left(  \rho_{0},u_{0}\right)  ,\left\Vert \rho_{0},\frac{1}{\rho_{0}%
}\right\Vert _{L^{\infty}}\right)  ,
\end{equation}
.
\end{theorem}

As far as we are aware, this is the first result treating the vanishing
capillarity limit. As a corollary of the previous Theorem we obtain via a
vanishing capillary process the existence of global strong solutions for the
$1D$ Navier-Stokes system with discontinuous initial densities . However,
besides the fact that $\rho_{0}$ and $\frac{1}{\rho_{0}}$ must be in
$L^{\infty}\left(  \mathbb{R}\right)  $ we also require that%
\begin{equation}
u_{0}+\partial_{x}\varphi(\rho_{0})\leq M_{0}\; \label{BVloc}%
\end{equation}
for some $M_{0}\in\mathbb{R}$, this is obviously an additional condition
compared to Theorem \ref{theo3}. In particular (\ref{BVloc}) implies that
$\varphi(\rho_{0})$ belongs to $BV_{loc}$ and that the only discontinuities in
$x$ on $\rho_{0}$ are such as the jump of $\varphi\left(  \rho_{0}\right)  $
at $x$ is negative:
\[
\lbrack\varphi\left(  \rho_{0}\right)  ]\left(  x\right)  =\lim_{h\rightarrow
0,h>0}\varphi\left(  \rho_{0}\left(  x+h\right)  \right)  -\lim_{h\rightarrow
0,h>0}\varphi\left(  \rho_{0}\left(  x-h\right)  \right)  \leq0.
\]

\begin{remarka}
We would like also emphasize that it is possible to obtain a similar theorem
if we replace the condition (\ref{mesure}) by:
\begin{equation}
u_{0}+\partial_{x}\varphi\left(  \rho_{0}\right)  \geq M_{0}\text{ in the
sense of measures,} \label{mesure1}%
\end{equation}
with $M_{0}\in\mathbb{R}$. We refer to Section
\ref{lower_bound_for_density_section} and Remark \ref{rem5} for more details
on these questions.\newline
\end{remarka}

\section{Proof of the main results}

\label{sec3}Let us give the general plan for the proof of the main results.
\textit{The main difficulty is to obtain a priori estimates assuring that the
density is bounded and bounded by below}. Once this is achieved, one may
follow the approach of D. Hoff in order to obtain the estimates necessary to
prove existence and uniqueness, see \cite{Hof87}, \cite{Hof98} for the
original approach or our more recent contributions \cite{BurHas2020},
\cite{BurHas2021}. Thus, we will not insist on these, by now well-understood points.

First, we recall the finite-time existence results of strong solutions and
some explosion criterion stating that the only way in which a classical
solution might blow-up is because of the appearance of vacuum regions or
because the density does not remain bounded, more precisely, the $L^{\infty}%
$-norm of $\frac{1}{\rho}$ or of $\rho$ blows up. These results were stated
and proved in \cite{BurHas2020} and \cite{BurHas2021} (see Theorem 3.1 from
these papers).

We consider initial data as in Theorem \ref{theo3} and Theorem \ref{theo1} and
we regularize it with a family of molifiers such as to fit in the scope of the
local existence, see Theorem \ref{Cons} below. For each mollified initial data
we consider the regular solution defined on a maximal time of existence. The
regularity of the solutions ensures that all the computations presented below
are justified. In particular, the regularity of the local-in-time solution
justifies the passage from the Eulerian formulation $\left(
\text{\ref{NSK_intro}}\right)  $ to the so-called mass-Lagrangian formulation,
see Section \ref{mass_lag_section}. The fact that this later formulation is
more adapted to obtain apriori estimates is known since the work of Kazhikhov
and Shelukhin \cite{KS77} and was used more recently by Germain and Le Floch
in \cite{Germain}.

A first important step is to show that the density \textit{is uniformly
bounded using only the uniform bound on the energy and the $L^{\infty}$ norm
of $\rho_{0}$}. This is the objective of Section
\ref{Upper_bound_for_density_section} . More precisely, we show that the local
solution (or any regular enough solution) for the NSK system verifies
\[
\rho\left(  t,x\right)  \leq C\left(  t,E_{c},\left\Vert \rho_{0}\right\Vert
_{L^{\infty}}\right)  \;\;\mbox{for all}\;x\in\mathbb{R}
\]
where $C$ is continuous and increasing with respect to $t$. This estimate is
crucial if one wants to consider the vanishing capillarity limit. The proof is
rather technical and uses a lot of tricks inspired by earlier works of the
Russian school of PDEs, see \cite{Ka},\cite{KS77}. We refer the reader to the
monographic \cite{AntKazMon} for a systematic treatment of the $1d$ NS system
with various boundary conditions. Loosely speaking, the idea is to work with a
primitive of the function defining the velocity's equation. The fact that we
work with functions defined on the whole real line adds another difficulty to
the proof because we need to localize the arguments. We point out that
Proposition \ref{R_auxilary_lemma_1.2} below turns out to be crucial in order
to carry out our proof.

In Sections \ref{bounds_effective_section} and
\ref{lower_bound_for_density_section} we show that the density part of such a
regular solution remains bounded and bounded below. In particular, this will
ensure that its maximal time of existence is $+\infty$. In the same way as
working with a primitive for the momentum equation leads to an upper bound for
the density under not to restrictive hypothesis, \textit{working with a
primitive of the equation of the effective velocity will enable us to obtain a
lower bound for the density with minimal regularity assumptions}. More
precisely, in the NSK case we show that the density is lower bounded by a
constant that depends continuously on time, initial energy, $\left\Vert
\rho_{0}\right\Vert _{L^{\infty}}$, $\left\Vert \frac{1}{\rho_{0}}\right\Vert
_{L^{\infty}}$ and on $M_{0}\in\mathbb{R}$ defined as the smallest constant
such that
\begin{equation}
v_{1|t=0}=u_{0}+r_{1}\partial_{x}\varphi\left(  \rho_{0}\right)  \leq
M_{0},\text{ a.e. in }\mathbb{R}. \label{plan1}%
\end{equation}
Relation $\left(  \text{\ref{plan1}}\right)  $ translates the fact the the
positive part of the effective velocity should be bounded but does not offer
any information on the negative part. In particular, this explains why
densities having negative jumps enter the framework of $\left(
\text{\ref{plan1}}\right)  $. In the case of the NS system, the information we
obtain is better. Indeed, we can say that the effective velocity is controlled
pointwise by the initial effective velocity plus a term that we control. This
allows to show that $\rho$ is lower bounded without assuming $\left(
\text{\ref{plan1}}\right)  $.

Finally, in Section \ref{Hoff_section}, armed with the uniform estimates for
the density, we show how to carry-out the Hoff-type estimates program. With
respect to our recent contribution \cite{BurHas2021}, the Hoff-estimates for
the effective velocity that degenerates to $u$, are shown to hold true
uniformly with respect to the capillarity coefficient $c$ which is important
for the proof of Theorem \ref{theo2}.

\subsection{Existence of strong solution in finite time}

\label{subsec1} In order to prove the existence of global strong solutions for
the NS and NSK systems, we start with recalling the following result.

\begin{theorem}
\label{Cons} Assume that the viscosity coefficients satisfy the assumptions
(\ref{H1})-(\ref{H4}), $s\geq3$ and $(\rho_{0}-1,u_{0})\in H^{s+1}\times
H^{s}(\mathbb{R})$ with $\frac{1}{\rho_{0}} \in L^{\infty}(\mathbb{R})$. Then
there exists $T^{*}>0$ such that there exists a strong solution $(\rho,u)$ of
the system (\ref{NSK_intro}) on $(0,T^{*})$ with $\forall T\in(0,T^{*})$:
\[
(\rho-1)\in C(0,T,H^{s+1}(\mathbb{R}))\cap L^{2}(0,T,H^{s+2}(\mathbb{R}%
)),\,u\in C(0,T,H^{s}(\mathbb{R}))\cap L^{2}(0,T,H^{s+1}(\mathbb{R})),
\]
and for all $t\in(0,T^{*})$:
\[
\|\frac{1}{\rho}(t,\cdot)\|_{L^{\infty}}\leq C(t),
\]
where $C(t)<+\infty$ if $t\in(0,T^{*})$. In addition, if:
\[
\sup_{t\in(0,T^{*})}[\|\frac{1}{\rho}(t,\cdot)\|_{L^{\infty}}+\|\rho
(t,\cdot)\|_{L^{\infty}}]\leq C<+\infty,
\]
then the solution can be continued beyond $(0,T^{*})$. \label{theo4}
\end{theorem}

The above result claims that the only way a strong solution might blow-up in
finite time $T^{\star}$ is if the $L^{\infty}$-norm of $\frac{1}{\rho}$ or
$\rho$ blows-up at time $T^{\ast}$. We refer to \cite{Chen1,Chen2} for the
proof of existence of a strong solution in finite time and we refer to
\cite{BurHas2021} for a proof of the blow-up criterion.\newline Since the
initial data in Theorem \ref{theo3} and Theorem \ref{theo1} are less regular
than in the Theorem \ref{Cons}, we cannot directly use Theorem \ref{Cons}.
Consider $\left(  \rho_{0},u_{0}\right)  \in L^{\infty}\left(  \mathbb{R}%
\right)  \times L^{2}\left(  \mathbb{R}\right)  $ such that $E_{c}\left(
\rho_{0},u_{0}\right)  <\infty$ and, in the case of the NSK system,
\[
u_{0}+r_{1}\left(  c\right)  \partial_{x}\varphi\left(  \rho_{0}\right)  \leq
M_{0}\text{ a.e. on }\mathbb{R}\text{.}%
\]

We regularize the initial data as follows:%
\begin{equation}
\left\{
\begin{array}
[c]{l}%
\rho_{0}^{n}=j_{n}\ast\rho_{0},\\
v_{1|t=0}^{n}=j_{n}\ast\left(  v_{1|t=0}\right)  =j_{n}\ast\left(  u_{0}%
+r_{1}\left(  c\right)  \partial_{x}\varphi\left(  \rho_{0}\right)  \right)
,\\
u_{0}^{n}=v_{1|t=0}^{n}-r_{1}\left(  c\right)  \partial_{x}\varphi\left(
\rho_{0}^{n}\right)  .
\end{array}
\right.  \label{regularization}%
\end{equation}
with $j_{n}$ a family of mollifiers, $j_{n}(y)=nj(ny)$ with $j\in C^{\infty
}(\mathbb{R})$ such that%
\[
0\leq j\leq1,%
{\displaystyle\int_{\mathbb{R}}}
j(y)dy=1\text{ and }\operatorname*{Supp}j\subset\lbrack-2,2].
\]
We deduce that $(\rho_{0}^{n}-1,v_{1|t=0}^{n})$ belong to all Sobolev spaces
$H^{s}(\mathbb{R})$ with $s\geq0$. Furthermore, by the composition theorem we
can prove that $\varphi(\rho_{0}^{n})-\varphi(1)$ belongs to $H^{k}%
(\mathbb{R})$ for any $k\geq0$ and consequently we obtain that $u_{0}^{n}\in
H^{k}(\mathbb{R})$ for $k\geq3$. However, the higher order Sobolev norms
explode when $n\rightarrow\infty$. The only informations that are uniform in
$n\in\mathbb{N}$, the regularization parameter, are the following:
\begin{equation}
\left\{
\begin{array}
[c]{l}%
0<\left\Vert \frac{1}{\rho_{0}}\right\Vert _{L^{\infty}}\leq\rho_{0}^{n}%
\leq\left\Vert \rho_{0}\right\Vert _{L^{\infty}}<+\infty,\\
E_{c}\left(  \rho_{0}^{n},u_{0}^{n}\right)  \leq2E_{c}\left(  \rho_{0}%
,u_{0}\right)  ,\\
v_{1|t=0}^{n}=u_{0}^{n}+r_{1}\left(  c\right)  \partial_{x}\varphi\left(
\rho_{0}^{n}\right)  \leq M_{0}.
\end{array}
\right.  \label{goo}%
\end{equation}

We can now apply the Theorem \ref{Cons} to the sequence of initial data
$(\rho_{0}^{n},u_{0}^{n})_{n\geq0}$ which provides us the existence of a
strong solution $(\rho^{n},u^{n})$ of the system $\left(
\text{\ref{NSK_intro}}\right)  $ on some finite time interval $(0,T_{n})$ with
$T_{n}>0$. We denote by $v_{0}^{n}=u^{n}+r_{0}\left(  c\right)  \partial
_{x}\varphi(\rho^{n})$ and by $v_{1}^{n}=u^{n}+r_{1}\left(  c\right)
\partial_{x}\varphi(\rho^{n})$ the effective velocities introduced in Section
\ref{algebrai_structure_Section}. Using $\left(  \text{\ref{goo}}\right)  $
and the regularity of the constructed solution we see that for any
$t\in(0,T_{n})$ and sufficiently large $n\in\mathbb{N}$ we have that%
\begin{equation}
\int_{\mathbb{R}}\left(  \frac{1}{2}\rho^{n}\left\vert u^{n}\right\vert
^{2}(t,x)+\rho^{n}e\left(  \rho^{n}\right)  (t,x)+c(\partial_{x}G\left(
\rho^{n}\right)  )^{2}\right)  dx+\int_{0}^{t}\int_{\mathbb{R}}\mu(\rho
^{n})(\partial_{x}u^{n})^{2}(s,x)dsdx\leq2E_{c}\left(  \rho_{0},u_{0}\right)
<+\infty\label{basic_energy8}%
\end{equation}%
\begin{equation}
\int_{\mathbb{R}}\left(  \frac{1}{2}\rho^{n}\left\vert v_{0}^{n}\right\vert
^{2}(t,x)+\rho^{n}e\left(  \rho^{n}\right)  (t,x)\right)  dx+(1-r_{0}%
(c))\int_{0}^{t}\int_{\mathbb{R}}\mu(\rho^{n})(\partial_{x}v_{0}^{n}%
)^{2}(s,x)dsdx\leq2E_{c}\left(  \rho_{0},u_{0}\right)  <+\infty.
\label{basic_energy7}%
\end{equation}

We anticipate that at the end of Section \ref{lower_bound_for_density_section}
when we will have shown that $\frac{1}{\rho_{n}}$ is bounded, we will be able
to conclude that for any $n\in\mathbb{N}$ we have $T_{n}=+\infty$, see
Proposition \ref{rho_is_bounded_by_bellow}.

The end of the proof of Theorem \ref{theo1} and \ref{theo3} will consist in
proving that $(\rho^{n},u^{n})_{n\in\mathbb{N}}$ converges up to a subsequence
to a global strong solution $(\rho,u)$ of the system $\left(
\text{\ref{NSK_intro}}\right)  $. In the sequel for simplicity we omit the
subscript $n\in\mathbb{N}$ and all the estimate are obtained on the time
interval $(0,T_{n})$.

\subsection{The Mass-Lagrangian formulation\label{mass_lag_section}}

It is well known that if
\[
\inf_{x\in\mathbb{R}}\rho_{0}\left(  x\right)  >0\text{ for all }%
x\in\mathbb{R}\text{,}%
\]
then the change of variables:
\[
\widetilde{\rho}(t,m)=\rho(t,X(t,Y(m)))\;\;\mbox{and}\;\;\widetilde
{u}(t,m)=u(t,X(t,Y(m)))
\]
with:
\[
\begin{aligned}
&X(t,x)=x+\int^t_0u(s,X(s,x) )ds\;\;\mbox{and}\;\;Y^{-1}(m)=\int^m_0\rho_0(z)dz.
\end{aligned}
\]
where $Y^{-1}$ is the inverse of $Y$ (we observe in particular that $Y^{-1}$
is strictly increasing on $\mathbb{R}$ since $\inf_{x\in\mathbb{R}}\rho
_{0}\left(  x\right)  >0$) transforms the system $\left(
\text{\ref{NSK_intro}}\right)  $ and $\left(  \text{\ref{equations_for_v}%
}\right)  $ into:
\begin{equation}
\left\{
\begin{array}
[c]{l}%
\partial_{t}\widetilde{\rho}+\widetilde{\rho}^{2}\partial_{m}\widetilde
{u}=0,\\
\partial_{t}\widetilde{u}-\partial_{m}\left(  \widetilde{\rho}\mu\left(
\widetilde{\rho}\right)  \partial_{m}\widetilde{u}\right)  +\partial
_{m}\widetilde{\rho}^{\gamma}=c\partial_{m}(\widetilde{\rho}\mu\left(
\widetilde{\rho}\right)  \partial_{mm}^{2}\psi\left(  \widetilde{\rho}\right)
),\\
\partial_{t}\widetilde{v}_{i}-\left(  1-r_{i}\right)  \partial_{m}\left(
\widetilde{\rho}\mu\left(  \widetilde{\rho}\right)  \partial_{m}\widetilde
{v}_{i}\right)  +\partial_{m}\widetilde{\rho}^{\gamma}=0,
\end{array}
\right.  \label{NSK_mass_Lagrangian}%
\end{equation}
where we denote%
\begin{equation}
\psi^{\prime}\left(  \rho\right)  =\rho\varphi^{\prime}\left(  \rho\right)
=\frac{\mu\left(  \rho\right)  }{\rho}, \label{def_psi}%
\end{equation}
and with $m\in\mathbb{R}$. We note that we have taken $a=1$ in the pressure
term only for simplifying the notations. The conditions $\left(
\text{\ref{conditions_at_infinity}}\right)  $ yield%
\[
\widetilde{\rho}\left(  t,m\right)  \rightarrow1\text{ and }\widetilde
{u}\left(  t,m\right)  \rightarrow0\text{ when }\left\vert m\right\vert
\rightarrow\infty.
\]

Obviously, along as the solutions are regular, for instance the regularity
assumed in Theorem \ref{Cons} is sufficient, then for all $t>0,$ the
application
\[
m\rightarrow X\left(  t,Y\left(  m\right)  \right)
\]
is a $C^{1}$-diffeormorfism from $\mathbb{R}$ to $\mathbb{R}$ with inverse
denoted by
\begin{equation}
x\rightarrow m\left(  t,x\right)  . \label{definitie_invers_iso_lag}%
\end{equation}
Obviously, we have that%
\begin{equation}
\int_{\mathbb{R}}F\left(  \rho\left(  t,x\right)  ,u\left(  t,x\right)
\right)  =\int_{\mathbb{R}}F\left(  \widetilde{\rho}\left(  t,m\right)
,\widetilde{u}\left(  t,m\right)  \right)  \frac{dm}{\widetilde{\rho}\left(
t,m\right)  }, \label{change_of_var_integral_formula}%
\end{equation}
for any function $F$. Also, we have that%
\begin{equation}
\partial_{x}G\left(  t,X\left(  t,Y\left(  m\right)  \right)  \right)
=\widetilde{\rho}\partial_{m}G\left(  t,m\right)
\label{change_of_var_derivative_formula}%
\end{equation}

We will use system $\left(  \text{\ref{NSK_mass_Lagrangian}}\right)  $ in
order to obtain a priori estimates as it is easier to manipulate than $\left(
\text{\ref{NSK_intro}}\right)  $. Also, using $\left(
\text{\ref{change_of_var_derivative_formula}}\right)  $, we see that for
$i\in\{0,1\}$ we have that
\begin{equation}
\widetilde{v}_{i}=\widetilde{\left(  u+r_{i}\left(  c\right)  \varphi^{\prime
}\left(  \rho\right)  \partial_{x}\rho\right)  }=\widetilde{u}+r_{i}\left(
c\right)  \widetilde{\rho}\varphi^{\prime}\left(  \widetilde{\rho}\right)
\partial_{m}\widetilde{\rho}=\widetilde{u}+r_{i}\left(  c\right)  \partial
_{m}\psi\left(  \widetilde{\rho}\right)  . \label{relation_v_i}%
\end{equation}
The mass equation rewrites in terms of the so-called specific volume:%
\begin{equation}
\frac{\partial}{\partial t}\frac{1}{\widetilde{\rho}}=\partial_{m}%
\widetilde{u}. \label{specific_volume_equation}%
\end{equation}

\subsection{Basic energy estimates and some properties of regular solutions}

The natural energy functional associated to system $\left(
\text{\ref{NSK_mass_Lagrangian}}\right)  $ is for any $t>0$
\begin{equation}
\int_{\mathbb{R}}\left(  \dfrac{\widetilde{u}^{2}}{2}+e\left(  \widetilde
{\rho}\right)  +\frac{c}{2}(\partial_{m}\psi\left(  \widetilde{\rho}\right)
)^{2}\right)  (t,x)dx+\int_{0}^{t}\int_{\mathbb{R}}\widetilde{\rho}\mu\left(
\widetilde{\rho}\right)  (\partial_{m}\widetilde{u})^{2}dsdx\leq E_{c}\left(
\rho_{0},u_{0}\right)  , \label{R_energy1}%
\end{equation}
where $E_{c}\left(  \rho_{0},u_{0}\right)  $ is defined in $\left(
\text{\ref{basic_energy}}\right)  $ while $e$ is defined by the relation
$\left(  \text{\ref{potential_energy}}\right)  $. Obviously, the functional in
$\left(  \text{\ref{R_energy1}}\right)  $ is nothing other than the Eulerian
energy transformed via the formulae $\left(
\text{\ref{change_of_var_integral_formula}}\right)  $ and $\left(
\text{\ref{change_of_var_derivative_formula}}\right)  $. Moreover we see that%
\begin{align}
&  \int_{\mathbb{R}}\left(  \dfrac{\widetilde{v}_{0}^{2}}{2}+e(\widetilde
{\rho})\right)  (t,x)dx+\left(  1-r_{0}(c)\right)  \int_{0}^{t}\int
_{\mathbb{R}}\widetilde{\rho}\mu\left(  \widetilde{\rho}\right)  (\partial
_{m}\widetilde{v}_{0})^{2}dsdx\label{R_weak_BD_entropy1}\\
&  \leq\int_{\mathbb{R}}\left(  \frac{1}{2}\widetilde{v}_{0}^{2}%
+e(\widetilde{\rho}_{0})\right)  (x)dx\leq2E_{c}\left(  \rho_{0},u_{0}\right)
.
\end{align}

\begin{remark}
In the rest of the paper, in order to ease the reading we will rather use the
notation $E_{c}$ instead of $E_{c}\left(  \rho_{0},u_{0}\right)  .$
\end{remark}

Using the uniform bounds ensured by the energy inequality $\left(
\text{\ref{R_energy1}}\right)  $, we deduce the following proposition.

\begin{proposition}
\label{R_auxilary_lemma_1}Consider $\ell\in\mathbb{R}$ and $\theta\in(0,1]$.
Then,
\begin{equation}
\left\{
\begin{array}
[c]{l}%
\frac{1}{2}\left(  2\left(  \gamma-1\right)  E_{c}+\gamma\right)
^{-\frac{\theta}{\gamma-1}}\leq%
{\displaystyle\int_{\ell}^{\ell+1}}
\dfrac{dq}{\widetilde{\rho}^{\theta}\left(  t,q\right)  }\leq\left(
E_{c}+\frac{\gamma}{\gamma-1}\right)  ^{\frac{1}{\theta}},\\%
{\displaystyle\int_{\ell}^{\ell+1}}
\widetilde{\rho}^{\gamma-1}\left(  t,q\right)  dq\leq\left(  \gamma-1\right)
E_{c}+\gamma.
\end{array}
\right.  \label{inequality_whole_space}%
\end{equation}

\end{proposition}

\textit{Proof of Proposition} \ref{R_auxilary_lemma_1}:

The upper bound is easy to obtain as we can write that%
\[
\int_{\ell}^{\ell+1}\frac{1}{\widetilde{\rho}}\leq\int_{\ell}^{\ell+1}%
\frac{\widetilde{\rho}^{\gamma-1}}{\gamma-1}+\frac{1}{\widetilde{\rho}}%
=\int_{\ell}^{\ell+1}e\left(  \widetilde{\rho}\right)  +\frac{\gamma}%
{\gamma-1}\leq E_{c}+\frac{\gamma}{\gamma-1}.
\]
Obviously, using Jensen's inequality, we have that%
\[
\int_{\ell}^{\ell+1}\frac{1}{\widetilde{\rho}^{\theta}}\leq\left(  \int_{\ell
}^{\ell+1}\frac{1}{\widetilde{\rho}}\right)  ^{\theta},
\]
from which we deduce the right hand side of the first inequality in $\left(
\text{\ref{inequality_whole_space}}\right)  $. The lower bound is deduced in
the following way: we suppose now that there exists $\varepsilon>0$ such that%
\begin{equation}
\int_{\ell}^{\ell+1}\frac{dq}{\widetilde{\rho}^{\theta}\left(  t,q\right)
}\leq\varepsilon. \label{utile}%
\end{equation}
We set now:
\[
A=\left\{  m\in\left[  \ell,\ell+1\right]  :\widetilde{\rho}\left(
t,m\right)  \geq\frac{1}{\left(  2\varepsilon\right)  ^{\frac{1}{\theta}}%
}\right\}  .
\]
Then using (\ref{utile}) we deduce that:
\[
2\varepsilon|^{c}A|\leq\int_{^{c}A}\frac{dq}{\widetilde{\rho}^{\theta}\left(
t,q\right)  }\leq\int_{\ell}^{\ell+1}\frac{dq}{\widetilde{\rho}^{\theta
}\left(  t,q\right)  }\leq\varepsilon.
\]
It implies that
\[
meas\left\{  m\in\left[  \ell,\ell+1\right]  :\widetilde{\rho}\left(
t,m\right)  \geq\frac{1}{\left(  2\varepsilon\right)  ^{\frac{1}{\theta}}%
}\right\}  \geq\frac{1}{2}.
\]
But then, we have that%
\[
E_{c}\geq\int_{\ell}^{\ell+1}e\left(  \widetilde{\rho}\right)  \geq\frac{1}%
{2}\left\{  \frac{1}{\gamma-1}\left(  \frac{1}{2\varepsilon}\right)
^{\frac{\gamma-1}{\theta}}-\frac{\gamma}{\gamma-1}\right\}  ,
\]
which implies that
\[
\varepsilon\geq\frac{1}{2}\left(  2\left(  \gamma-1\right)  E_{c}%
+\gamma\right)  ^{-\frac{\theta}{\gamma-1}}.
\]
This concludes the proof of Proposition \ref{R_auxilary_lemma_1} for the first
estimate in (\ref{inequality_whole_space}). The second one is a direct
consequence of the energy estimate (\ref{R_energy1}). $\Box$

\begin{proposition}
\label{R_auxilary_lemma_1.2}There exists a constant $C\left(  E_{c}\right)  $
that depends only on the initial energy $E_{c}$ such that the following holds
true. For any $t>0$ and $\ell\in\mathbb{R}$ there exists a point $m\left(
t,\ell\right)  \in\left[  \ell,\ell+1\right]  $ such that%
\[
\sup_{s\in\left[  0,t\right]  }\Xi\left(  \widetilde{\rho}\left(  s,m\left(
t,\ell\right)  \right)  \right)  \leq\sup_{m\in\left[  \ell,\ell+1\right]
}\Xi\left(  \rho_{0}\left(  m\right)  \right)  +C\left(  E_{c}\right)  t.
\]

\end{proposition}

\textit{Proof of Proposition \ref{R_auxilary_lemma_1.2}:} We recall that the
mass equation is%
\[
\partial_{t}\widetilde{\rho}\left(  t,m\right)  =-\widetilde{\rho}^{2}\left(
t,m\right)  \partial_{m}\widetilde{u}%
\]
and as such, we get that%
\[
\Xi\left(  \widetilde{\rho}\left(  s,m\right)  \right)  =\Xi\left(
\widetilde{\rho}_{0}\left(  m\right)  \right)  -\int_{0}^{s}\widetilde{\rho
}^{\frac{\gamma-1}{2}}(\widetilde{\rho}\mu\left(  \widetilde{\rho}\right)
)^{\frac{1}{2}}\partial_{m}\widetilde{u}.
\]
From the last identity we get that%
\[
\sup_{s\in\left[  0,t\right]  }\Xi\left(  \widetilde{\rho}\left(  s,m\right)
\right)  \leq\Xi\left(  \widetilde{\rho}_{0}\left(  m\right)  \right)
+\int_{0}^{t}\widetilde{\rho}^{\frac{\gamma-1}{2}}(\widetilde{\rho}\mu\left(
\widetilde{\rho}\right)  )^{\frac{1}{2}}\left\vert \partial_{m}\widetilde
{u}\right\vert .
\]
Consider $\ell\in\mathbb{R}$ and integrate the above relation and using the
second inequality from Proposition $\text{\ref{R_auxilary_lemma_1}} $ we
obtain%
\begin{align*}
\int_{\ell}^{\ell+1}\sup_{s\in\left[  0,t\right]  }\Xi\left(  \widetilde{\rho
}\left(  s,m\right)  \right)  dm  &  \leq\int_{\ell}^{\ell+1}\Xi\left(
\widetilde{\rho}_{0}\left(  m\right)  \right)  dm+\int_{0}^{t}\int_{\ell
}^{\ell+1}\widetilde{\rho}^{\frac{\gamma-1}{2}}(\widetilde{\rho}\mu\left(
\widetilde{\rho}\right)  )^{\frac{1}{2}}\left\vert \partial_{m}\widetilde
{u}\right\vert \\
&  \leq\int_{\ell}^{\ell+1}\Xi\left(  \widetilde{\rho}_{0}\left(  m\right)
\right)  dm+\int_{0}^{t}\int_{\ell}^{\ell+1}\widetilde{\rho}^{\gamma-1}%
+\int_{0}^{t}\int_{\ell}^{\ell+1}\widetilde{\rho}\mu\left(  \widetilde{\rho
}\right)  \left\vert \partial_{m}\widetilde{u}\right\vert ^{2}\\
&  \leq\sup_{m\in\left[  \ell,\ell+1\right]  }\Xi\left(  \rho_{0}\left(
m\right)  \right)  +t\left(  \left(  \gamma-1\right)  E_{c}+\gamma\right)
+E_{c}.
\end{align*}
From this we infer that $\int_{\ell}^{\ell+1}\sup_{s\in\left[  0,t\right]
}\Xi\left(  \widetilde{\rho}\left(  s,m\right)  \right)  dm$ is controlled by
a linear in time function the coefficients of which are controlled uniformly
with respect to the initial data. It then follows that there is a point
$m\left(  t,\ell\right)  \in\left[  \ell,\ell+1\right]  $ where
\[
\sup_{s\in\left[  0,t\right]  }\Xi\left(  \widetilde{\rho}\left(  s,m\left(
t,\ell\right)  \right)  \right)  \leq\sup_{m\in\left[  \ell,\ell+1\right]
}\Xi\left(  \rho_{0}\left(  m\right)  \right)  +C\left(  E_{c}\right)  t.
\]
This concludes the proof of Proposition \ref{R_auxilary_lemma_1.2}. $\Box$

\begin{corollary}
\label{Corolar_points}Suppose that $\mu\left(  \rho\right)  $ verifies
hypothesis $\left(  \text{\ref{H5}}\right)  $. There exists a constant
$C\left(  E_{c},\left\Vert \rho_{0}\right\Vert _{L^{\infty}},t\right)  $ that
depends only on the initial energy $E_{c}$, $\left\Vert \rho_{0}\right\Vert
_{L^{\infty}}$ and time such that the following holds true. For any $t>0$ and
any $\ell\in\mathbb{R}$ there exists a point $m\left(  t,\ell\right)
\in\left[  \ell,\ell+1\right]  $ such that%
\[
\sup_{s\in\left[  0,t\right]  }\widetilde{\rho}\left(  s,m\left(
t,\ell\right)  \right)  \leq C\left(  E_{c},\left\Vert \rho_{0}\right\Vert
_{L^{\infty}},t\right)  .
\]

\end{corollary}

Obviously, \ref{Corolar_points} is a consequence of Proposition
\ref{R_auxilary_lemma_1.2} and of the hypothesis $\left(  \text{\ref{H5}%
}\right)  $ assumed for the function $\Xi$.

In the next section we show how to use these two basic estimates in order to
obtain upper and lower bounds for the density.

\subsection{Upper bound for the density\label{Upper_bound_for_density_section}%
}

The first step is to recover an upper bound for the density that \textit{does
not degenerate when the capillarity coefficient }$c$\textit{ goes to }$0$. It
transpires that for this purpose we cannot use the so called BD-entropy.

\begin{proposition}
\label{rho_is_bounded_in_Linfty}Assume that $\mu$ verifies the hypothesis
$\left(  \text{\ref{H1}}\right)  $, $\left(  \text{\ref{H2}}\right)  $,
$\left(  \text{\ref{H3}}\right)  $. Consider $\left(  \widetilde{\rho
},\widetilde{u}\right)  $ regular enough solutions for the $\left(
\text{\ref{NSK_mass_Lagrangian}}\right)  $ system verifying the energy
estimate:%
\[
\int_{\mathbb{R}}\left(  \dfrac{\widetilde{u}^{2}}{2}+e\left(  \widetilde
{\rho}\right)  +\frac{c}{2}(\partial_{m}\psi\left(  \widetilde{\rho}\right)
)^{2}\right)  (t,x)dx+\int_{0}^{t}\int_{\mathbb{R}}\widetilde{\rho}\mu\left(
\widetilde{\rho}\right)  (\partial_{m}\widetilde{u})^{2}dsdx\leq E_{c}.
\]
Also, we assume that
\[
\widetilde{\rho}_{|t=0}=\rho_{0}\in L^{\infty}\left(  \mathbb{R}\right)  .
\]

\begin{enumerate}
\item The general case: Assume that $\mu$ also verifies $\left(
\text{\ref{H5}}\right)  $. Then, there exists a continuous function
$C:[0,\infty)\times\lbrack0,\infty)\times\lbrack0,\infty)\rightarrow
\lbrack0,\infty)$ which is increasing w.r.t. to the first variable and such
that
\begin{equation}
\widetilde{\rho}\left(  t,m\right)  \leq C\left(  t,E_{c},\left\Vert \rho
_{0}\right\Vert _{L^{\infty}}\right)  . \label{rho_borne_1}%
\end{equation}

\item The NS system: If $c=0$, i.e. we are considering the Navier-Stokes
system, the same conclusion as in $\left(  \text{\ref{rho_borne_1}}\right)  $
holds true without assuming hypothesis $\left(  \text{\ref{H5}}\right)  $ for
the viscosity function $\mu$.
\end{enumerate}
\end{proposition}

\begin{remark}
The results of Proposition \ref{R_auxilary_lemma_1}, Proposition
\ref{R_auxilary_lemma_1.2} and Corollary \ref{Corolar_points} hold true for
the pair $\left(  \widetilde{\rho},\widetilde{u}\right)  $ of Proposition
\ref{rho_is_bounded_in_Linfty}.
\end{remark}

\begin{remark}
In the case of the Navier-Stokes system for coefficients of the form
$\mu\left(  \rho\right)  =\rho^{\alpha},$ the only restriction that we need
is
\[
0<\alpha\leq\gamma.
\]

\end{remark}

\begin{remark}
\label{remarca_rho}Owing to the fact that
\[
\rho\left(  t,x\right)  =\widetilde{\rho}\left(  t,m\left(  t,x\right)
\right)  ,
\]
with $\left(  \text{\ref{definitie_invers_iso_lag}}\right)  $, we also have
that%
\begin{equation}
\rho\left(  t,x\right)  \leq C\left(  t,E_{c},\left\Vert \rho_{0}\right\Vert
_{L^{\infty}}\right)  . \label{rho_borne2}%
\end{equation}

\end{remark}

\begin{remark}
The proof of the above proposition is rather technical and uses a lot of
tricks inspired by earlier works of the Russian school of PDEs, see
\cite{Ka},\cite{KS77} or the monography \cite{AntKazMon} for a systematic
treatment of the NS system with various boundary conditions. Loosely speaking,
the idea is to work with a primitive of the function defining the velocity's
equation. Taking in account that the viscosity term can be written as the time
derivative of $\psi\left(  \widetilde{\rho}\right)  $, in the case of the NS
system, one obtains a dumped differential equation for $\psi\left(
\widetilde{\rho}\right)  $ while the source terms are integrals of the
unknowns for which it turns out that we can provide estimates using only
$E_{c}$. The same conclusion holds true for the NSK system although things
become more complicated. In this case, $\psi\left(  \widetilde{\rho}\right)  $
verifies a parabolic equation such that we need to use a maximum principle.
The fact that we work with functions defined on the whole real line adds
another difficulty to the proof because we need to localize the arguments. We
point out that Proposition \ref{R_auxilary_lemma_1.2} turns out to be crucial
in order to carry out our proof.
\end{remark}

\textit{Proof of Proposition \ref{rho_is_bounded_in_Linfty}}:

Consider $m,q,\ell$ $\in\mathbb{R}$ fixed arbitrarily. We are going to use the
following identity which is obtained by integrating the mass-Lagrangian
velocity's equation (\ref{NSK_mass_Lagrangian}) between $m$ and $q$ :%

\begin{equation}
\frac{d}{dt}\left\{  \int_{q}^{m}\widetilde{u}\right\}  +\left.  \left(
-\widetilde{\rho}\mu\left(  \widetilde{\rho}\right)  \partial_{m}\widetilde
{u}+\widetilde{\rho}^{\gamma}-c\widetilde{\rho}\mu\left(  \widetilde{\rho
}\right)  \partial_{mm}^{2}\psi\left(  \widetilde{\rho}\right)  \right)
\right\vert _{q}^{m}=0 \label{difference}%
\end{equation}
We put the above relation under the form
\begin{align*}
&  \frac{d}{dt}\left\{  \int_{\ell}^{m}\widetilde{u}\right\}  +\left(
-\widetilde{\rho}\mu\left(  \widetilde{\rho}\right)  \partial_{m}\widetilde
{u}+\widetilde{\rho}^{\gamma}-c\widetilde{\rho}\mu\left(  \widetilde{\rho
}\right)  \partial_{mm}^{2}\psi\left(  \widetilde{\rho}\right)  \right)
\left(  t,m\right) \\
&  =\frac{d}{dt}\left\{  \int_{\ell}^{q}\widetilde{u}\right\}  +\left(
-\widetilde{\rho}\mu\left(  \widetilde{\rho}\right)  \partial_{m}\widetilde
{u}+\widetilde{\rho}^{\gamma}-c\widetilde{\rho}\mu\left(  \widetilde{\rho
}\right)  \partial_{mm}^{2}\psi\left(  \widetilde{\rho}\right)  \right)
\left(  t,q\right)  ,
\end{align*}
Next we multiply the above equation with $\dfrac{1}{\widetilde{\rho}\left(
t,q\right)  }$ and we integrate between $\ell$ and $\ell+1$ w.r.t. $q$. We
infer that%
\begin{align*}
&  \left(  \int_{\ell}^{\ell+1}\frac{dq}{\widetilde{\rho}\left(  t,q\right)
}\right)  \left\{  \frac{d}{dt}\left\{  \int_{\ell}^{m}\widetilde{u}\right\}
+\left(  - \widetilde{\rho}\mu\left(  \widetilde{\rho}\right)  \partial
_{m}\widetilde{u}+\widetilde{\rho}^{\gamma}-c\widetilde{\rho}\mu\left(
\widetilde{\rho}\right)  \partial_{mm}^{2}\psi\left(  \widetilde{\rho}\right)
\right)  \left(  t,m\right)  \right\} \\
&  =\int_{\ell}^{\ell+1}\left\{  \frac{1}{\widetilde{\rho}\left(  t,q\right)
}\frac{d}{dt}\left\{  \int_{\ell}^{q}\widetilde{u}\right\}  \right\}
dq+\int_{\ell}^{\ell+1}\left(  -\mu\left(  \widetilde{\rho}\right)
\partial_{m}\widetilde{u}+\widetilde{\rho}^{\gamma-1}-c\mu\left(
\widetilde{\rho}\right)  \partial_{mm}^{2}\psi\left(  \widetilde{\rho}\right)
\right)  \left(  t,q\right)  dq.
\end{align*}
Now, using (\ref{specific_volume_equation}) and a simple integration by parts
gives us:
\begin{align*}
\int_{\ell}^{\ell+1}\left\{  \frac{1}{\widetilde{\rho}\left(  t,q\right)
}\frac{d}{dt}\left\{  \int_{\ell}^{q}\widetilde{u}\right\}  \right\}  dq  &
=\frac{d}{dt}\int_{\ell}^{\ell+1}\left(  \frac{1}{\widetilde{\rho}\left(
t,q\right)  }\int_{\ell}^{q}\widetilde{u}\right)  dq-\int_{\ell}^{\ell
+1}\left(  \frac{d}{dt}\frac{1}{\widetilde{\rho}\left(  t,q\right)  }%
\int_{\ell}^{q}\widetilde{u}\right)  dq\\
&  =\frac{d}{dt}\int_{\ell}^{\ell+1}\left(  \frac{1}{\widetilde{\rho}\left(
t,q\right)  }\int_{\ell}^{q}\widetilde{u}\right)  dq-\int_{\ell}^{\ell
+1}\left(  \partial_{m}\widetilde{u}\left(  t,q\right)  \int_{\ell}%
^{q}\widetilde{u}\right)  dq\\
&  =\frac{d}{dt}\int_{\ell}^{\ell+1}\left(  \frac{1}{\widetilde{\rho}\left(
t,q\right)  }\int_{\ell}^{q}\widetilde{u}\right)  dq-\widetilde{u}\left(
t,\ell+1\right)  \int_{\ell}^{\ell+1}\widetilde{u}\left(  t,q\right)
dq+\int_{\ell}^{\ell+1}\widetilde{u}^{2}\left(  t,q\right)  dq.
\end{align*}
The above identity is rearranged to give%
\begin{align*}
&  \left(  \int_{\ell}^{\ell+1}\frac{dq}{\widetilde{\rho}\left(  t,q\right)
}\right)  \left\{  \frac{d}{dt}\left\{  \int_{\ell}^{m}\widetilde{u}\right\}
+\left(  -\widetilde{\rho}\mu\left(  \widetilde{\rho}\right)  \partial
_{m}\widetilde{u}+\widetilde{\rho}^{\gamma}-c\widetilde{\rho}\mu\left(
\widetilde{\rho}\right)  \partial_{mm}^{2}\psi\left(  \widetilde{\rho}\right)
\right)  \left(  t,m\right)  \right\} \\
&  =\frac{d}{dt}\int_{\ell}^{\ell+1}\left(  \frac{1}{\widetilde{\rho}\left(
t,q\right)  }\int_{\ell}^{q}\widetilde{u}\right)  dq-\widetilde{u}\left(
t,\ell+1\right)  \int_{\ell}^{\ell+1}\widetilde{u}\left(  t,q\right)
dq+\int_{\ell}^{\ell+1}\widetilde{u}^{2}\left(  t,q\right)  dq\\
&  +\int_{\ell}^{\ell+1}\left(  -\mu\left(  \widetilde{\rho}\right)
\partial_{m}\widetilde{u}+\widetilde{\rho}^{\gamma-1}-c\mu\left(
\widetilde{\rho}\right)  \partial_{mm}^{2}\psi\left(  \widetilde{\rho}\right)
\right)  \left(  t,q\right)  dq.
\end{align*}
Recalling that $c=r_{1}(1-r_{1})$, and using that%
\begin{align*}
-\widetilde{\rho}\mu\left(  \widetilde{\rho}\right)  \partial_{m}\widetilde
{u}-c\widetilde{\rho}\mu\left(  \widetilde{\rho}\right)  \partial_{mm}^{2}%
\psi\left(  \widetilde{\rho}\right)   &  =-r_{1}\widetilde{\rho}\mu\left(
\widetilde{\rho}\right)  \partial_{m}\widetilde{u}-\left(  1-r_{1}\right)
\widetilde{\rho}\mu\left(  \widetilde{\rho}\right)  \partial_{m}\widetilde
{u}-c\widetilde{\rho}\mu\left(  \widetilde{\rho}\right)  \partial_{mm}^{2}%
\psi\left(  \widetilde{\rho}\right) \\
&  =\frac{d}{dt}r_{1}\psi\left(  \widetilde{\rho}\right)  -\left(
1-r_{1}\right)  \widetilde{\rho}\mu\left(  \widetilde{\rho}\right)
\partial_{m}\left(  \widetilde{u}+r_{1}\partial_{m}\psi\left(  \widetilde
{\rho}\right)  \right) \\
&  =\frac{d}{dt}r_{1}\psi\left(  \widetilde{\rho}\right)  -\left(
1-r_{1}\right)  \widetilde{\rho}\mu\left(  \widetilde{\rho}\right)
\partial_{mm}\left(  \int_{\ell}^{m}\widetilde{u}+r_{1}\psi\left(
\widetilde{\rho}\right)  \right)  ,
\end{align*}
the left hand side of the above identity is arranged as to obtain%
\begin{align}
&  \frac{d}{dt}\left\{  \int_{\ell}^{m}\widetilde{u}+r_{1}\psi(\widetilde
{\rho})\left(  t,m\right)  \right\}  -\left(  1-r_{1}\right)  \widetilde{\rho
}\mu\left(  \widetilde{\rho}\right)  \partial_{mm}^{2}\left\{  \int_{\ell}%
^{m}\widetilde{u}+r_{1}\psi(\widetilde{\rho})\right\}  +\widetilde{\rho
}^{\gamma}\left(  t,m\right) \nonumber\\
&  =\left(  \int_{\ell}^{\ell+1}\frac{dq}{\widetilde{\rho}\left(  t,q\right)
}\right)  ^{-1}\left\{  \frac{d}{dt}\int_{\ell}^{\ell+1}\left(  \frac
{1}{\widetilde{\rho}\left(  t,q\right)  }\int_{\ell}^{q}\widetilde{u}\right)
dq-\widetilde{u}\left(  t,\ell+1\right)  \int_{\ell}^{\ell+1}\widetilde
{u}\left(  t,q\right)  dq\right. \nonumber\\
&  \left.  \text{
\ \ \ \ \ \ \ \ \ \ \ \ \ \ \ \ \ \ \ \ \ \ \ \ \ \ \ \ \ \ \ \ \ }+\int
_{\ell}^{\ell+1}\left(  \widetilde{u}^{2}-\mu\left(  \widetilde{\rho}\right)
\partial_{m}\widetilde{u}+\widetilde{\rho}^{\gamma-1}-c\mu\left(
\widetilde{\rho}\right)  \partial_{mm}^{2}\psi\left(  \widetilde{\rho}\right)
\right)  \left(  t,q\right)  dq\right\}
\label{R_principe_max_upper_density_1}%
\end{align}
From (\ref{specific_volume_equation}), we get
\begin{align}
&  \left(  \int_{\ell}^{\ell+1}\frac{dq}{\widetilde{\rho}\left(  t,q\right)
}\right)  ^{-1}\frac{d}{dt}\int_{\ell}^{\ell+1}\left(  \frac{1}{\widetilde
{\rho}\left(  t,q\right)  }\int_{\ell}^{q}\widetilde{u}\right)  dq=\frac
{d}{dt}\left\{  \left(  \int_{\ell}^{\ell+1}\frac{dq}{\widetilde{\rho}\left(
t,q\right)  }\right)  ^{-1}\int_{\ell}^{\ell+1}\left(  \frac{1}{\widetilde
{\rho}\left(  t,q\right)  }\int_{\ell}^{q}\widetilde{u}\right)  dq\right\}
\nonumber\\
&  +\int_{\ell}^{\ell+1}\left(  \frac{1}{\widetilde{\rho}\left(  t,q\right)
}\int_{\ell}^{q}\widetilde{u}\right)  dq\left(  \int_{\ell}^{\ell+1}\frac
{dq}{\widetilde{\rho}\left(  t,q\right)  }\right)  ^{-2}\left(  \int_{\ell
}^{\ell+1}\frac{d}{dt}\frac{dq}{\widetilde{\rho}\left(  t,q\right)  }\right)
\nonumber\\
&  =\frac{d}{dt}\left\{  \left(  \int_{\ell}^{\ell+1}\frac{dq}{\widetilde
{\rho}\left(  t,q\right)  }\right)  ^{-1}\int_{\ell}^{\ell+1}\left(  \frac
{1}{\widetilde{\rho}\left(  t,q\right)  }\int_{\ell}^{q}\widetilde{u}\right)
dq\right\}  +\int_{\ell}^{\ell+1}\left(  \frac{1}{\widetilde{\rho}\left(
t,q\right)  }\int_{\ell}^{q}\widetilde{u}\right)  dq\left(  \int_{\ell}%
^{\ell+1}\frac{dq}{\widetilde{\rho}\left(  t,q\right)  }\right)  ^{-2}\left(
\widetilde{u}\left(  t,\ell+1\right)  -\widetilde{u}\left(  t,\ell\right)
\right)  \label{R_principe_max_upper_density_2}%
\end{align}
Owing to $\left(  \text{\ref{R_principe_max_upper_density_1}}\right)  $ and
$\left(  \text{\ref{R_principe_max_upper_density_2}}\right)  $ we get that
\begin{align*}
&  \frac{d}{dt}\left\{  \int_{\ell}^{m}\widetilde{u}+r_{1}\psi(\widetilde
{\rho})\left(  t,m\right)  \right\}  -\left(  1-r_{1}\right)  \widetilde{\rho
}\mu\left(  \widetilde{\rho}\right)  \partial_{mm}^{2}\left\{  \int_{\ell}%
^{m}\widetilde{u}+r_{1}\psi(\widetilde{\rho})\left(  t,m\right)  \right\}
+\widetilde{\rho}^{\gamma}\left(  t,m\right) \\
&  =\frac{d}{dt}\left\{  \left(  \int_{\ell}^{\ell+1}\frac{dq}{\widetilde
{\rho}\left(  t,q\right)  }\right)  ^{-1}\int_{\ell}^{\ell+1}\left(  \frac
{1}{\widetilde{\rho}\left(  t,q\right)  }\int_{\ell}^{q}\widetilde{u}\right)
dq\right\} \\
&  +\int_{\ell}^{\ell+1}\left(  \frac{1}{\widetilde{\rho}\left(  t,q\right)
}\int_{\ell}^{q}\widetilde{u}\right)  dq\left(  \int_{\ell}^{\ell+1}\frac
{dq}{\widetilde{\rho}\left(  t,q\right)  }\right)  ^{-2}\left(  \widetilde
{u}\left(  t,\ell+1\right)  -\widetilde{u}\left(  t,\ell\right)  \right) \\
&  +\left(  \int_{\ell}^{\ell+1}\frac{dq}{\widetilde{\rho}\left(  t,q\right)
}\right)  ^{-1}\text{\ }\left\{  -\widetilde{u}\left(  t,\ell+1\right)
\int_{\ell}^{q}\widetilde{u}\left(  t,q^{\prime}\right)  dq^{\prime}%
+\int_{\ell}^{\ell+1}\left(  \widetilde{u}^{2}-\mu(\widetilde{\rho}%
)\partial_{m}\widetilde{u}+\widetilde{\rho}^{\gamma-1}-c\mu(\widetilde{\rho
})\partial_{mm}^{2}\psi(\widetilde{\rho})\right)  \left(  t,q\right)
dq\right\}
\end{align*}
Now, fix $N\in\mathbb{Z}$ and recalling that the above identity holds for all
$\ell\in\mathbb{R}$, we integrate w.r.t. $\ell$ on $\left[  N,N+1\right]  $ in
order to obtain that%
\begin{align}
&  \frac{d}{dt}\left\{  \int_{N}^{N+1}\left(  \int_{\ell}^{m}\widetilde
{u}\right)  d\ell+r_{1}\psi(\widetilde{\rho})\left(  t,m\right)  \right\}
-\left(  1-r_{1}\right)  \widetilde{\rho}\mu\left(  \widetilde{\rho}\right)
\partial_{mm}^{2}\left\{  \int_{N}^{N+1}\left(  \int_{\ell}^{m}\widetilde
{u}\right)  d\ell+r_{1}\psi(\widetilde{\rho})\left(  t,m\right)  \right\}
+\widetilde{\rho}^{\gamma}\left(  t,m\right) \nonumber\\
&  =\frac{d}{dt}\left\{  \int_{N}^{N+1}\left(  \int_{\ell}^{\ell+1}\frac
{dq}{\widetilde{\rho}\left(  t,q\right)  }\right)  ^{-1}\big(\int_{\ell}%
^{\ell+1}\left(  \frac{1}{\widetilde{\rho}\left(  t,q\right)  }\int_{\ell}%
^{q}\widetilde{u}\right)  dq \big)d\ell\right\} \nonumber\\
&  +\int_{N}^{N+1}\left[  \int_{\ell}^{\ell+1}\left(  \frac{1}{\widetilde
{\rho}\left(  t,q\right)  }\int_{\ell}^{q}\widetilde{u}\right)  dq\left(
\int_{\ell}^{\ell+1}\frac{dq}{\widetilde{\rho}\left(  t,q\right)  }\right)
^{-2}\left(  \widetilde{u}\left(  t,\ell+1\right)  -\widetilde{u}\left(
t,\ell\right)  \right)  \right]  d\ell\nonumber\\
&  +\int_{N}^{N+1}\left(  \int_{\ell}^{\ell+1}\frac{dq}{\widetilde{\rho
}\left(  t,q\right)  }\right)  ^{-1}\text{\ }\left\{  -\widetilde{u}\left(
t,\ell+1\right)  \int_{\ell}^{q}\widetilde{u}\left(  t,q^{\prime}\right)
dq^{\prime}+\int_{\ell}^{\ell+1}\left(  \widetilde{u}^{2}-\mu(\widetilde{\rho
})\partial_{m}\widetilde{u}+\widetilde{\rho}^{\gamma-1}-c\mu(\widetilde{\rho
})\partial_{mm}^{2}\psi(\widetilde{\rho})\right)  \left(  t,q\right)
dq\right\}  d\ell\label{Relation_long}%
\end{align}
Observe that the term appearing in the RHS depends only on time. In the
following lines, we analyze the different terms appearing in the right hand
side of the previous inequality and we want to prove that $\left(
\text{\ref{Relation_long}}\right)  $ implies that
\begin{gather*}
\frac{d}{dt}\left\{  \int_{N}^{N+1}\left(  \int_{\ell}^{m}\widetilde{u}%
d\ell\right)  +r_{1}\psi(\widetilde{\rho})\left(  t,m\right)  \right\}
-\left(  1-r_{1}\right)  \widetilde{\rho}\mu\left(  \widetilde{\rho}\right)
\partial_{mm}^{2}\left\{  \int_{N}^{N+1}\left(  \int_{\ell}^{m}\widetilde
{u}d\ell\right)  +r_{1}\psi(\widetilde{\rho})\left(  t,m\right)  \right\} \\
+\widetilde{\rho}^{\gamma}\left(  t,m\right)  \leq\frac{d\Gamma}{dt}\left(
t\right)  ,
\end{gather*}
where $\Gamma$ is such that%
\begin{equation}
\Gamma\left(  t\right)  \leq C\left(  E_{c}\right)  \left(  1+t\right)  ,
\label{estimGa}%
\end{equation}
for some constant that depends only on $E_{c}$. Let us write that%
\[
c\mu(\widetilde{\rho})\partial_{mm}^{2}\psi(\widetilde{\rho})=r_{1}%
\mu(\widetilde{\rho})r_{0}\partial_{mm}^{2}\psi(\widetilde{\rho})=-r_{1}%
\mu(\widetilde{\rho})\partial_{m}\widetilde{u}+r_{1}\mu(\widetilde{\rho
})\partial_{m}\widetilde{v}_{0}%
\]
such that we get using Proposition \ref{R_auxilary_lemma_1}%
\begin{align}
&  \int_{N}^{N+1}\int_{\ell}^{\ell+1}\left(  \widetilde{u}^{2}-\mu
(\widetilde{\rho})\partial_{m}\widetilde{u}+\widetilde{\rho}^{\gamma-1}%
-c\mu(\widetilde{\rho})\partial_{mm}^{2}\psi(\widetilde{\rho})\right)  \left(
t,q\right)  dqd\ell\nonumber\\
&  =\int_{N}^{N+1}\int_{\ell}^{\ell+1}\left(  \widetilde{u}^{2}+\widetilde
{\rho}^{\gamma-1}\right)  \left(  t,q\right)  dq-\left(  1-r_{1}\right)
\int_{N}^{N+1}\int_{\ell}^{\ell+1}\mu(\widetilde{\rho})\partial_{m}%
\widetilde{u}d\ell-r_{1}\int_{N}^{N+1}\int_{\ell}^{\ell+1}\mu(\widetilde{\rho
})\partial_{m}\widetilde{v}_{0}d\ell\nonumber\\
&  \leq(\gamma+1) E_{c}+\gamma+\left(  1-r_{1}\right)  \int_{N}^{N+1}%
\int_{\ell}^{\ell+1}\left\vert \mu(\widetilde{\rho})\partial_{m}\widetilde
{u}\right\vert d\ell+r_{1}\int_{N}^{N+1}\int_{\ell}^{\ell+1}\left\vert
\mu(\widetilde{\rho})\partial_{m}\widetilde{v}_{0}\right\vert d\ell.
\label{inegalitati_termeni_sursa}%
\end{align}
Using Proposition \ref{R_auxilary_lemma_1} and hypothesis $\left(
\text{\ref{H3}}\right)  $ we get that
\begin{align*}
\left(  \int_{\ell}^{\ell+1}\frac{dq}{\widetilde{\rho}\left(  t,q\right)
}\right)  ^{-1}\int_{\ell}^{\ell+1}\left\vert \mu(\widetilde{\rho}%
)\partial_{m}\widetilde{u}\right\vert  &  \leq\left(  \int_{\ell}^{\ell
+1}\frac{dq}{\widetilde{\rho}\left(  t,q\right)  }\right)  ^{-1}\left(
\int_{\ell}^{\ell+1}\widetilde{\rho}^{-1}\mu\left(  \widetilde{\rho}\right)
\right)  ^{\frac{1}{2}}\left(  \int_{\ell}^{\ell+1}\widetilde{\rho}%
\mu(\widetilde{\rho})(\partial_{m}\widetilde{u})^{2}\right)  ^{\frac{1}{2}}\\
&  \leq C\left(  E_{c},\gamma\right)  +\int_{\ell}^{\ell+1}\widetilde{\rho}%
\mu(\widetilde{\rho})(\partial_{m}\widetilde{u})^{2}.
\end{align*}
Thus, we get that%
\[
\int_{N}^{N+1}\left(  \int_{\ell}^{\ell+1}\frac{dq}{\widetilde{\rho}\left(
t,q\right)  }\right)  ^{-1}\int_{\ell}^{\ell+1}\left\vert \mu(\widetilde{\rho
})\partial_{m}\widetilde{u}\right\vert d\ell\leq C\left(  E_{c},\gamma\right)
+\int_{N}^{N+1}\left(  \int_{\ell}^{\ell+1}\widetilde{\rho}\mu(\widetilde
{\rho})(\partial_{m}\widetilde{u})^{2}\right)  d\ell
\]
We may treat in the same manner the last term appearing in $\left(
\text{\ref{inegalitati_termeni_sursa}}\right)  $. Thus we get that
\begin{gather}
\int_{N}^{N+1}\left\{  \left(  \int_{\ell}^{\ell+1}\frac{dq}{\widetilde{\rho
}\left(  t,q\right)  }\right)  ^{-1}\int_{\ell}^{\ell+1}\left(  \widetilde
{u}^{2}-\mu(\widetilde{\rho})\partial_{m}\widetilde{u}+\widetilde{\rho
}^{\gamma-1}-c\mu(\widetilde{\rho})\partial_{mm}^{2}\psi(\widetilde{\rho
})\right)  \left(  t,q\right)  dq\right\}  d\ell\nonumber\\
\leq C\left(  E_{c},\gamma\right)  +\int_{N}^{N+1}\left(  \int_{\ell}^{\ell
+1}\widetilde{\rho}\mu(\widetilde{\rho})(\partial_{m}\widetilde{u}%
)^{2}\right)  d\ell+\int_{N}^{N+1}\left(  \int_{\ell}^{\ell+1}\widetilde{\rho
}\mu(\widetilde{\rho})(\partial_{m}\widetilde{v}_{0})^{2}\right)  d\ell.
\label{R_rho_upper_Term1}%
\end{gather}

Next, using Proposition \ref{R_auxilary_lemma_1} we get that
\begin{align}
&  \left(  \int_{\ell}^{\ell+1}\frac{dq}{\widetilde{\rho}\left(  t,q\right)
}\right)  ^{-1}\left\{  \widetilde{u}\left(  t,\ell+1\right)  \int_{\ell
}^{\ell+1}\widetilde{u}\left(  t,q^{\prime}\right)  dq^{\prime}\right\}
\nonumber\\
&  \leq4\widetilde{u}\left(  t,\ell+1\right)  E_{c}^{\frac{1}{2}}\left(
2\left(  \gamma-1\right)  E_{c}+\gamma\right)  ^{\frac{1}{\gamma-1}}\leq
C\left(  E_{c},\gamma\right)  \left\vert \widetilde{u}\left(  t,\ell+1\right)
\right\vert , \label{R_rho_upper_Term2}%
\end{align}
respectively%
\begin{align}
&  \int_{\ell}^{\ell+1}\left(  \frac{1}{\widetilde{\rho}\left(  t,q\right)
}\int_{\ell}^{q}\widetilde{u}\right)  dq\left(  \int_{\ell}^{\ell+1}\frac
{dq}{\widetilde{\rho}\left(  t,q\right)  }\right)  ^{-2}\left(  \widetilde
{u}\left(  t,\ell+1\right)  -\widetilde{u}\left(  t,\ell\right)  \right)
\nonumber\\
&  \leq C\left(  E_{c},\gamma\right)  \left(  \left\vert \widetilde{u}\left(
t,\ell+1\right)  \right\vert +\left\vert \widetilde{u}\left(  t,\ell\right)
\right\vert \right)  . \label{R_rho_upper_Term3}%
\end{align}
We thus get that
\begin{align}
&  \int_{N}^{N+1}\left\{  \left(  \int_{\ell}^{\ell+1}\frac{dq}{\widetilde
{\rho}\left(  t,q\right)  }\right)  ^{-1}\left(  - \widetilde{u}\left(
t,\ell+1\right)  \int_{\ell}^{\ell+1}\widetilde{u}\left(  t,q^{\prime}\right)
dq^{\prime}\right)  \right\}  d\ell+\nonumber\\
&  +\int_{N}^{N+1}\left\{  \left[  \int_{\ell}^{\ell+1}\left(  \frac
{1}{\widetilde{\rho}\left(  t,q\right)  }\int_{\ell}^{q}\widetilde{u}\right)
dq\right]  \left(  \int_{\ell}^{\ell+1}\frac{dq}{\widetilde{\rho}\left(
t,q\right)  }\right)  ^{-2}\left(  \widetilde{u}\left(  t,\ell+1\right)
-\widetilde{u}\left(  t,\ell\right)  \right)  \right\}  d\ell\nonumber\\
&  \leq C\left(  E_{c},\gamma\right)  \int_{N}^{N+1}\left(  \left\vert
\widetilde{u}\left(  t,\ell+1\right)  \right\vert +\left\vert \widetilde
{u}\left(  t,\ell\right)  \right\vert \right)  d\ell\leq C\left(  E_{c}%
,\gamma\right)  . \label{R_rho_upper_Term2_and3_integrated}%
\end{align}

Obviously, for all $t\geq0$ we have that%
\begin{equation}
\left(  \int_{\ell}^{\ell+1}\frac{dq}{\widetilde{\rho}\left(  t,q\right)
}\right)  ^{-1}\int_{\ell}^{\ell+1}\left(  \frac{1}{\widetilde{\rho}\left(
t,q\right)  }\int_{\ell}^{q}\widetilde{u}dq\right)  dq\leq C\left(
E_{c},\gamma\right)  . \label{R_rho_upper_Term4}%
\end{equation}
Gathering the estimates $\left(  \text{\ref{R_rho_upper_Term1}}\right)  $,
$\left(  \text{\ref{R_rho_upper_Term2_and3_integrated}}\right)  $, $\left(
\text{\ref{R_rho_upper_Term4}}\right)  $ along with the maximum principle we
see that%

\begin{gather}
\frac{\partial}{\partial t}\left\{  \int_{N}^{N+1}\left(  \int_{\ell}%
^{m}\widetilde{u}d\ell\right)  +r_{1}\psi(\widetilde{\rho})\left(  t,m\right)
\right\}  \text{
\ \ \ \ \ \ \ \ \ \ \ \ \ \ \ \ \ \ \ \ \ \ \ \ \ \ \ \ \ \ \ \ \ \ \ \ \ \ \ \ \ \ \ \ \ \ \ }%
\nonumber\\
\text{ \ \ \ \ \ \ \ \ \ \ \ \ }-\left(  1-r_{1}\right)  \widetilde{\rho}%
\mu\left(  \widetilde{\rho}\right)  \partial_{mm}^{2}\left\{  \int_{N}%
^{N+1}\left(  \int_{\ell}^{m}\widetilde{u}d\ell\right)  +r_{1}\psi
(\widetilde{\rho})\left(  t,m\right)  \right\}  +\widetilde{\rho}^{\gamma
}\left(  t,m\right)  \leq\frac{d\Gamma}{dt}\left(  t\right)
\label{ineg_Gamma}%
\end{gather}
We denote by
\[
\Theta\left(  t,m\right)  =\int_{N}^{N+1}\left(  \int_{\ell}^{m}\widetilde
{u}d\ell\right)  +r_{1}\psi(\widetilde{\rho})\left(  t,m\right)
-\Gamma\left(  t\right)
\]
and we see that $\left(  \text{\ref{ineg_Gamma}}\right)  $ implies that%
\[
\frac{\partial\Theta}{\partial t}\left(  t,m\right)  -\left(  1-r_{1}\right)
\widetilde{\rho}\mu\left(  \widetilde{\rho}\right)  \partial_{mm}^{2}%
\Theta\left(  t,m\right)  \leq0.
\]
Recall the conclusion of Corollary \ref{Corolar_points}, let $m_{1}\left(
t,N\right)  \in\left[  N-1,N\right]  $ and $m_{2}\left(  t,N\right)
\in\left[  N+1,N+2\right]  $ such that
\[
i\in\left\{  1,2\right\}  \forall s\in\left[  0,t\right]  :\rho\left(
s,m_{i}\left(  t,N\right)  \right)  \leq C\left(  t,E_{c},\left\Vert \rho
_{0}\right\Vert _{L^{\infty}}\right)  .
\]
It then follows that%
\begin{equation}
\Theta\left(  s,m_{i}\left(  t,N\right)  \right)  \leq\bar{C}\left(
t,E_{c},\left\Vert \rho_{0}\right\Vert _{L^{\infty}}\right)
\label{control_au_bord}%
\end{equation}
Using a maximum-principle, we aim at showing that $\Theta$ is bounded on the
whole interval $\forall m\in\left[  m_{1}\left(  t,N\right)  ,m_{2}\left(
t,N\right)  \right]  $, namely%
\[
\forall s\in\left[  0,t\right]  :\Theta\left(  s,m\right)  \leq\max\left\{
\sup_{\left[  N-1,N+2\right]  }\Theta\left(  0,m\right)  ,\bar{C}\left(
t,E_{c},\left\Vert \rho_{0}\right\Vert _{L^{\infty}}\right)  \right\}  .
\]
consider
\begin{equation}
A=\left\{  s\in(0,t):\sup_{m\in\left[  m_{1}\left(  t,N\right)  ,m_{2}\left(
t,N\right)  \right]  }\Theta\left(  s,m\right)  >\max\left\{  \sup_{\left[
N-1,N+2\right]  }\Theta\left(  0,m\right)  ,\bar{C}\left(  t,E_{c},\left\Vert
\rho_{0}\right\Vert _{L^{\infty}}\right)  \right\}  \right\}  \label{set}%
\end{equation}
where $\bar{C}\left(  t,E_{c},\left\Vert \rho_{0}\right\Vert _{L^{\infty}%
}\right)  $ is the same constant as in $\left(  \text{\ref{control_au_bord}%
}\right)  $. The function
\[
s\rightarrow\sup_{m\in\left[  m_{1}\left(  t,N\right)  ,m_{2}\left(
t,N\right)  \right]  }\Theta\left(  s,m\right)
\]
is continuous on $[0,t]$ and therefore the set $A$ appearing in $\left(
\text{\ref{set}}\right)  $ is open and therefore it is the union of a at most
countable union of open intervals $I$. Now, there are two ingredients that put
together imply that the set $A$ is empty:

\begin{itemize}
\item In each of the endpoint of such an interval, say $I=\left(
t_{I,g},t_{I,d}\right)  $ we have that
\begin{equation}
\sup_{m\in\left[  m_{1}\left(  t,N\right)  ,m_{2}\left(  t,N\right)  \right]
}\Theta\left(  t_{I,g},m\right)  =\left\{
\begin{array}
[c]{l}%
\sup\limits_{m\in\left[  m_{1}\left(  t,N\right)  ,m_{2}\left(  t,N\right)
\right]  }\Theta\left(  0,m\right)  \text{ if }t_{I,g}=0,\\
\max\left\{  \sup_{\left[  N-1,N+2\right]  }\Theta\left(  0,m\right)  ,\bar
{C}\left(  t,E_{c},\left\Vert \rho_{0}\right\Vert _{L^{\infty}}\right)
\right\}  \text{ if }t_{I,g}>0
\end{array}
\right.  \label{bord_gauche}%
\end{equation}
while
\[
\sup_{m\in\left[  m_{1}\left(  t,N\right)  ,m_{2}\left(  t,N\right)  \right]
}\Theta\left(  t_{I,d},m\right)  =\max\left\{  \sup_{\left[  N-1,N+2\right]
}\Theta\left(  0,m\right)  ,\bar{C}\left(  t,E_{c},\left\Vert \rho
_{0}\right\Vert _{L^{\infty}}\right)  \right\}  .
\]
This fact and (\ref{control_au_bord}) imply that for all $s\in I$ the supremum
over $\left[  m_{1}\left(  t,N\right)  ,m_{2}\left(  t,N\right)  \right]  $ of
$\Theta\left(  s,m\right)  $ is achieved in the interior of the segment
$\left[  m_{1}\left(  t,N\right)  ,m_{2}\left(  t,N\right)  \right]  $ i.e.
\[
\forall s\in I\text{ }\exists m\left(  s\right)  \in(m_{1}\left(  t,N\right)
,m_{2}\left(  t,N\right)  )\text{ such that }\Theta\left(  s,m\left(
s\right)  \right)  =\sup\limits_{m\in\left[  m_{1}\left(  t,N\right)
,m_{2}\left(  t,N\right)  \right]  }\Theta\left(  s,m\right)  .
\]

\item Classical considerations lead to the fact that the function
$s\rightarrow$ $\Theta\left(  s,m\left(  s\right)  \right)  $ is
differentiable almost everywhere on $I$ w.r.t. $s$ and that%
\[
\frac{d}{ds}[\Theta\left(  s,m\left(  s\right)  \right)  ] =\frac
{\partial\Theta}{\partial s}\left(  s,m\left(  s\right)  \right)  .
\]
The fact that $m\left(  s\right)  $ is a maximum point achieved in the
interior of $\left[  m_{1}\left(  t,N\right)  ,m_{2}\left(  t,N\right)
\right]  $ implies that%
\[
\partial_{mm}^{2}\Theta\left(  s,m\left(  s\right)  \right)  \leq0.
\]

\end{itemize}

Using the above relations we see that $s\rightarrow\Theta\left(  s,m\left(
s\right)  \right)  $ is non-increassing on $I$ and therefore it is controlled
by the value achieved in $t_{I,g}$ which is, in turn controlled by
$\max\left\{  \sup\limits_{\left[  N-1,N+2\right]  }\Theta\left(  0,m\right)
,\bar{C}\left(  t,E_{c},\left\Vert \rho_{0}\right\Vert _{L^{\infty}}\right)
\right\}  $, see $\left(  \text{\ref{bord_gauche}}\right)  $. We thus obtain a
contradiction of the very definition of $I$. The contradiction comes from the
fact that we assumed the set $A$ defined in $\left(  \text{\ref{set}}\right)
$ in non-empty. Therefore, for all $s\in\left[  0,t\right]  $ and all
$m\in\left[  m_{1}\left(  t,N\right)  ,m_{2}\left(  t,N\right)  \right]  $ we
have that%
\[
\Theta\left(  s,m\right)  \leq\max\left\{  \sup_{\left[  N-1,N+2\right]
}\Theta\left(  0,m\right)  ,\bar{C}\left(  t,E_{c},\left\Vert \rho
_{0}\right\Vert _{L^{\infty}}\right)  \right\}  .
\]
Since $\left[  N,N+1\right]  \subset\left[  m_{1}\left(  t,N\right)
,m_{2}\left(  t,N\right)  \right]  $ and $N\in\mathbb{Z}$ and $t>0$ where
chosen arbitrarily we conclude that
\[
\Theta\left(  t,m\right)  \leq C\left(  t,E_{c},\left\Vert \rho_{0}\right\Vert
_{L^{\infty}}\right)  \;\;\mbox{for all}\;\;m\in\mathbb{R}.
\]
Owing to the definition of $\Theta\left(  t,m\right)  $, the estimate
(\ref{estimGa}) on $\Gamma$, the energy estimates and the hypothesis $\left(
\text{\ref{H1}}\right)  $ and $\left(  \text{\ref{H2}}\right)  $ we conclude
that there exists a constant such that%
\begin{equation}
\widetilde{\rho}\left(  t,m\right)  \leq C\left(  t,E_{c},\left\Vert \rho
_{0}\right\Vert _{L^{\infty}}\right)  \;\;\mbox{for all}\;\;m\in\mathbb{R}.
\label{rho_is_bounded_conclusion}%
\end{equation}

In order to prove the second part of Proposition
\ref{rho_is_bounded_in_Linfty}, let us recall that relation $\left(
\text{\ref{Relation_long}}\right)  $ when $c=0$ reads%
\begin{align}
&  \frac{d}{dt}\left\{  \int_{N}^{N+1}\left(  \int_{\ell}^{m}\widetilde
{u}\right)  d\ell+r_{1}\psi(\widetilde{\rho})\left(  t,m\right)  \right\}
+\widetilde{\rho}^{\gamma}\left(  t,m\right) \nonumber\\
&  =\frac{d}{dt}\left\{  \int_{N}^{N+1}\left(  \int_{\ell}^{\ell+1}\frac
{dq}{\widetilde{\rho}\left(  t,q\right)  }\right)  ^{-1}\big( \int_{\ell
}^{\ell+1}\left(  \frac{1}{\widetilde{\rho}\left(  t,q\right)  }\int_{\ell
}^{q}\widetilde{u}\right)  dq \big) d\ell\right\} \nonumber\\
&  +\int_{N}^{N+1}\left[  \int_{\ell}^{\ell+1}\left(  \frac{1}{\widetilde
{\rho}\left(  t,q\right)  }\int_{\ell}^{q}\widetilde{u}\right)  dq\left(
\int_{\ell}^{\ell+1}\frac{dq}{\widetilde{\rho}\left(  t,q\right)  }\right)
^{-2}\left(  \widetilde{u}\left(  t,\ell+1\right)  -\widetilde{u}\left(
t,\ell\right)  \right)  \right]  d\ell\nonumber\\
&  +\int_{N}^{N+1}\left(  \int_{\ell}^{\ell+1}\frac{dq}{\widetilde{\rho
}\left(  t,q\right)  }\right)  ^{-1}\text{\ }\left\{  -\widetilde{u}\left(
t,\ell+1\right)  \int_{\ell}^{q}\widetilde{u}\left(  t,q^{\prime}\right)
dq^{\prime}+\int_{\ell}^{\ell+1}\left(  \widetilde{u}^{2}-\mu(\widetilde{\rho
})\partial_{m}\widetilde{u}+\widetilde{\rho}^{\gamma-1}\right)  \left(
t,q\right)  dq\right\}  d\ell\label{for_the_NS}%
\end{align}
for all $N\in\mathbb{Z}$ and for all $m\in\left[  N,N+1\right]  $. Hypothesis
$\left(  \text{\ref{H5}}\right)  $ is needed in order to have the validity of
Proposition \ref{R_auxilary_lemma_1.2} in order to apply the maximum
principle. However, as one can see from relation $\left(
\text{\ref{for_the_NS}}\right)  $, the situation for the NS system is simpler:
we just integrate in time \ref{R_auxilary_lemma_1.2} and proceed exactly as in
estimates $\left(  \text{\ref{R_rho_upper_Term1}}\right)  $,$\left(
\text{\ref{R_rho_upper_Term2}}\right)  $, $\left(
\text{\ref{R_rho_upper_Term4}}\right)  $ in order to conclude that the
validity of $\left(  \text{\ref{rho_is_bounded_conclusion}}\right)  $.

\subsection{Bounds for the effective
velocities\label{bounds_effective_section}}

In the following, we will deduce a one sided inequality concerning
$\widetilde{v}_{1}$, the equation of which can be put under the form%
\begin{equation}
\partial_{t}\widetilde{v}_{1}-\left(  1-r_{1}\right)  \partial_{m}\left(
\widetilde{\rho}\mu\left(  \widetilde{\rho}\right)  \partial_{m}\widetilde
{v}_{1}\right)  +\frac{\gamma\widetilde{\rho}^{\gamma}}{r_{1}\mu\left(
\widetilde{\rho}\right)  }\left(  \widetilde{v}_{1}-\widetilde{u}\right)  =0.
\label{R_equation_for_v1_2}%
\end{equation}

\begin{lemma}
\label{lemma_v_eff_is_bounded}Assume that $\mu$ satisfies the hypothesis
$\left(  \text{\ref{H1}}\right)  $, $\left(  \text{\ref{H2}}\right)  $,
$\left(  \text{\ref{H3}}\right)  $ and $\left(  \text{\ref{H6}}\right)  $.
Then, the following estimates hold true:

\begin{enumerate}
\item In the general NSK case: we also suppose that $\left(  \text{\ref{H5}%
}\right)  $ holds true. Then, we have that for all $m\in\mathbb{R}$
\begin{align*}
\widetilde{v}_{1}\left(  t,m\right)   &  \leq\max\left\{  \sup_{m\in
\mathbb{R}}\widetilde{v}_{1}\left(  0,m\right)  ,0\right\}  +\left(
1+t\right)  C\left(  E_{c}\right) \\
&  +C\left(  E_{c},\gamma\right)  \int_{0}^{t}\sup_{q\in\mathbb{R}}\Phi\left(
\frac{1}{\widetilde{\rho}\left(  s,q\right)  }\right)  \left(  \int_{-\infty
}^{+\infty}\widetilde{\rho}\left(  s\right)  \mu\left(  \widetilde{\rho
}\left(  s\right)  \right)  \left\vert \partial_{m}\widetilde{u}\left(
s\right)  \right\vert ^{2}\right)  ^{\frac{1}{2}}ds
\end{align*}

\item In the case of the NS system: we do not need to impose $\left(
\text{\ref{H5}}\right)  $ and, moreover, the following more precise statement
holds true for all $m\in\mathbb{R}$:%
\begin{align*}
\widetilde{v}_{1}\left(  t,m\right)   &  \leq\widetilde{v}_{1}\left(
0,m\right)  +\left(  1+t\right)  C\left(  E_{c}\right) \\
&  +C\left(  E_{c},\gamma\right)  \int_{0}^{t}\sup_{q\in\mathbb{R}}\Phi\left(
\frac{1}{\widetilde{\rho}\left(  s,q\right)  }\right)  \left(  \int_{-\infty
}^{+\infty}\widetilde{\rho}\left(  s\right)  \mu\left(  \widetilde{\rho
}\left(  s\right)  \right)  \left\vert \partial_{m}\widetilde{u}\left(
s\right)  \right\vert ^{2}\right)  ^{\frac{1}{2}}ds
\end{align*}

\end{enumerate}
\end{lemma}


\begin{corollary}
\label{corollary_v_eff_is_bounded_in_eulerian}Assume that $\mu$ satisfies the
hypothesis $\left(  \text{\ref{H1}}\right)  $, $\left(  \text{\ref{H2}%
}\right)  $, $\left(  \text{\ref{H3}}\right)  $ and $\left(  \text{\ref{H6}%
}\right)  $. Then, the following estimates hold true:

\begin{enumerate}
\item In the general NSK case: we also suppose that $\left(  \text{\ref{H5}%
}\right)  $ holds true. For all $x\in\mathbb{R}$ we have that%
\begin{align}
v_{1}\left(  t,x\right)   &  =u\left(  t,x\right)  +r_{1}\partial_{x}%
\varphi\left(  \rho\right)  \left(  t,x\right) \nonumber\\
&  \leq\max\left\{  \sup_{x\in\mathbb{R}}\left(  u_{0}\left(  x\right)
+r_{1}\partial_{x}\varphi\left(  \rho_{0}\left(  x\right)  \right)  \right)
,0\right\}  +C\left(  E_{c}\right)  \left(  1+t\right) \nonumber\\
&  +C\left(  E_{c},\gamma\right)  \int_{0}^{t}\sup_{z\in\mathbb{R}}\Phi\left(
\frac{1}{\rho\left(  s,z\right)  }\right)  \left(  \int_{-\infty}^{+\infty}%
\mu\left(  \rho\left(  s,z\right)  \right)  \left\vert \partial_{x}u\left(
s,z\right)  \right\vert ^{2}dz\right)  ^{\frac{1}{2}}ds.
\label{R_inegalite_v1}%
\end{align}

\item In the case of the NS system: we do not need to impose $\left(
\text{\ref{H5}}\right)  $ and for all $x\in\mathbb{R}$ we have that%
\begin{align}
v_{1}\left(  t,x\right)   &  =u\left(  t,x\right)  +r_{1}\partial_{x}%
\varphi\left(  \rho\right)  \left(  t,x\right) \nonumber\\
&  \leq u_{0}\left(  x\right)  +r_{1}\partial_{x}\varphi\left(  \rho
_{0}\right)  \left(  x\right)  +C\left(  E_{c}\right)  (1+t)\nonumber\\
&  +C\left(  E_{c},\gamma\right)  \int_{0}^{t}\sup_{z\in\mathbb{R}}\Phi\left(
\frac{1}{\rho\left(  s,z\right)  }\right)  \left(  \int_{-\infty}^{+\infty}%
\mu\left(  \rho\left(  s,z\right)  \right)  \left\vert \partial_{x}u\left(
s,z\right)  \right\vert ^{2}dz\right)  ^{\frac{1}{2}}ds. \label{Estimate_NS}%
\end{align}

\end{enumerate}
\end{corollary}

\begin{remark}
We mention that this is the key estimate leading to the fact that the density
is bounded away from vacuum. In some sense, in the same way as working with a
primitive for the momentum equation leads to an upper bound for the density,
working with a primitive of the equation of the effective velocity will enable
us to obtain a lower bound for the density.
\end{remark}

\begin{remark}
The local nature of the estimate $\left(  \text{\ref{Estimate_NS}}\right)  $
with respect to the initial data allows us to show that, in the case of the NS
system, the density is bounded away from vacuum without using any information
on the derivatives of $\rho_{0}$.
\end{remark}

\textit{Proof of Lemma \ref{lemma_v_eff_is_bounded}:} Let us notice that for
all $t>0$ and $m,q\in\left[  N,N+1\right]  $ we have that%
\begin{align}
\widetilde{u}\left(  t,m\right)  -\widetilde{u}\left(  t,q\right)   &
=\int_{q}^{m}\partial_{m}\widetilde{u}\leq\left(  \int_{q}^{m}\frac
{1}{\widetilde{\rho}\mu\left(  \widetilde{\rho}\right)  }\right)  ^{\frac
{1}{2}}\left(  \int_{q}^{m}\widetilde{\rho}\mu\left(  \widetilde{\rho}\right)
\left\vert \partial_{m}\widetilde{u}\right\vert ^{2}\right)  ^{\frac{1}{2}%
}\nonumber\\
&  =\left(  \int_{q}^{m}\frac{1}{\widetilde{\rho}}\Lambda\left(  \frac
{1}{\widetilde{\rho}}\right)  \right)  ^{\frac{1}{2}}\left(  \int_{q}%
^{m}\widetilde{\rho}\mu\left(  \widetilde{\rho}\right)  \left\vert
\partial_{m}\widetilde{u}\right\vert ^{2}\right)  ^{\frac{1}{2}}\nonumber\\
&  \leq\left(  \int_{q}^{m}\frac{1}{\widetilde{\rho}}\right)  ^{\frac{1}{2}%
}\left(  \int_{q}^{m}\widetilde{\rho}\mu\left(  \widetilde{\rho}\right)
\left\vert \partial_{m}\widetilde{u}\right\vert ^{2}\right)  ^{\frac{1}{2}%
}\left(  1+\sup\limits_{\ell\in\mathbb{R}}\Phi\left(  \frac{1}{\widetilde
{\rho}\left(  t,\ell\right)  }\right)  \right)  , \label{lower_bound_v1_1}%
\end{align}
where we recall that $\Lambda$ and $\Phi$ are defined in $\left(
\text{\ref{Var_Phi}}\right)  $ and that hypothesis $\left(  \text{\ref{H6}%
}\right)  $ reads:%
\[
\forall\tau>0:\Lambda\left(  \tau\right)  =\frac{1}{\mu\left(  \dfrac{1}{\tau
}\right)  }\leq\left(  1+\Phi\left(  \tau\right)  \right)  ^{2}.
\]
Integrating the relation $\left(  \text{\ref{lower_bound_v1_1}}\right)  $
w.r.t. $q$ on $\left[  N,N+1\right]  $ and using Proposition
\ref{R_auxilary_lemma_1} yields%
\begin{align}
\widetilde{u}\left(  t,m\right)   &  \leq\int_{N}^{N+1}\left\vert
\widetilde{u}\left(  t,q\right)  \right\vert dq+\left(  1+\sup\limits_{q\in
\mathbb{R}}\Phi\left(  \frac{1}{\widetilde{\rho}\left(  t,q\right)  }\right)
\right)  \int_{N}^{N+1}\left(  \int_{q}^{m}\frac{1}{\widetilde{\rho}}\right)
^{\frac{1}{2}}\left(  \int_{q}^{m}\widetilde{\rho}\mu\left(  \widetilde{\rho
}\right)  \left\vert \partial_{m}\widetilde{u}\right\vert ^{2}\right)
^{\frac{1}{2}}dq\nonumber\\
&  \leq C\left(  E_{c}\right)  +C\left(  E_{c},\gamma\right)  \left(
1+\sup\limits_{q\in\mathbb{R}}\Phi\left(  \frac{1}{\widetilde{\rho}\left(
t,q\right)  }\right)  \right)  \left(  \int_{-\infty}^{+\infty}\widetilde
{\rho}\mu\left(  \widetilde{\rho}\right)  \left\vert \partial_{m}\widetilde
{u}\right\vert ^{2}\right)  ^{\frac{1}{2}}. \label{upper_u}%
\end{align}
Since $m\in\left[  N,N+1\right]  $, the RHS of the above estimate is
independent of $N$ and $N$ was arbitrarly fixed we conclude that $\left(
\text{\ref{upper_u}}\right)  $ holds true for all $m\in\mathbb{R}$. Using
$\left(  \text{\ref{upper_u}}\right)  $ in $\left(
\text{\ref{R_equation_for_v1_2}}\right)  $ along with the fact that
$\widetilde{\rho}$ is bounded and (\ref{H3}), we obtain that%
\begin{align}
&  \partial_{t}\widetilde{v}_{1}-\left(  1-r_{1}\right)  \partial_{m}\left(
\widetilde{\rho}\mu\left(  \widetilde{\rho}\right)  \partial_{m}\widetilde
{v}_{1}\right)  +\frac{\gamma\widetilde{\rho}^{\gamma}}{r_{1}\mu\left(
\widetilde{\rho}\right)  }\widetilde{v}_{1}\nonumber\\
&  \leq C\left(  E_{c}\right)  +C\left(  E_{c},\gamma\right)  \left(
1+\sup\limits_{q\in\mathbb{R}}\Phi\left(  \frac{1}{\widetilde{\rho}\left(
t,q\right)  }\right)  \right)  \left(  \int_{-\infty}^{+\infty}\widetilde
{\rho}\mu\left(  \widetilde{\rho}\right)  \left\vert \partial_{m}\widetilde
{u}\right\vert ^{2}\right)  ^{\frac{1}{2}}. \label{upper_v}%
\end{align}
We denote by
\[
\Psi\left(  t\right)  =tC\left(  E_{c}\right)  +C\left(  E_{c},\gamma\right)
\int_{0}^{t}\left(  1+\sup\limits_{q\in\mathbb{R}}\Phi\left(  \frac
{1}{\widetilde{\rho}\left(  s,q\right)  }\right)  \right)  \left(
\int_{-\infty}^{+\infty}\widetilde{\rho}\mu\left(  \widetilde{\rho}\right)
\left\vert \partial_{m}\widetilde{u}\right\vert ^{2}\right)  ^{\frac{1}{2}}
ds,
\]
the primitive vanishing at $t=0$ of the RHS term appearing in $\left(
\text{\ref{upper_v}}\right)  $. The inequality $\left(  \text{\ref{upper_v}%
}\right)  $ can be put under the following form:%
\[
\frac{d}{dt}\left(  \widetilde{v}_{1}-\Psi\right)  -\left(  1-r_{1}\right)
\partial_{m}\left(  \widetilde{\rho}\mu\left(  \widetilde{\rho}\right)
\partial_{m}(\widetilde{v}_{1}-\Psi)\right)  +\frac{\gamma\widetilde{\rho
}^{\gamma}}{r_{1}\mu\left(  \widetilde{\rho}\right)  }(\widetilde{v}_{1}%
-\Psi)\leq-\frac{\gamma\widetilde{\rho}^{\gamma}}{r_{1}\mu\left(
\widetilde{\rho}\right)  }\Psi\leq0.
\]
We remark that owing to the fact that
\[
\int_{\mathbb{R}}\widetilde{v}_{1}^{2}\left(  t,m\right)  dm<\infty,
\]
we have that for all $t\geq0$:%
\[
\lim_{\left\vert m\right\vert \rightarrow\infty}\left(  \widetilde{v}%
_{1}\left(  t,m\right)  -\Psi\left(  t\right)  \right)  =-\Psi\left(
t\right)  \leq0.
\]
The maximum principle then implies that%
\[
\widetilde{v}_{1}\left(  t,m\right)  -\Psi\left(  t\right)  \leq\max\left\{
\sup_{q\in\mathbb{R}}\widetilde{v}_{1}\left(  0,q\right)  ,\sup_{s\in\left[
0,t\right]  }\left(  -\Psi\left(  s\right)  \right)  \right\}  \leq
\max\left\{  \sup_{q\in\mathbb{R}}\widetilde{v}_{1}\left(  0,q\right)
,0\right\}
\]
and accordingly we get that%
\begin{align}
\widetilde{v}_{1}\left(  t,m\right)   &  \leq\max\left\{  \sup_{q\in
\mathbb{R}}\widetilde{v}_{1}\left(  0,q\right)  ,0\right\}  +\Psi\left(
t\right) \nonumber\\
&  \leq\max\left\{  \sup_{q\in\mathbb{R}}\widetilde{v}_{1}\left(  0,q\right)
,0\right\}  +C\left(  t,E_{c}\right) \nonumber\\
&  + C\left(  t,E_{c}\right)  \int_{0}^{t}\sup_{q\in\mathbb{R}}\Phi\left(
\frac{1}{\widetilde{\rho}\left(  s,q\right)  }\right)  \left(  \int_{-\infty
}^{+\infty}\widetilde{\rho}\left(  s\right)  \mu\left(  \widetilde{\rho
}\left(  s\right)  \right)  \left\vert \partial_{m}\widetilde{u}\left(
s\right)  \right\vert ^{2}\right)  ^{\frac{1}{2}}ds. \label{R_v1_inf}%
\end{align}

In order to prove the second part of Lemma \ref{lemma_v_eff_is_bounded}, we
remark that in the case of the Navier-Stokes system i.e. $r_{1}\left(
c\right)  =1,$ with $c=0$ we can get a much more precise estimate since, in
this case, $\left(  \text{\ref{R_equation_for_v1_2}}\right)  $ does not have a
diffusion operator. Thus, we see that there is no need to invoke the maximum
principle and the estimate is localized: for all $m\in\mathbb{R}$ we get using
(\ref{R_equation_for_v1_2}), (\ref{rho_is_bounded_conclusion}), (\ref{H3}),
(\ref{upper_u}) that%
\begin{align*}
\widetilde{v}_{1}\left(  t,m\right)   &  \leq\widetilde{v}_{1}\left(
0,m\right)  +\int_{0}^{t}\frac{\gamma\widetilde{\rho}^{\gamma}\left(
s,m\right)  }{r_{1}\mu\left(  \widetilde{\rho}\left(  s,m\right)  \right)
}\widetilde{u}\left(  s,m\right)  \exp\left(  -\int_{s}^{t}\frac
{\gamma\widetilde{\rho}^{\gamma}}{\mu\left(  \widetilde{\rho}\right)
}\right)  ds\\
&  \leq\widetilde{v}_{1}\left(  0,m\right)  +C\left(  t,E_{c}\right)  \int
_{0}^{t}\left\Vert \widetilde{u}\left(  s\right)  \right\Vert _{L^{\infty}%
}ds\\
&  \leq\widetilde{v}_{1}\left(  0,m\right)  +C\left(  t,E_{c}\right)
+C\left(  E_{c}\right)  \int_{0}^{t}\left(  1+\sup_{q\in\mathbb{R}}\Phi\left(
\frac{1}{\widetilde{\rho}\left(  s,q\right)  }\right)  \right)  \left(
\int_{-\infty}^{+\infty}\widetilde{\rho}\left(  s\right)  \mu\left(
\widetilde{\rho}\left(  s\right)  \right)  \left\vert \partial_{m}%
\widetilde{u}\left(  s\right)  \right\vert ^{2}\right)  ^{\frac{1}{2}}ds\\
&  \leq\widetilde{v}_{1}\left(  0,m\right)  +C\left(  t,E_{c}\right) \\
&  +C\left(  E_{c},\gamma\right)  \int_{0}^{t}\sup_{q\in\mathbb{R}}\Phi\left(
\frac{1}{\widetilde{\rho}\left(  s,q\right)  }\right)  \left(  \int_{-\infty
}^{+\infty}\widetilde{\rho}\left(  s\right)  \mu\left(  \widetilde{\rho
}\left(  s\right)  \right)  \left\vert \partial_{m}\widetilde{u}\left(
s\right)  \right\vert ^{2}\right)  ^{\frac{1}{2}}ds.
\end{align*}
This concludes the proof of Lemma \ref{lemma_v_eff_is_bounded}$\Box$

The result announced in Corollary is obtained just by returning back to the
Eulerian variables, taking in consideration that for all $x\in\mathbb{R}$ it
holds true that%
\[
v_{1}\left(  t,x\right)  =\widetilde{v}_{1}\left(  t,m\left(  t,x\right)
\right)
\]
where $m\left(  t,x\right)  $ is the inverse of the function $m\rightarrow
X(t,Y(m))$.

In the next section we show how to use the estimates provided by Corollary
\ref{corollary_v_eff_is_bounded_in_eulerian} in order to obtain the fact that
the density is bounded away from vacuum.

\subsection{Lower bound for the density\label{lower_bound_for_density_section}%
}

In the following, we show how it is possible to find a lower bound for the
density with the use of the inequality $\left(  \text{\ref{R_v1_inf}}\right)
$. At this point, we need an Eulerian-equivalent for Proposition
\ref{R_auxilary_lemma_1}. We claim that:

\begin{proposition}
\label{R_auxilary_lemma_2}Consider $\ell\in\mathbb{R}$. Then, for all $t\geq0$
there exists a point $y\left(  t,\ell\right)  \in\left[  \ell,\ell
+2E_{0}\right]  $ such that
\[
\rho\left(  t,y\left(  t,\ell\right)  \right)  \geq\pi^{-1}\left(  1/2\right)
,
\]
where we recall that $\pi\left(  \rho\right)  =\rho e\left(  \rho\right)  $.
\end{proposition}

\textit{Proof of Proposition} \ref{R_auxilary_lemma_2}. Recall that according
to the energy inequality for all $t\geq0$ we have%
\[
\int_{\ell}^{\ell+2E_{0}}\rho\left(  t,x\right)  e\left(  \rho\left(
t,x\right)  \right)  dx\leq E_{0}.
\]
Assume that for all $x\in\left[  \ell,\ell+2E_{0}\right]  $ one has%
\[
\rho\left(  t,x\right)  <\varepsilon,
\]
with $0<\varepsilon\leq1$ then, as the function $\pi$ is strictly decreasing
on $\left[  0,1\right]  $ we find that%
\[
\pi\left(  \rho\left(  t,x\right)  \right)  \geq\pi\left(  \varepsilon\right)
.
\]
Integrating the previous inequality on $\left[  \ell,\ell+2E_{0}\right]  $
yields%
\[
E_{0}\geq2\pi\left(  \varepsilon\right)  E_{0},
\]
such that
\[
\varepsilon\geq\pi^{-1}\left(  1/2\right)  ,
\]
with $\pi^{-1}\left(  1/2\right)  $ is here the unique element $\alpha$ of
$[0,1]$ such that $\pi(\alpha)=\frac{1}{2}$. The conclusion follows. $\Box$
\newline\newline We are now in the position of obtaining a lower bound for the density.

\begin{proposition}
\label{rho_is_bounded_by_bellow} Assume that $\mu$ satisfies the hypothesis
$\left(  \text{\ref{H1}}\right)  $, $\left(  \text{\ref{H2}}\right)  $,
$\left(  \text{\ref{H3}}\right)  $, $\left(  \text{\ref{H4}}\right)  $,
$\left(  \text{\ref{H6}}\right)  $ and $\left(  \text{\ref{H7}}\right)  $.

\begin{enumerate}
\item In the general NSK case we also suppose that $\left(  \text{\ref{H5}%
}\right)  $ holds true. Suppose that%
\[
a.e.\text{ }x\in\mathbb{R}\text{ : }u_{0}\left(  x\right)  +\partial
_{x}\varphi\left(  \rho_{0}\left(  x\right)  \right)  \leq M_{0},
\]
with $M_{0}\in\mathbb{R}$. Then, we have that for all $x\in\mathbb{R}$:%
\begin{equation}
C\left(  t,E_{c},M_{0},\left\Vert (\rho_{0},\frac{1}{\rho_{0}})\right\Vert
_{L^{\infty}}\right)  \leq\rho\left(  t,x\right)  . \label{rho_lower_NSK}%
\end{equation}

\item In the case when $c=0$, i.e. for the NS system, we have that for all
$x\in\mathbb{R}$:
\begin{equation}
C\left(  t,E_{c},\left\Vert (\rho_{0},\frac{1}{\rho_{0}})\right\Vert
_{L^{\infty}}\right)  \leq\rho\left(  t,x\right)  . \label{rho_lower_NS}%
\end{equation}

\end{enumerate}
\end{proposition}

\begin{remarka}
At this level, it is important to mention that we can obtain similar estimates
if we assume that%
\[
M_{1}\leq u_{0}+r_{1}\left(  c\right)  \partial_{x}\varphi(\rho_{0})
\]
for some $M_{1}\in\mathbb{R}$. It suffices to repeat the same argument except
that we take $y\in\lbrack N-2E_{0}(c),N]$ and we integrate between $x$ and
$y$. \label{rem5}
\end{remarka}

\begin{remarka}
\label{remarca_continuare}In the case of the NSK system, we also have that%
\begin{equation}
u\left(  t,x\right)  +r_{1}\left(  c\right)  \partial_{x}\varphi\left(
\rho\left(  t,x\right)  \right)  \leq C\left(  t,E_{c},M_{0},\left\Vert
(\rho_{0},\frac{1}{\rho_{0}})\right\Vert _{L^{\infty}}\right)
\label{propagation_one_sided_for_v1}%
\end{equation}

\end{remarka}

\textit{Proof of Proposition \ref{rho_is_bounded_by_bellow}}: As before, let
$N\in\mathbb{Z}$, consider any $x\in\left[  N,N+1\right]  $ and $y\in\left[
N+1,N+1+2E_{c}\right]  $ (we have in particular $x\leq y$) and integrate
$\left(  \text{\ref{R_inegalite_v1}}\right)  $ between $x$ and $y$ in order to
obtain that%
\begin{align}
\int_{x}^{y}u\left(  t,z\right)  dz+  &  r_{1}(\varphi\left(  \rho(t,y\right)
)-\varphi\left(  \rho\left(  t,x\right)  \right) \nonumber\\
&  \leq(y-x)\biggl( \max(M_{0},0)+C\left(  E_{c}\right)  (1+t)\nonumber\\
&  +C\left(  E_{c},\gamma\right)  \int_{0}^{t}\sup_{z\in\mathbb{R}}\Phi\left(
\frac{1}{\rho\left(  s,z\right)  }\right)  \left(  \int_{-\infty}^{+\infty}%
\mu\left(  \rho\left(  s,z\right)  \right)  \left\vert \partial_{x}u\left(
s,z\right)  \right\vert ^{2}dz\right)  ^{\frac{1}{2}}ds\biggl)\nonumber\\
&  \leq(2E_{c}+1)\max(M_{0},0)+C\left(  E_{c}\right)  \left(  1+t\right)
\nonumber\\
&  +C\left(  E_{c},\gamma\right)  \int_{0}^{t}\sup_{z\in\mathbb{R}}\Phi\left(
\frac{1}{\rho\left(  s,z\right)  }\right)  \left(  \int_{-\infty}^{+\infty}%
\mu\left(  \rho\left(  s,z\right)  \right)  \left\vert \partial_{x}u\left(
s,z\right)  \right\vert ^{2}dz\right)  ^{\frac{1}{2}}ds. \label{estimation}%
\end{align}
Using Proposition \ref{R_auxilary_lemma_2}, there exists $y_{N}\left(
t\right)  \in$ $\left[  N+1,N+1+2E_{c}\right]  $ with the property that%
\[
\rho\left(  t,y_{N}\left(  t\right)  \right)  \geq\pi^{-1}\left(  1/2\right)
.
\]
We take $y=$ $y_{N}\left(  t\right)  $ in $\left(  \text{\ref{estimation}%
}\right)  $ and recalling that%
\[
\Phi\left(  \tau\right)  =-\varphi\left(  \frac{1}{\tau}\right)
\]
we get that
\begin{align}
r_{1}\Phi\left(  \frac{1}{\rho\left(  t,x\right)  }\right)  \leq &  r_{1}%
\Phi\left(  \frac{1}{\rho\left(  t,y_{N}\left(  t\right)  \right)  }\right)
+(2E_{c}+1)\max(M_{0},0) +C\left(  E_{c},\gamma\right)  \left(  1+t\right)
\nonumber\\
&  -\int_{x}^{y_{N}\left(  t\right)  }u\left(  t,z\right)  dz+C\left(
E_{c},\gamma\right)  \int_{0}^{t}\sup_{z\in\mathbb{R}}\Phi\left(  \frac
{1}{\rho\left(  t,z\right)  }\right)  \left(  \int_{-\infty}^{+\infty}%
\mu\left(  \rho\right)  \left\vert \partial_{x}u\right\vert ^{2}\right)
^{\frac{1}{2}}. \label{R_bound_inf_rho}%
\end{align}
The only term that needs to be treated in $\left(  \text{\ref{R_bound_inf_rho}%
}\right)  $ is the integral of $u$ over $\left[  x,y_{N}\left(  t\right)
\right]  $. This is done in the following lines using hypothesis $\left(
\text{\ref{H7}}\right)  $ we infer that for all $\varepsilon>0$ we have:
\begin{align*}
-\int_{x}^{y_{N}(t)}u\left(  t,z\right)  dz  &  =-\int_{x}^{y_{N}\left(
t\right)  }\frac{1}{\sqrt{\rho\left(  t,z\right)  }}\sqrt{\rho\left(
t,z\right)  }u\left(  t,z\right)  dz\\
&  \leq\int_{x}^{y_{N}\left(  t\right)  }\frac{1}{\sqrt{\rho\left(
t,z\right)  }}\left\vert \sqrt{\rho\left(  t,z\right)  }u\left(  t,z\right)
\right\vert dz\\
&  \leq\int_{x}^{y_{N}\left(  t\right)  }\left(  1+\Phi\left(  \frac{1}%
{\rho\left(  t,z\right)  }\right)  \right)  ^{1-\eta}\left\vert \sqrt
{\rho\left(  t,z\right)  }u\left(  t,z\right)  \right\vert dz\\
&  \leq\left(  1+\sup_{z\in\mathbb{R}}\Phi\left(  \frac{1}{\rho\left(
t,z\right)  }\right)  \right)  ^{1-\eta}\int_{x}^{y_{N}\left(  t\right)
}\left\vert \sqrt{\rho\left(  t,z\right)  }u\left(  t,z\right)  \right\vert
dz\\
&  \leq\left(  1+\sup_{z\in\mathbb{R}}\Phi\left(  \frac{1}{\rho\left(
t,z\right)  }\right)  \right)  ^{1-\eta}\left(  2E_{c}+1\right)  ^{\frac{1}%
{2}} \sqrt{E_{c}}\\
&  \leq\frac{C\left(  E_{c},\eta\right)  }{4\varepsilon}+\varepsilon\sup
_{z\in\mathbb{R}}\Phi\left(  \frac{1}{\rho\left(  t,z\right)  }\right)  .
\end{align*}
Thus, we get that for all $x\in\mathbb{R}$ and all $\varepsilon>0$ we have
that%
\begin{align}
r_{1}\Phi\left(  \frac{1}{\rho\left(  t,x\right)  }\right)   &  \leq C\left(
t,E_{c},M_{0},\left\Vert (\rho_{0},\frac{1}{\rho_{0}})\right\Vert _{L^{\infty
}}\right) \nonumber\\
+  &  \frac{C\left(  E_{c},\eta\right)  }{4\varepsilon}+\varepsilon\sup
_{z\in\mathbb{R}}\Phi\left(  \frac{1}{\rho\left(  t,z\right)  }\right)
\nonumber\\
&  +C\left(  E_{c},\gamma\right)  \int_{0}^{t}\sup_{z\in\mathbb{R}}\Phi\left(
\frac{1}{\rho\left(  t,z\right)  }\right)  \left(  \int_{-\infty}^{+\infty}%
\mu\left(  \rho\right)  \left\vert \partial_{x}u\right\vert ^{2}\right)
^{\frac{1}{2}}.
\end{align}
Taking the supremum w.r.t. $x\in\mathbb{R}$ in the LHS, taking $\varepsilon
=r_{1}/2$ and using Gr\"{o}nwall's inequality we obtain that%
\[
\Phi\left(  \frac{1}{\rho\left(  t,x\right)  }\right)  \leq C\left(
t,E_{c},M_{0},\left\Vert (\rho_{0},\frac{1}{\rho_{0}})\right\Vert _{L^{\infty
}}\right)  .
\]
Thus using the hypothesis $\left(  \text{\ref{H4}}\right)  $ we get that
$\rho$ is bounded below and that we have for any $(t,x)\in\mathbb{R}^{+}%
\times\mathbb{R}$%
\begin{equation}
C\left(  t,E_{c},M_{0},\left\Vert (\rho_{0},\frac{1}{\rho_{0}})\right\Vert
_{L^{\infty}}\right)  \leq\rho\left(  t,x\right)  . \label{R_bound_inf_rho_2}%
\end{equation}

The second part of Proposition \textit{\ref{rho_is_bounded_by_bellow}} follows
owing to the fact that in the case of the Navier-Stokes system, a more precise
information can be obtained (we recall that $r_{1}=1$ in this case). Recall
that according to $\left(  \text{\ref{Estimate_NS}}\right)  $, for all
$x\in\mathbb{R}$ we have that
\begin{align}
v_{1}\left(  t,x\right)   &  =u\left(  t,x\right)  +r_{1}\partial_{x}%
\varphi\left(  \rho\right)  \left(  t,x\right) \nonumber\\
&  \leq u_{0}\left(  x\right)  +r_{1}\partial_{x}\varphi\left(  \rho
_{0}\right)  \left(  x\right)  +C(E_{c},\gamma)\left(  1+t\right) \\
&  +C\left(  E_{c},\gamma\right)  \int_{0}^{t}\sup_{z\in\mathbb{R}}\Phi\left(
\frac{1}{\rho\left(  t,z\right)  }\right)  \left(  \int_{-\infty}^{+\infty}%
\mu\left(  \rho\right)  \left\vert \partial_{x}u\right\vert ^{2}\right)
^{\frac{1}{2}}.
\end{align}
Choosing as above $x\in\left[  N,N+1\right]  ,$ $y_{N}\left(  t\right)
\in\lbrack N+1,N+1+2E_{c}]$ with%
\[
\rho\left(  t,y_{N}\left(  t\right)  \right)  \geq\pi^{-1}\left(  1/2\right)
\]
and integrating we get that%
\begin{align*}
&  \int_{x}^{y_{N}\left(  t\right)  }u\left(  t,z\right)  dz+r_{1}\Phi\left(
\frac{1}{\rho\left(  t,x\right)  }\right)  -r_{1}\Phi\left(  \frac{1}%
{\rho\left(  t,y_{N}\left(  t\right)  \right)  }\right) \\
&  \leq\int_{x}^{y_{N}\left(  t\right)  }\left(  u_{0}\left(  z\right)
+r_{1}\partial_{x}\varphi\left(  \rho_{0}\right)  \left(  z\right)  \right)
dz\\
&  C\left(  E_{c}\right)  \left(  1+t\right)  +C\left(  E_{c},\gamma\right)
\int_{0}^{t}\sup_{x\in\mathbb{R}}\Phi\left(  \frac{1}{\rho\left(  t,x\right)
}\right)  \left(  \int_{-\infty}^{+\infty}\mu\left(  \rho\right)  \left\vert
\partial_{x}u\right\vert ^{2}\right)  ^{\frac{1}{2}}\\
&  \leq\int_{x}^{y_{N}\left(  t\right)  }u_{0}\left(  z\right)  dz+r_{1}%
\left(  \varphi\left(  \rho_{0}\left(  y_{N}\left(  t\right)  \right)
-\varphi(\rho_{0}(x)\right)  \right)  +C\left(  E_{c},\gamma\right)  \int
_{0}^{t}\sup_{z\in\mathbb{R}}\Phi\left(  \frac{1}{\rho\left(  t,z\right)
}\right)  \left(  \int_{-\infty}^{+\infty}\mu\left(  \rho\right)  \left\vert
\partial_{x}u\right\vert ^{2}\right)  ^{\frac{1}{2}}.
\end{align*}
We deduce that%
\begin{align*}
r_{1}\Phi\left(  \frac{1}{\rho\left(  t,x\right)  }\right)   &  \leq r_{1}%
\Phi\left(  \frac{1}{\rho\left(  t,y_{N}\left(  t\right)  \right)  }\right) \\
&  +\int_{x}^{y_{N}\left(  t\right)  }u_{0}\left(  z\right)  dz+r_{1}\left\{
\varphi(\rho_{0}(y_{N}\left(  t\right)  ))-\varphi\left(  \rho_{0}(x)\right)
\right\} \\
&  -\int_{x}^{y_{N}\left(  t\right)  }u\left(  t,z\right)  dz\\
&  +C\left(  E_{c},\gamma\right)  \int_{0}^{t}\sup_{x\in\mathbb{R}}\Phi\left(
\frac{1}{\rho\left(  t,x\right)  }\right)  \left(  \int_{-\infty}^{+\infty}%
\mu\left(  \rho\right)  \left\vert \partial_{x}u\right\vert ^{2}\right)
^{\frac{1}{2}},
\end{align*}
from which it yields%
\begin{align*}
r_{1}\Phi\left(  \frac{1}{\rho\left(  t,x\right)  }\right)   &  \leq r_{1}%
\Phi\left(  \frac{1}{\rho\left(  t,y_{N}\left(  t\right)  \right)  }\right)
+(2E_{0}+1)^{\frac{1}{2}}\left\Vert u_{0}\right\Vert _{L^{2}}+2\sup
_{p\in\lbrack\frac{1}{\left\Vert \rho_{0}\right\Vert _{L^{\infty}}},\left\Vert
\rho_{0}\right\Vert _{L^{\infty}}]}\left\vert \varphi\left(  p\right)
\right\vert \\
&  -\int_{x}^{y_{N}\left(  t\right)  }u\left(  t,z\right)  dz\\
&  +C\left(  E_{c},\gamma\right)  \int_{0}^{t}\sup_{x\in\mathbb{R}}\Phi\left(
\frac{1}{\rho\left(  t,x\right)  }\right)  \left(  \int_{-\infty}^{+\infty}%
\mu\left(  \rho\right)  \left\vert \partial_{x}u\right\vert ^{2}\right)
^{\frac{1}{2}}%
\end{align*}
and since the terms on the RHS can be controlled we use Gr\"{o}nwall's lemma
in order to obtain that%
\[
\Phi\left(  \frac{1}{\rho\left(  t,x\right)  }\right)  \leq C\left(
t,E_{c},\left\Vert (\rho_{0},\frac{1}{\rho_{0}})\right\Vert _{L^{\infty}%
}\right)  ,
\]
from which we can deduce that%
\[
C\left(  t,E_{c},\left\Vert (\rho_{0},\frac{1}{\rho_{0}})\right\Vert
_{L^{\infty}}\right)  \leq\rho\left(  t,x\right)  .
\]
This concludes the proof of Proposition \ref{rho_is_bounded_by_bellow}. $\Box$

The estimate from $\left(  \text{\ref{propagation_one_sided_for_v1}}\right)  $
follows by combining the conclusion of Proposition
\ref{rho_is_bounded_by_bellow} with the inequality $\left(
\text{\ref{R_inegalite_v1}}\right)  $.

\subsection{Global existence of the approximate solutions}

Now we can come back to the sequence of solution $(\rho^{n},u^{n}%
)_{n\in\mathbb{N}}$ of the Navier-Stokes Korteweg system (\ref{NSK_intro}) and
of the compressible Navier-Stokes system (\ref{Navier_Stokes_1d}) constructed
in the subsection \ref{subsec1}. We recall that these are defined on a finite
time interval $(0,T_{n})$ for each $n\in\mathbb{N}$.

In view of Proposition \ref{rho_is_bounded_in_Linfty} and the estimate
$\left(  \text{\ref{rho_borne2}}\right)  $ from Remark \ref{remarca_rho} that
follows we obtain that:
\begin{equation}
\Vert\rho^{n}(t,\cdot)\Vert_{L^{\infty}(\mathbb{R})}\leq C\left(
t,E_{c}\left(  \rho_{0}^{n},u_{0}^{n}\right)  ,\left\Vert \rho_{0}%
^{n}\right\Vert _{L^{\infty}}\right)  <+\infty.
\end{equation}
for any $n\in\mathbb{N}$ and $t\in(0,T_{m})$. Next, owing to Proposition
\textit{\ref{rho_is_bounded_by_bellow},} and (\ref{rho_is_bounded_conclusion})
we know that for all $t\in(0,T_{n})$ and $x\in\mathbb{R}$ we have that:
\begin{equation}
\begin{aligned} &\left\Vert \frac{1}{\rho^{n}\left( t,\cdot\right) }\right\Vert _{L^{\infty }\left( \mathbb{R}\right) }\leq C\left( t,E_{c}\left( \rho_{0}^{n},u_{0}^{n}\right) ,M_{0}^{n},\left\Vert \left( \rho_{0}^{n},\frac{1}{\rho_{0}^{n}}\right) \right\Vert _{L^{\infty}}\right) ,\\ &\left\Vert\rho^{n}\left( t,\cdot\right) \right\Vert _{L^{\infty }\left( \mathbb{R}\right) }\leq C\left( t,E_{c}\left( \rho_{0}^{n},u_{0}^{n}\right) ,\left\Vert \rho_{0}^{n} \right\Vert _{L^{\infty}}\right) \end{aligned} \label{cru}%
\end{equation}
where%
\[
v_{1|t=0}^{n}\leq M_{0}^{n},
\]
with $C$ a continuous function for the NSK system (there is no $M_{0}^{n}$ in
the case of the NS system). Of course, owing to subsection \ref{subsec1} we
know that%
\[
E_{c}\left(  \rho_{0}^{n},u_{0}^{n}\right)  ,M_{0}^{n},\left\Vert \left(
\rho_{0}^{n},\frac{1}{\rho_{0}^{n}}\right)  \right\Vert _{L^{\infty}}%
\]
are uniformly bounded in $n$ and we thus deduce, that, for a larger $C$ if
necessary, we have that
\begin{equation}
\Vert\rho^{n}(t,\cdot)\Vert_{L^{\infty}(\mathbb{R})}\leq C\left(
t,E_{c}\left(  \rho_{0},u_{0}\right)  ,\left\Vert \rho_{0}\right\Vert
_{L^{\infty}}\right)  <+\infty. \label{cru2}%
\end{equation}
and that%
\begin{equation}
\left\Vert \frac{1}{\rho^{n}\left(  t,\cdot\right)  }\right\Vert _{L^{\infty
}\left(  \mathbb{R}\right)  }\leq C\left(  t,E_{c}\left(  \rho_{0}%
,u_{0}\right)  ,M_{0},\left\Vert \left(  \rho_{0},\frac{1}{\rho_{0}}\right)
\right\Vert _{L^{\infty}}\right)  \label{cru1}%
\end{equation}
The blow-up criterion of the Theorem \ref{Cons} and $\left(  \text{\ref{cru2}%
}\right)  $, $\left(  \text{\ref{cru1}}\right)  $ imply that for any
$n\in\mathbb{N}$ we have $T_{n}=+\infty$ and that $\left(  \text{\ref{cru2}%
}\right)  $ and $\left(  \text{\ref{cru1}}\right)  $ hold true in fact for all
$t\geq0$.

In the case of the NSK system, for all $n\in\mathbb{N}$ and all $t\geq0$ we
also have
\begin{equation}
u^{n}+r_{1}\left(  c\right)  \partial_{x}\varphi\left(  \rho_{0}^{n}\right)
\leq C\left(  t,E_{c}\left(  \rho_{0},u_{0}\right)  ,M_{0},\left\Vert \left(
\rho_{0},\frac{1}{\rho_{0}}\right)  \right\Vert _{L^{\infty}}\right)
\label{cru3}%
\end{equation}

\subsection{Total variation estimates for solutions of the NSK
system\label{total_variation_section}}

Using the estimate $\left(  \text{\ref{cru1}}\right)  $ we deduce immediately
for any $n\in\mathbb{N}$, an estimate for the total variation of $\varphi
(\rho_{n})$. Indeed recall that $\left(  \text{\ref{cru3}}\right)  $ states
that%
\begin{equation}
u^{n}\left(  t,x\right)  +r_{1}\left(  c\right)  \partial_{x}\varphi\left(
\rho^{n}\left(  t,x\right)  \right)  \leq C\left(  t,E_{c}\left(  \rho
_{0},u_{0}\right)  ,M_{0},\left\Vert \left(  \rho_{0},\frac{1}{\rho_{0}%
}\right)  \right\Vert _{L^{\infty}}\right)  \label{cru4}%
\end{equation}
This implies that for all $n\in\mathbb{N}$, for all $t\geq0$, the function
\[
x\rightarrow\int_{0}^{x}u^{n}\left(  t,z\right)  dz+r_{1}\left(  c\right)
\varphi\left(  \rho^{n}\left(  t,x\right)  \right)  -C\left(  t,E_{c}\left(
\rho_{0},u_{0}\right)  ,M_{0},\left\Vert \left(  \rho_{0},\frac{1}{\rho_{0}%
}\right)  \right\Vert _{L^{\infty}}\right)  x
\]
is deacreasing and, taking in account $\left(  \text{\ref{cru1}}\right)  $, it
is of bounded variation on any compact. Using the energy estimate along with
$\left(  \text{\ref{cru1}}\right)  $ we obtain that:
\begin{equation}
TV\left(  \left[  -L,L\right]  ,\varphi\left(  \rho^{n}\left(  t,x\right)
\right)  \right)  \leq C\left(  t,E_{c}\left(  \rho_{0},u_{0}\right)
,M_{0},\left\Vert \left(  \rho_{0},\frac{1}{\rho_{0}}\right)  \right\Vert
_{L^{\infty}}\right)  \left(  1+L\right)  , \label{R_rho_bv}%
\end{equation}
for all $t\geq0$ and $L>0$ with $TV\left(  \left[  -L,L\right]  ,\varphi
\left(  \rho^{n}\left(  t,\cdot\right)  \right)  \right)  $ the total
variation of $\varphi^{n}(t,\cdot)$ on $[-L,L]$. Moreover, we have that%
\[
\partial_{t}\varphi(\rho^{n})=-u^{n}\partial_{x}\varphi(\rho^{n})+\frac
{\mu\left(  \rho^{n}\right)  }{\rho^{n}}\partial_{x}u^{n}%
\]
such that%
\[
\int_{0}^{T}\int_{-L}^{L}\left\vert \partial_{t}\varphi(\rho^{n})\right\vert
\leq\left\Vert u^{n}\right\Vert _{L_{t}^{1}L_{x}^{\infty}}\sup_{t\in
\lbrack0,T]}TV\left(  \left[  -L,L\right]  ,\varphi(\rho^{n}(t,\cdot))\right)
+\left\Vert \frac{\mu^{\frac{1}{2}}(\rho^{n})}{\rho^{n}}\right\Vert
_{L^{\infty}(\mathbb{R}^{+}\times\mathbb{R})}\left\Vert \mu^{\frac{1}{2}}%
(\rho^{n})\partial_{x}u^{n}\right\Vert _{L^{1}([0,T]\times\lbrack-L,L])}%
\]
Using $\left(  \text{\ref{cru1}}\right)  $, $\left(  \text{\ref{cru2}}\right)
$, $\left(  \text{\ref{R_rho_bv}}\right)  $ and the energy estimate, we deduce
that $\varphi(\rho_{n})$ belongs uniformly to $BV_{loc}\left(  [0,\infty
)\times\mathbb{R}\right)  $ namely%
\begin{equation}
\left\Vert \varphi(\rho^{n})\right\Vert _{BV(\left[  0,t\right]  \times\left[
-L,L\right]  )}\leq C\left(  t,E_{c}\left(  \rho_{0},u_{0}\right)
,M_{0},\left\Vert \left(  \rho_{0},\frac{1}{\rho_{0}}\right)  \right\Vert
_{L^{\infty}}\right)  . \label{BV_loc_phi}%
\end{equation}
From $\left(  \text{\ref{cru1}}\right)  $, $\left(  \text{\ref{cru2}}\right)
$ we deduce that for any $T,L>0$ we have:%
\begin{equation}
\left\Vert \rho^{n}\right\Vert _{BV(\left[  0,T\right]  \times\left[
-L,L\right]  )}\leq C\left(  t,L,E_{c}\left(  \rho_{0},u_{0}\right)
,M_{0},\left\Vert \left(  \rho_{0},\frac{1}{\rho_{0}}\right)  \right\Vert
_{L^{\infty}}\right)  . \label{BV}%
\end{equation}
with $C(T,L)$ independent on $n$.

\section{Estimates in the Hoff regularity class\label{Hoff_section}}

We recall the Eulerian effective velocity's equation with $\dot{v}%
_{0}=\partial_{t} v_{0}+u\partial_{x} v_{0}$:
\begin{equation}
\rho\dot{v}_{0}-\left(  1-r_{0}\right)  \partial_{x}(\mu\left(  \rho\right)
\partial_{x}v_{0})+\partial_{x}\rho^{\gamma}=0. \label{equation_hoff}%
\end{equation}
In the computations that follow, we drop the $0-$lowerscript of $v_{0}$ and
$r_{0}$ in order to ease readability. As it is by now well-known, once we have
uniform bounds for the density we can obtain some higher order estimates by
using the energy-type functionals first considered by Hoff in \cite{Hof87}
(see also \cite{Hof98}). These are gathered in the following proposition:

\begin{proposition}
\label{Hoff_estimates}The following estimates hold true:%
\[
A\left(  \rho,v_{0}\right)  =\int_{0}^{t}\int_{\mathbb{R}}\sigma\rho\left\vert
\dot{v}_{0}\right\vert ^{2}+\frac{(1-r_{0})}{2}\sigma\left(  t\right)
\int_{\mathbb{R}}\mu\left(  \rho\left(  t\right)  \right)  (\partial_{x}%
v_{0}\left(  t\right)  )^{2}\leq C\left(  t,E_{c}\left(  \rho_{0}%
,u_{0}\right)  ,M_{0},\left\Vert \left(  \rho_{0},\frac{1}{\rho_{0}}\right)
\right\Vert _{L^{\infty}}\right)  ,
\]%
\[
B\left(  \rho,v_{0}\right)  =\frac{1}{2}\int_{\mathbb{R}}\sigma^{2}\left(
t\right)  \rho\left(  t\right)  \dot{v}_{0}^{2}\left(  t\right)  +\left(
1-r_{0}\right)  \int_{0}^{t}\int_{\mathbb{R}}\sigma^{2}\mu\left(  \rho\right)
(\partial_{x}\dot{v}_{0})^{2}\leq C\left(  t,E_{c}\left(  \rho_{0}%
,u_{0}\right)  ,M_{0},\left\Vert \left(  \rho_{0},\frac{1}{\rho_{0}}\right)
\right\Vert _{L^{\infty}}\right)  ,
\]
where $\sigma\left(  t\right)  =\min\left\{  1,t\right\}  $. As consequence we
deduce that%
\[
\int_{0}^{t}\sigma\left(  t\right)  \left\Vert \partial_{x}v_{0}\right\Vert
_{L^{\infty}}^{2}+\sigma\left(  t\right)  ^{\frac{1}{2}}\left\Vert
\partial_{x}v_{0}(t)\right\Vert _{L^{\infty}}\leq C\left(  t,E_{c}\left(
\rho_{0},u_{0}\right)  ,M_{0},\left\Vert \left(  \rho_{0},\frac{1}{\rho_{0}%
}\right)  \right\Vert _{L^{\infty}}\right)  .
\]

\end{proposition}

\begin{remark}
\label{remarca_hoff}In the case of the Navier-Stokes system $\left(
\text{\ref{Navier_Stokes_1d}}\right)  $ the functionals $A\left(
\rho,u\right)  $ and $B\left(  \rho,u\right)  $ are bounded by functions of
the form $C\left(  t,E_{c}\left(  \rho_{0},u_{0}\right)  ,\left\Vert \left(
\rho_{0},\frac{1}{\rho_{0}}\right)  \right\Vert _{L^{\infty}}\right)  $.
\end{remark}

\subsection{The first Hoff energy}

Multiplying the equation $\left(  \text{\ref{equation_hoff}}\right)  $ with
$\dot{v}$ and after some cumbersome yet straightforward computations we get
that%
\begin{align*}
&  \int_{\mathbb{R}}\rho\left\vert \dot{v}\right\vert ^{2}+\frac{d}%
{dt}\left\{  \frac{(1-r)}{2}\int_{\mathbb{R}}\mu\left(  \rho\right)
(\partial_{x}v)^{2}-\int_{\mathbb{R}}\rho^{\gamma}\partial_{x}v\right\} \\
&  =\frac{1-r}{2}\int_{\mathbb{R}}\left(  \mu\left(  \rho\right)  +\rho
\mu^{\prime}\left(  \rho\right)  \right)  \partial_{x}u\left(  \partial
_{x}v\right)  ^{2}+\gamma\int_{\mathbb{R}}\rho^{\gamma}\partial_{x}%
u\partial_{x}v.
\end{align*}
We multiply the previous identity with $\sigma\left(  t\right)  =\min\left\{
1,t\right\}  $ and, by time integration we get that%
\begin{align*}
&  \int_{0}^{t}\int_{\mathbb{R}}\sigma\rho\left\vert \dot{v}\right\vert
^{2}+\frac{(1-r)}{2}\sigma\left(  t\right)  \int_{\mathbb{R}}\mu\left(
\rho\left(  t\right)  \right)  (\partial_{x}v\left(  t\right)  )^{2}\\
&  =\sigma\left(  t\right)  \int_{\mathbb{R}}(\rho\left(  t\right)  ^{\gamma
}-1)\partial_{x}v\left(  t\right)  +\int_{0}^{\min\left\{  1,t\right\}
}\left\{  \frac{(1-r)}{2}\int_{\mathbb{R}}\mu\left(  \rho\right)
(\partial_{x}v)^{2}-\int_{\mathbb{R}}(\rho^{\gamma}-1)\partial_{x}v\right\} \\
&  +\frac{1-r}{2}\int_{0}^{t}\int_{\mathbb{R}}\sigma\left(  \mu\left(
\rho\right)  +\rho\mu^{\prime}\left(  \rho\right)  \right)  \partial
_{x}u\left(  \partial_{x}v\right)  ^{2}+\gamma\int_{0}^{t}\int_{\mathbb{R}%
}\sigma\rho^{\gamma}\partial_{x}u\partial_{x}v.
\end{align*}
Also, we write that%
\begin{align*}
\frac{1}{2}\int_{0}^{t}\int_{\mathbb{R}}\sigma\left(  \mu\left(  \rho\right)
+\rho\mu^{\prime}\left(  \rho\right)  \right)  \partial_{x}u\left(
\partial_{x}v\right)  ^{2}  &  =\frac{1}{2}\int_{0}^{t}\int_{\mathbb{R}}%
\sigma\frac{\left(  \mu\left(  \rho\right)  +\rho\mu^{\prime}\left(
\rho\right)  \right)  }{(\mu\left(  \rho\right)  )^{2}}\partial_{x}u\left(
\mu\left(  \rho\right)  \partial_{x}v-(\rho^{\gamma}-1)\right)  ^{2}\\
&  +\frac{1}{2}\int_{0}^{t}\int_{\mathbb{R}}\sigma\frac{\left(  \mu\left(
\rho\right)  +\rho\mu^{\prime}\left(  \rho\right)  \right)  }{(\mu\left(
\rho\right)  )^{2}}(\rho^{\gamma}-1)^{2}\partial_{x}u\\
&  +\int_{0}^{t}\int_{\mathbb{R}}\sigma\frac{\left(  \mu\left(  \rho\right)
+\rho\mu^{\prime}\left(  \rho\right)  \right)  }{(\mu\left(  \rho\right)
)^{2}}(\rho^{\gamma}-1)\partial_{x}u\partial_{x}v.
\end{align*}
Combining the last two identities we get that%
\begin{align}
&  A\left(  \rho,v\right)  =\int_{0}^{t}\int_{\mathbb{R}}\sigma\rho\left\vert
\dot{v}\right\vert ^{2}+\frac{(1-r)}{2}\sigma\left(  t\right)  \int
_{\mathbb{R}}\mu\left(  \rho\left(  t\right)  \right)  (\partial_{x}v\left(
t\right)  )^{2}\nonumber\\
&  =\sigma\left(  t\right)  \int_{\mathbb{R}}(\rho\left(  t\right)  ^{\gamma
}-1)\partial_{x}v\left(  t\right)  +\int_{0}^{\min\left\{  1,t\right\}
}\left\{  \frac{(1-r)}{2}\int_{\mathbb{R}}\mu\left(  \rho\right)
(\partial_{x}v)^{2}-\int_{\mathbb{R}}(\rho^{\gamma}-1)\partial_{x}v\right\}
\nonumber\\
&  +\int_{0}^{t}\int_{\mathbb{R}}\sigma\gamma\rho^{\gamma}\partial
_{x}u\partial_{x}v+\left(  1-r\right)  \int_{0}^{t}\int_{\mathbb{R}}%
\sigma\frac{\left(  \mu\left(  \rho\right)  +\rho\mu^{\prime}\left(
\rho\right)  \right)  }{(\mu\left(  \rho\right)  )^{2}}\left(  \rho^{\gamma
}-1\right)  \partial_{x}u\partial_{x}v\nonumber\\
&  +\frac{1-r}{2}\int_{0}^{t}\int_{\mathbb{R}}\sigma\frac{\left(  \mu\left(
\rho\right)  +\rho\mu^{\prime}\left(  \rho\right)  \right)  }{(\mu\left(
\rho\right)  )^{2}}(\rho^{\gamma}-1)^{2}\partial_{x}u\nonumber\\
&  +\frac{1-r}{2}\int_{0}^{t}\int_{\mathbb{R}}\sigma\frac{\left(  \mu\left(
\rho\right)  +\rho\mu^{\prime}\left(  \rho\right)  \right)  }{(\mu\left(
\rho\right)  )^{2}}\left(  \mu\left(  \rho\right)  \partial_{x}v-(\rho
^{\gamma}-1)\right)  ^{2}\partial_{x}u. \label{first_hoff_energy_1}%
\end{align}

Owing to the basic energy estimate and using (\ref{cru1}), (\ref{cru2}), the
following inequality holds true:
\begin{gather}
\int_{0}^{\min\left\{  1,t\right\}  }\left\{  \frac{(1-r)}{2}\int_{\mathbb{R}%
}\mu\left(  \rho\right)  (\partial_{x}v)^{2}-\int_{\mathbb{R}}(\rho^{\gamma
}-1)\partial_{x}v\right\} \nonumber\\
+\int_{0}^{t}\int_{\mathbb{R}}\sigma\gamma\rho^{\gamma}\partial_{x}%
u\partial_{x}v+\left(  1-r\right)  \int_{0}^{t}\int_{\mathbb{R}}\sigma
\frac{\left(  \mu\left(  \rho\right)  +\rho\mu^{\prime}\left(  \rho\right)
\right)  }{(\mu\left(  \rho\right)  )^{2}}\left(  \rho^{\gamma}-1\right)
\partial_{x}u\partial_{x}v\nonumber\\
+\frac{1-r}{2}\int_{0}^{t}\int_{\mathbb{R}}\sigma\frac{\left(  \mu\left(
\rho\right)  +\rho\mu^{\prime}\left(  \rho\right)  \right)  }{(\mu\left(
\rho\right)  )^{2}}(\rho^{\gamma}-1)^{2}\partial_{x}u\leq C\left(  t\right)  .
\label{first_hoff_energy_2}%
\end{gather}

Next using again energy estimate and ,(\ref{cru1}), (\ref{cru2}), we have that%
\begin{align}
\sigma\left(  t\right)  \int_{\mathbb{R}}(\rho\left(  t\right)  ^{\gamma
}-1)\partial_{x}v\left(  t\right)   &  \leq\varepsilon\frac{1-r}{2}%
\sigma\left(  t\right)  \int_{\mathbb{R}}\mu\left(  \rho\left(  t\right)
\right)  (\partial_{x}v\left(  t\right)  )^{2}+\frac{1}{2\varepsilon\left(
1-r\right)  }\sigma\left(  t\right)  \int_{\mathbb{R}}\frac{(\rho\left(
t\right)  ^{\gamma}-1)^{2}}{\mu\left(  \rho\left(  t\right)  \right)
}\nonumber\\
&  \leq\varepsilon\frac{1-r}{2}\sigma\left(  t\right)  \int_{\mathbb{R}}%
\mu\left(  \rho\left(  t\right)  \right)  (\partial_{x}v\left(  t\right)
)^{2}+C\left(  t\right)  . \label{first_hoff_energy_3}%
\end{align}

Finally, let us observe that combining (\ref{cru1}), (\ref{cru2}), energy
estimate and Young inequality we get for any $\varepsilon>0$%
\begin{align}
&  \frac{1}{2}\int_{0}^{t}\int_{\mathbb{R}}\sigma\frac{\left(  \mu\left(
\rho\right)  +\rho\mu^{\prime}\left(  \rho\right)  \right)  }{(\mu\left(
\rho\right)  )^{2}}\left(  \mu\left(  \rho\right)  \partial_{x}v-(\rho
^{\gamma}-1)\right)  ^{2}\partial_{x}u\nonumber\\
&  \leq\frac{1}{2}\int_{0}^{t}\int_{\mathbb{R}}\sigma\frac{\left(  \mu\left(
\rho\right)  +\rho\mu^{\prime}\left(  \rho\right)  \right)  }{(\mu\left(
\rho\right)  )^{2}}\left\vert \partial_{x}u\right\vert \left\vert \mu\left(
\rho\right)  \partial_{x}v-(\rho^{\gamma}-1)\right\vert \left\Vert \mu\left(
\rho\right)  \partial_{x}v-(\rho^{\gamma}-1)\right\Vert _{L^{\infty}%
}\nonumber\\
&  \leq\frac{1}{2}\int_{0}^{t}\int_{\mathbb{R}}\sigma\frac{\left(  \mu\left(
\rho\right)  +\rho\mu^{\prime}\left(  \rho\right)  \right)  }{(\mu\left(
\rho\right)  )^{2}}\left\vert \partial_{x}u\right\vert \left\vert \mu\left(
\rho\right)  \partial_{x}v-(\rho^{\gamma}-1)\right\vert \left\Vert \mu\left(
\rho\right)  \partial_{x}v-(\rho^{\gamma}-1)\right\Vert ^{\frac{1}{2}}_{L^{2}%
}\left\Vert \partial_{x}\left(  \mu\left(  \rho\right)  \partial_{x}%
v-\rho^{\gamma}\right)  \right\Vert ^{\frac{1}{2}}_{L^{2}}\nonumber\\
&  \leq\frac{1}{2}C\left(  t\right)  \int_{0}^{t}\sigma\left\Vert \rho\dot
{v}\right\Vert _{L^{2}}^{\frac{1}{2}}\left\Vert \partial_{x}u\right\Vert
_{L^{2}}\left\Vert \mu\left(  \rho\right)  \partial_{x}v-(\rho^{\gamma
}-1)\right\Vert _{L^{2}}^{\frac{3}{2}}\nonumber\\
&  \leq\frac{\varepsilon}{2}\int_{0}^{t}\sigma\left\Vert \rho\dot
{v}\right\Vert _{L^{2}}^{2}+\frac{C\left(  t\right)  }{2\varepsilon}\int
_{0}^{t}\left\Vert \partial_{x}u\right\Vert _{L^{2}}^{\frac{4}{3}}%
\sigma\left\Vert \mu\left(  \rho\right)  \partial_{x}v\right\Vert _{L^{2}}%
^{2}. \label{first_hoff_energy_4}%
\end{align}

Using in the identity $\left(  \text{\ref{first_hoff_energy_1}}\right)  $ the
estimates $\left(  \text{\ref{first_hoff_energy_2}}\right)  $, $\left(
\text{\ref{first_hoff_energy_3}}\right)  $ and $\left(
\text{\ref{first_hoff_energy_4}}\right)  $ along with Gronwall's lemma we get
that%
\begin{equation}
\int_{0}^{t}\int_{\mathbb{R}}\sigma\rho\left\vert \dot{v}_{0}\right\vert
^{2}+\frac{(1-r)}{2}\sigma\left(  t\right)  \int_{\mathbb{R}}\mu\left(
\rho\left(  t\right)  \right)  (\partial_{x}v_{0}\left(  t\right)  )^{2}\leq
C\left(  t\right)  . \label{first_hoff_energy_conclusion}%
\end{equation}
As it is well known from previous works \cite{Hof87}, \cite{BurHas2020} we can
already recover that%
\[
\int_{0}^{t}\left\Vert \partial_{x}v_{0}\right\Vert _{L^{\infty}}\leq\left(
\int_{0}^{t}\sigma^{-\frac{1}{2}}\left(  \tau\right)  d\tau\right)  ^{\frac
{1}{2}}\left(  \int_{0}^{t}\sigma\left(  t\right)  \left\Vert \partial
_{x}v_{0}\right\Vert _{L^{\infty}}^{2}\right)  ^{\frac{1}{2}}.
\]
Furthermore for the NSK system we have similar estimate for $v_{1}$ but these
estimates are not uniform in $c>0$.

\subsection{The second Hoff energy}

In the following lines, we aim at obtaining higher order estimates. In order
to do so, we apply the operator $\partial_{t}+u\partial_{x}$ to the equation
of $v_{0}$ (again we drop the $0-$lower script) in order to obtain that%
\[
\left(  \partial_{t}+u\partial_{x}\right)  \left(  \rho\dot{v}\right)
-\left(  1-r\right)  \left(  \partial_{t}+u\partial_{x}\right)  \partial
_{x}\left(  \mu\left(  \rho\right)  \partial_{x}v\right)  +\left(
\partial_{t}+u\partial_{x}\right)  \partial_{x}\rho^{\gamma}=0.
\]
Multiplying with $\dot{v}$ we see that%
\[
\int_{\mathbb{R}}\partial_{t}(\rho\dot{v})\dot{v}=\int_{\mathbb{R}}%
\partial_{t}\rho\dot{v}^{2}+\int_{\mathbb{R}}\rho\frac{d}{dt}\frac{\dot{v}%
^{2}}{2}=\frac{1}{2}\frac{d}{dt}\left\{  \int_{\mathbb{R}}\rho\dot{v}%
^{2}\right\}  +\frac{1}{2}\int_{\mathbb{R}}\partial_{t}\rho\dot{v}^{2}%
\]%
\[
\int_{\mathbb{R}}\left(  \partial_{t}+u\partial_{x}\right)  \left(  \rho
\dot{v}\right)  \dot{v}=\frac{1}{2}\frac{d}{dt}\left\{  \int_{\mathbb{R}}%
\rho\dot{v}^{2}\right\}  -\int_{\mathbb{R}}\rho\dot{v}^{2}\partial_{x}u.
\]
Next, we have that%
\begin{align*}
&  -\int_{\mathbb{R}}\left(  \partial_{t}+u\partial_{x}\right)  \partial
_{x}\left(  \mu\left(  \rho\right)  \partial_{x}v\right)  \dot{v}\\
&  =\int_{\mathbb{R}}\partial_{t}\mu\left(  \rho\right)  \partial_{x}%
v\partial_{x}\dot{v}+\int_{\mathbb{R}}\mu\left(  \rho\right)  \partial
_{x}\partial_{t}v\partial_{x}\dot{v}-\int_{\mathbb{R}}\mu\left(  \rho\right)
\partial_{x}v\partial_{xx}^{2}\left(  u\dot{v}\right)
\end{align*}
Let us deal with the last term%
\begin{align*}
&  -\int_{\mathbb{R}}\mu\left(  \rho\right)  \partial_{x}v\partial_{xx}%
(u\dot{v})\\
&  =-\int_{\mathbb{R}}\mu\left(  \rho\right)  \partial_{x}v(u\partial_{xx}%
^{2}\dot{v}+2\partial_{x}u\partial_{x}\dot{v}+\dot{v}\partial_{xx}^{2}u)\\
&  =\int_{\mathbb{R}}\partial_{x}(\mu\left(  \rho\right)  \partial
_{x}vu)\partial_{x}\dot{v}-2\int_{\mathbb{R}}\mu\left(  \rho\right)
\partial_{x}v\partial_{x}u\partial_{x}\dot{v}-\int_{\mathbb{R}}\mu\left(
\rho\right)  \dot{v}\partial_{x}v\partial_{xx}^{2}u\\
&  =\int_{\mathbb{R}}u\partial_{x}(\mu\left(  \rho\right)  )\partial
_{x}v\partial_{x}\dot{v}+\int_{\mathbb{R}}\mu(\rho)\partial_{x}(u\partial
_{x}v)\partial_{x}\dot{v}-2\int_{\mathbb{R}}\mu\left(  \rho\right)
\partial_{x}v\partial_{x}u\partial_{x}\dot{v}-\int_{\mathbb{R}}\mu\left(
\rho\right)  \dot{v}\partial_{x}v\partial_{xx}^{2}u
\end{align*}
Putting togeather the last two relations, we get that%
\begin{align*}
&  -\int_{\mathbb{R}}\left(  \partial_{t}+u\partial_{x}\right)  \partial
_{x}\left(  \mu\left(  \rho\right)  \partial_{x}v\right)  \dot{v}\\
&  =\int_{\mathbb{R}}\partial_{t}\mu\left(  \rho\right)  \partial_{x}%
v\partial_{x}\dot{v}+\int_{\mathbb{R}}\mu\left(  \rho\right)  \partial
_{x}\partial_{t}v\partial_{x}\dot{v}\\
&  +\int_{\mathbb{R}}u\partial_{x}(\mu\left(  \rho\right)  )\partial
_{x}v\partial_{x}\dot{v}+\int_{\mathbb{R}}\mu(\rho)\partial_{x}(u\partial
_{x}v)\partial_{x}\dot{v}-2\int_{\mathbb{R}}\mu\left(  \rho\right)
\partial_{x}v\partial_{x}u\partial_{x}\dot{v}-\int_{\mathbb{R}}\mu\left(
\rho\right)  \dot{v}\partial_{x}v\partial_{xx}^{2}u\\
&  =\left(  \int_{\mathbb{R}}\partial_{t}\mu\left(  \rho\right)  \partial
_{x}v\partial_{x}\dot{v}+\int_{\mathbb{R}}u\partial_{x}(\mu\left(
\rho\right)  )\partial_{x}v\partial_{x}\dot{v}\right)  +\left(  \int
_{\mathbb{R}}\mu\left(  \rho\right)  \partial_{x}\partial_{t}v\partial_{x}%
\dot{v}+\int_{\mathbb{R}}\mu(\rho)\partial_{x}(u\partial_{x}v)\partial_{x}%
\dot{v}\right) \\
&  -2\int_{\mathbb{R}}\mu\left(  \rho\right)  \partial_{x}v\partial
_{x}u\partial_{x}\dot{v}-\int_{\mathbb{R}}\mu\left(  \rho\right)  \dot
{v}\partial_{x}v\partial_{xx}^{2}u\\
&  =-\int_{\mathbb{R}}\rho\mu^{\prime}\left(  \rho\right)  \partial
_{x}u\partial_{x}v\partial_{x}\dot{v}+\int_{\mathbb{R}}\mu\left(  \rho\right)
(\partial_{x}\dot{v})^{2}-2\int_{\mathbb{R}}\mu\left(  \rho\right)
\partial_{x}v\partial_{x}u\partial_{x}\dot{v}-\int_{\mathbb{R}}\mu\left(
\rho\right)  \dot{v}\partial_{x}v\partial_{xx}^{2}u\\
&  =-\int_{\mathbb{R}}\rho\mu^{\prime}\left(  \rho\right)  \partial
_{x}u\partial_{x}v\partial_{x}\dot{v}+\int_{\mathbb{R}}\mu\left(  \rho\right)
(\partial_{x}\dot{v})^{2}-\int_{\mathbb{R}}\mu\left(  \rho\right)
\partial_{x}v\partial_{x}u\partial_{x}\dot{v}-\int_{\mathbb{R}}\mu\left(
\rho\right)  \partial_{x}\left(  \dot{v}\partial_{x}u\right)  \partial_{x}v\\
&  =-\int_{\mathbb{R}}\rho\mu^{\prime}\left(  \rho\right)  \partial
_{x}u\partial_{x}v\partial_{x}\dot{v}+\int_{\mathbb{R}}\mu\left(  \rho\right)
(\partial_{x}\dot{v})^{2}-\int_{\mathbb{R}}\mu\left(  \rho\right)
\partial_{x}v\partial_{x}u\partial_{x}\dot{v}+\int_{\mathbb{R}}\dot{v}%
\partial_{x}u\partial_{x}(\mu\left(  \rho\right)  \partial_{x}v)
\end{align*}

Next, let us take a look at%
\begin{align*}
\int_{\mathbb{R}}\left(  \partial_{x}\partial_{t}\rho^{\gamma}+u\partial
_{xx}^{2}\rho^{\gamma}\right)  \dot{v}  &  =-\int_{\mathbb{R}}\partial_{t}%
\rho^{\gamma}\partial_{x}\dot{v}+\int_{\mathbb{R}}u\partial_{xx}\rho^{\gamma
}\dot{v}\\
&  =\int_{\mathbb{R}}u\partial_{x}\rho^{\gamma}\partial_{x}\dot{v}+\gamma
\int_{\mathbb{R}}\rho^{\gamma}\partial_{x}u\partial_{x}\dot{v}+\int
_{\mathbb{R}}u\partial_{xx}\rho^{\gamma}\dot{v}\\
&  =\gamma\int_{\mathbb{R}}\rho^{\gamma}\partial_{x}u\partial_{x}\dot{v}%
+\int_{\mathbb{R}}u\partial_{x}(\dot{v}\partial_{x}\rho^{\gamma})\\
&  =\gamma\int_{\mathbb{R}}\rho^{\gamma}\partial_{x}u\partial_{x}\dot{v}%
-\int_{\mathbb{R}}\dot{v}(\partial_{x}\rho^{\gamma})(\partial_{x}u)
\end{align*}
We gather the last three relations in order to obtain%
\begin{align*}
0  &  =\frac{1}{2}\frac{d}{dt}\left\{  \int_{\mathbb{R}}\rho\dot{v}%
^{2}\right\}  -\int_{\mathbb{R}}\rho\dot{v}^{2}\partial_{x}u\\
&  -\left(  1-r\right)  \int_{\mathbb{R}}\rho\mu^{\prime}\left(  \rho\right)
\partial_{x}u\partial_{x}v\partial_{x}\dot{v}+\left(  1-r\right)
\int_{\mathbb{R}}\mu\left(  \rho\right)  (\partial_{x}\dot{v})^{2}\\
&  -\left(  1-r\right)  \int_{\mathbb{R}}\mu\left(  \rho\right)  \partial
_{x}v\partial_{x}u\partial_{x}\dot{v}+\left(  1-r\right)  \int_{\mathbb{R}%
}\dot{v}\partial_{x}u\partial_{x}(\mu\left(  \rho\right)  \partial_{x}v)\\
&  +\gamma\int_{\mathbb{R}}\rho^{\gamma}\partial_{x}u\partial_{x}\dot{v}%
-\int_{\mathbb{R}}\dot{v}(\partial_{x}\rho^{\gamma})(\partial_{x}u)\\
&  =\frac{1}{2}\frac{d}{dt}\left\{  \int_{\mathbb{R}}\rho\dot{v}^{2}\right\}
+\left(  1-r\right)  \int_{\mathbb{R}}\mu\left(  \rho\right)  (\partial
_{x}\dot{v})^{2}\\
&  -\left(  1-r\right)  \int_{\mathbb{R}}(\rho\mu^{\prime}\left(  \rho\right)
+\mu\left(  \rho\right)  )\partial_{x}u\partial_{x}v\partial_{x}\dot{v}%
-\int_{\mathbb{R}}\rho\dot{v}^{2}\partial_{x}u+\gamma\int_{\mathbb{R}}%
\rho^{\gamma}\partial_{x}u\partial_{x}\dot{v}\\
&  +\int_{\mathbb{R}}\dot{v}\partial_{x}u\partial_{x}(\left(  1-r\right)
\mu\left(  \rho\right)  \partial_{x}v-\rho^{\gamma}).
\end{align*}
We end up with%
\begin{align*}
&  \frac{1}{2}\frac{d}{dt}\left\{  \int_{\mathbb{R}}\rho\dot{v}^{2}\right\}
+\left(  1-r\right)  \int_{\mathbb{R}}\mu\left(  \rho\right)  (\partial
_{x}\dot{v})^{2}\\
&  =\left(  1-r\right)  \int_{\mathbb{R}}(\rho\mu^{\prime}\left(  \rho\right)
+\mu\left(  \rho\right)  )\partial_{x}u\partial_{x}v\partial_{x}\dot{v}%
+\gamma\int_{\mathbb{R}}\rho^{\gamma}\partial_{x}u\partial_{x}\dot{v}.
\end{align*}
We multiply the last relation with $\sigma^{2}\left(  t\right)  $ and we
integrate in time in order to obtain%
\begin{align}
&  B\left(  \rho,v\right)  =\frac{1}{2}\int_{\mathbb{R}}\sigma^{2}\left(
t\right)  \rho\left(  t\right)  \dot{v}^{2}\left(  t\right)  +\left(
1-r\right)  \int_{0}^{t}\int_{\mathbb{R}}\sigma^{2}\mu\left(  \rho\right)
(\partial_{x}\dot{v})^{2}\nonumber\\
&  =\frac{1}{2}\int_{0}^{1}\int_{\mathbb{R}}\sigma\rho\dot{v}^{2}+\left(
1-r\right)  \int_{0}^{t}\int_{\mathbb{R}}\sigma^{2}(\rho\mu^{\prime}\left(
\rho\right)  +\mu\left(  \rho\right)  )\partial_{x}u\partial_{x}v\partial
_{x}\dot{v}+\gamma\int_{0}^{t}\int_{\mathbb{R}}\sigma^{2}\rho^{\gamma}%
\partial_{x}u\partial_{x}\dot{v} \label{second_hoff_energy_1}%
\end{align}
So, we see that owing to $\left(  \text{\ref{first_hoff_energy_conclusion}%
}\right)  $ we have that
\begin{equation}
\frac{1}{2}\int_{0}^{1}\int_{\mathbb{R}}\sigma\rho\dot{v}^{2}\leq C\left(
t\right)  . \label{second_hoff_energy_2}%
\end{equation}
The second term in $\left(  \text{\ref{second_hoff_energy_1}}\right)  $, we
see that%
\begin{align}
&  \left(  1-r\right)  \int_{0}^{t}\int_{\mathbb{R}}\sigma^{2}(\rho\mu
^{\prime}\left(  \rho\right)  +\mu\left(  \rho\right)  )\partial_{x}%
u\partial_{x}v\partial_{x}\dot{v}\nonumber\\
&  =\int_{0}^{t}\int_{\mathbb{R}}\sigma^{2}\frac{(\rho\mu^{\prime}\left(
\rho\right)  +\mu\left(  \rho\right)  )}{\mu\left(  \rho\right)  }\left(
1-r\right)  \left(  \mu\left(  \rho\right)  \partial_{x}v-(\rho^{\gamma
}-1)\right)  \partial_{x}u\partial_{x}\dot{v}\nonumber\\
&  +\int_{0}^{t}\int_{\mathbb{R}}\sigma^{2}\frac{(\rho\mu^{\prime}\left(
\rho\right)  +\mu\left(  \rho\right)  )}{\mu\left(  \rho\right)  }%
(\rho^{\gamma}-1)\partial_{x}u\partial_{x}\dot{v}.
\label{second_hoff_energy_3}%
\end{align}
We observe that using (\ref{cru1}), (\ref{cru2}) and Young inequality%
\begin{align}
&  \int_{0}^{t}\int_{\mathbb{R}}\sigma^{2}\frac{(\rho\mu^{\prime}\left(
\rho\right)  +\mu\left(  \rho\right)  )}{\mu\left(  \rho\right)  }%
(\rho^{\gamma}-1)\partial_{x}u\partial_{x}\dot{v}+\gamma\int_{0}^{t}%
\int_{\mathbb{R}}\sigma^{2}\rho^{\gamma}\partial_{x}u\partial_{x}\dot
{v}\nonumber\\
&  \leq\varepsilon\int_{0}^{t}\int_{\mathbb{R}}\sigma^{2}\mu\left(
\rho\right)  (\partial_{x}\dot{v})^{2}+C\left(  t\right)  .
\label{second_hoff_energy_4}%
\end{align}
Also, we have for any $\varepsilon>0$ by using again (\ref{cru1}),
(\ref{cru2}, Young inequality) and Gagliardo-Nirenberg inequality%
\begin{align}
&  \int_{0}^{t}\int_{\mathbb{R}}\sigma^{2}\frac{(\rho\mu^{\prime}\left(
\rho\right)  +\mu\left(  \rho\right)  )}{\mu\left(  \rho\right)  }\left(
1-r\right)  \left(  \mu\left(  \rho\right)  \partial_{x}v-(\rho^{\gamma
}-1)\right)  \partial_{x}u\partial_{x}\dot{v}\nonumber\\
&  \leq\int_{0}^{t}\int_{\mathbb{R}}\sigma^{2}\left\Vert \mu\left(
\rho\right)  \partial_{x}v-(\rho^{\gamma}-1)\right\Vert _{L^{2}}^{\frac{1}{2}%
}\left\Vert \rho\dot{v}\right\Vert _{L^{2}}^{\frac{1}{2}}\left\Vert
\partial_{x}u\right\Vert _{L^{2}}\left\Vert \partial_{x}\dot{v}\right\Vert
_{L^{2}}\nonumber\\
&  \leq C\left(  t\right)  \int_{0}^{t}\sigma^{2}\left\Vert \partial
_{x}u\right\Vert _{L^{2}}^{2}\left\Vert \mu\left(  \rho\right)  \partial
_{x}v-(\rho^{\gamma}-1)\right\Vert _{L^{2}}\left\Vert \rho^{\frac{1}{2}}%
\dot{v}\right\Vert _{L^{2}}+\varepsilon\int_{0}^{t}\int_{\mathbb{R}}\sigma
^{2}\mu\left(  \rho\right)  (\partial_{x}\dot{v})^{2}\nonumber\\
&  \leq C\left(  t\right)  \int_{0}^{t}\sigma^{2}\left\Vert \partial
_{x}u\right\Vert _{L^{2}}^{2}\left\Vert \rho^{\frac{1}{2}}\dot{v}\right\Vert
_{L^{2}}^{2}+C\left(  t\right)  \int_{0}^{t}\sigma^{2}\left\Vert \partial
_{x}u\right\Vert _{L^{2}}^{2}\left\Vert \mu\left(  \rho\right)  \partial
_{x}v-(\rho^{\gamma}-1)\right\Vert _{L^{2}}^{2}+\varepsilon\int_{0}^{t}%
\int_{\mathbb{R}}\sigma^{2}\mu\left(  \rho\right)  (\partial_{x}\dot{v})^{2}
\label{second_hoff_energy_5}%
\end{align}
Putting together $\left(  \text{\ref{first_hoff_energy_conclusion}}\right)  ,$
$\left(  \text{\ref{second_hoff_energy_1}}\right)  $, $\left(
\text{\ref{second_hoff_energy_2}}\right)  $, $\left(
\text{\ref{second_hoff_energy_3}}\right)  ,$ $\left(
\text{\ref{second_hoff_energy_4}}\right)  , $ and $\left(
\text{\ref{second_hoff_energy_5}}\right)  $ we get that%
\[
\frac{1}{2}\int_{\mathbb{R}}\sigma^{2}\left(  t\right)  \rho\left(  t\right)
\dot{v}^{2}\left(  t\right)  +\left(  1-r\right)  \int_{0}^{t}\int
_{\mathbb{R}}\sigma^{2}\mu\left(  \rho\right)  (\partial_{x}\dot{v})^{2}\leq
C\left(  t\right)  .
\]

The control over $\left\Vert \frac{1}{\rho}\right\Vert _{L^{\infty}}$,
$A\left(  \rho,u\right)  $ and $B\left(  \rho,u\right)  $ gives us
\begin{equation}
\sigma\left(  t\right)  ^{\frac{1}{2}}\left\Vert \partial_{x}v_{0}%
(t)\right\Vert _{L^{\infty}}\leq C\left(  t\right)  , \label{u_x_in_linfty}%
\end{equation}
for any $t\geq0$. The same estimate are again true for $v_{1}$ for the NSK
system but are not uniform in $c>0$. This concludes the proof of Proposition
\ref{Hoff_estimates}

\subsection{End of proof of the main results}

\label{sec4}

The sequence of solutions constructed in Section \ref{subsec1} $\left(
\rho^{n},u^{n}\right)  _{n\in\mathbb{N}}$ are globally defined owing to the
intermediary results proved in Proposition \ref{rho_is_bounded_in_Linfty},
Proposition \ref{rho_is_bounded_by_bellow}, the estimates from Proposition
\ref{Hoff_estimates} and the remark that follows. Using this, we may extract a
subsequence that converges to a solution of the NS system $\left(
\text{\ref{Navier_Stokes_1d}}\right)  $. The fact that the solution itself
verifies the estimates announced in Theorem \ref{theo3} is a consequence of
the Fatou lemma. We skip the technical details and we refer the reader to
\cite{MV,Jiu} for more details on this subject. It proves in particular the
Theorem \ref{theo3}.

\bigskip

Regarding Theorem \ref{theo1}, we skip the details showing that the
approximate solution $(\rho^{n},u^{n})_{n\in\mathbb{N}}$ converge up to a
subsequence to a unique solution $(\rho,u)$ of the system $\left(
\text{\ref{NSK_intro}}\right)  $ on $[0,T]$. Again, the fact that the solution
itself verifies the estimates announced in Theorem \ref{theo1} is a
consequence of the Fatou lemma. One can consult in particular
\cite{BurHas2020}.

In what Theorem \ref{theo2} is concerned, we only mention that the only not
so-classical estimates $\left(  \text{\ref{BV_loc_phi}}\right)  $ and $\left(
\text{\ref{BV}}\right)  $ are used in order to justify the fact that the
condition%
\[
u^{c}\left(  t,x\right)  +r_{1}\left(  c\right)  \partial_{x}\varphi\left(
\rho^{c}\left(  t,x\right)  \right)  \leq C\left(  t,E_{c}\left(  \rho
_{0},u_{0}\right)  ,M_{0},\left\Vert \left(  \rho_{0},\frac{1}{\rho_{0}%
}\right)  \right\Vert _{L^{\infty}}\right)  ,
\]
is obtained also in the limit $c\rightarrow0$, namely%
\[
u+\partial_{x}\varphi\left(  \rho\right)  \leq C\left(  t,E_{c}\left(
\rho_{0},u_{0}\right)  ,M_{0},\left\Vert \left(  \rho_{0},\frac{1}{\rho_{0}%
}\right)  \right\Vert _{L^{\infty}}\right)  ,
\]
in the sense of measures. The arguments concerning the weak convergence up to
subsequence of the sequence $(\rho_{c},u_{c})_{c>0}$ are the same as those in
the proof of the Theorem \ref{theo2}. It is important to point out that in
fact $(\rho_{c},u_{c})_{c>0}$ converges weakly and not only up to a
subsequence because $(\rho_{c},u_{c})_{c>0}$ has a unique accumulation point,
indeed the limit $(\rho,u)$ when $c$ goes to $0$ is solution of the NS system
and this solution is unique as we will see.

From the Hoff type estimates, we can deduce that $\partial_{x}u$ belongs to
$L^{2}([0,T],L^{\infty}(\mathbb{R}))$ for any $T>0$. The uniqueness of the
solution $(\rho,u)$ is then a consequence of the fact that on $[0,T]$,
$\partial_{x}u$ belongs to $L^{1}([0,T_{1}],L^{\infty})$. We refer to
\cite{BurHas2020} for the details of the proof in the case of the compressible
Navier-Stokes system where we use a Lagrangian formulation. The adaptation of
the proof to the Navier-Stokes Korteweg system is direct.

\section*{Acknowledgements}

CB has been partially funded by the projects SingFlows (ANR-18-CE40-0027-01)
and CRISIS (ANR-20-CE40-0020-01) both operated by the French National Research
Agency (ANR). BH has been partially funded by the ANR project INFAMIE ANR-15-CE40-0011.

\bibliographystyle{alpha}
\bibliography{BibliografieGeneralaCompleta}

\end{document}